%% file: kdv_mur_gap.tex
\documentclass[11pt]{article}
\usepackage{amssymb,amsmath,amsthm,mathtools}
\usepackage{graphicx,booktabs,array,longtable,float,alltt}
\usepackage[dvipsnames]{xcolor}
\usepackage[hidelinks]{hyperref}
\usepackage{geometry}
\allowdisplaybreaks
\numberwithin{equation}{section}
\newtheorem{theorem}{Theorem}[section]
\newtheorem{lemma}[theorem]{Lemma}
\newtheorem{proposition}[theorem]{Proposition}
\newtheorem{corollary}[theorem]{Corollary}

\theoremstyle{remark}
\newtheorem{remark}[theorem]{Remark}
\theoremstyle{plain}

\DeclareMathOperator{\Rea}{Re}
\DeclareMathOperator{\Ima}{Im}
\DeclareMathOperator{\Log}{Log}
\newcommand{\norm}[1]{\left\lVert #1\right\rVert}
\newcommand{\abs}[1]{\left\lvert #1\right\rvert}
\newcommand{\C}{\mathbb C}
\newcommand{\R}{\mathbb R}
\newcommand{\D}{\mathfrak D}
\newcommand{\Cd}{\mathfrak C}
\newcommand{\lean}[1]{{\small\nolinkurl{#1}}}

\begin{document}
\title{A Large Gap between Galerkin Errors with and without Mesh-Uniform
Regularity for the KdV Equation with Boundary Dissipation\thanks{This work
is supported by the National Natural Science Foundation of China (grant
No. 12171431).\medskip}}
\author{Kangsheng Liu\thanks{Academy for Advanced Interdisciplinary Studies,
Northeast Normal University, Changchun 130024, Jilin, China; School of
Mathematical Sciences, Zhejiang University, Hangzhou, 310058, China
(email: ksliu@zju.edu.cn).}}
\date{}
\maketitle

\begin{center}
\begin{minipage}{6in}
{\bf Abstract.}
We compare two cubic Hermite Galerkin semidiscretizations of the
Korteweg--de Vries equation on a bounded interval with boundary
dissipation. They use the same trial space. The first is the weighted
scheme considered in our previous paper~\cite{Main}, called the W-scheme. It has mesh-uniform
regularity: the discrete semigroups of the linearized equation inherit
Crandall--Pazy regularity and Kato smoothing, with bounds uniform in the
mesh size $h$. The second is the standard Galerkin method in $L^2$,
called the S-scheme. It has no mesh-uniform regularity. For the
W-scheme, the full error and the ordering error converge at an algebraic
rate in $h$, uniformly on $[t_0,\infty)$ times small balls of $L^2$ for any
fixed $t_0>0$. They are bounded by
$Ch^p\varepsilon$ and $Ch^p\varepsilon^2$, respectively, where $\varepsilon$ is the
size of the datum and $p>0$. We prove that neither bound holds for the
S-scheme by using exact high-frequency spurious eigenpairs of the S-scheme. The
existence of such eigenpairs for meshes with at least $512$ elements and
the failure of both bounds have proofs checked by Lean~4. The only premise
is the mesh-uniform Crandall--Pazy estimate of the W-scheme
(arXiv:2610.39269).
Computations compare the two
schemes on unit-norm data built from these eigenfunctions, on $512$,
$1024$ and $2048$ elements. The full error of the
S-scheme stays at $4.6\times10^{-2}$. That of the W-scheme is at most
$2.2\times10^{-9}$. The ordering error of the S-scheme stays at
$9.6\times10^{-6}$. That of the W-scheme falls from $2.3\times10^{-9}$ to
$1.4\times10^{-10}$. Their ratio grows from $4.2\times10^3$ to
$6.6\times10^4$.
\vskip 3mm

\vskip 4mm

\vskip 3mm
\end{minipage}
\end{center}

\section{Introduction}\label{sec:intro}

How much does mesh-uniform regularity improve the error estimates of a
Galerkin method? We answer this question for the Korteweg--de Vries (KdV)
equation on an interval with boundary dissipation,
\begin{equation}\label{eq:pde}
 u_t+u_{xxx}+uu_x=0,\qquad u(0,t)=u(L,t)=u_x(L,t)=0,\qquad u(x,0)=u_0(x).
\end{equation}
In~\cite{Main} we studied a weighted cubic Hermite Galerkin scheme for
\eqref{eq:pde}, called the W-scheme below. For this scheme, the discrete
semigroups of the linearized equation inherit exponential stability,
Crandall--Pazy regularity of index $3/2$ and Kato smoothing, with bounds
uniform in the mesh size. For the continuous semigroup the index $3/2$ is
exact~\cite[Theorem~2.4]{Main}. A regularity property inherited with
bounds uniform in the mesh size is called \emph{mesh-uniform regularity}
(Section~\ref{sec:mur}). Crandall--Pazy regularity and Kato smoothing are
two kinds of it. From them we derived error estimates in which the error
and the datum are measured in the same norm, uniformly on balls. They hold
uniformly on $[t_0,\infty)$ for any fixed $t_0>0$, and for the ordering
error for all time. Here we compare the W-scheme with the standard
Galerkin method in $L^2(0,L)$, called the S-scheme, on the same trial
space. Both schemes are contractive. The S-scheme has no mesh-uniform
regularity. We show that the gap between the two error estimates is
large. We also quantify it by computation. The weight is not what matters. It
is one way to obtain mesh-uniform regularity. The gap comes from this
regularity.

Table~\ref{tab:intro-numbers} shows the gap on one datum. The datum $f_h$
has unit $L^2$ norm and depends on the mesh. It is described below. On
the meshes with $512$, $1024$ and $2048$ elements, the computed full error
$u-u_h$ of the W-scheme is at most $2.2\times10^{-9}$. The computed full
error $u-\mathring u_h$ of the S-scheme stays at $4.6\times10^{-2}$. The
ordering error, introduced in~\cite{Main}, measures how the nonlinearity
changes the discretization error. It is $\Cd_h=(u-u_h)-(v-v_h)$ for the
W-scheme and $\mathring\Cd_h=(u-\mathring u_h)-(v-\mathring v_h)$ for the
S-scheme. Here $v$, $v_h$ and $\mathring v_h$ are the linear solutions.
For the W-scheme the computed ordering error is divided by about $4$ per
halving of $h$. For the S-scheme it stays at $9.6\times10^{-6}$. The ratio
grows from $4.2\times10^3$ to $6.6\times10^4$.

\begin{table}[H]
\centering\small
\caption{Errors on the datum $\varepsilon f_h$.}\label{tab:intro-numbers}
\begin{tabular}{@{}lccc@{}}
\toprule
number of elements & $512$ & $1024$ & $2048$\\
\midrule
full error, W-scheme & $2.2\times10^{-9}$ & $2.0\times10^{-9}$ & $3.2\times10^{-10}$\\
full error, S-scheme & $4.594\times10^{-2}$ & $4.596\times10^{-2}$ & $4.597\times10^{-2}$\\
ordering error, W-scheme & $2.31\times10^{-9}$ & $5.78\times10^{-10}$ & $1.44\times10^{-10}$\\
ordering error, S-scheme & $9.590\times10^{-6}$ & $9.592\times10^{-6}$ & $9.593\times10^{-6}$\\
ratio of ordering errors, S-scheme/W-scheme &$4.2\times10^{3}$ & $1.7\times10^{4}$ & $6.6\times10^{4}$\\
\bottomrule
\end{tabular}
\end{table}

In Table~\ref{tab:intro-numbers}, $L=1$, $\varepsilon=1/8$ and
$t=1/\mathcal C_S\approx6.3\times10^{-4}$. Here $\mathcal C_S=1575.38\ldots$ is
the limit of $-\Rea\zeta_h$ for $L=1$, where $\zeta_h$ are the eigenvalues of
the high-frequency spurious eigenpairs described below. The errors are $L^2$ norms. The
values are numerical estimates obtained by the procedures of
Appendix~\ref{app:num}. The
W-scheme values are uncertain by up to about $6\%$
(Section~\ref{app:fh}).

Table~\ref{tab:intro-theory} compares the theorems. The left column is
from~\cite{Main}. The right column is proved here. The rows on errors come
first. In the table, $L=1$ and $t>0$ is fixed. We write
$t_\wedge=\min\{1,t\}$ and $a=\min\{1,ht_\wedge^{-3/2}\}$. The decay rate
satisfies $0<\nu<3\pi^2/4$. The norms are those of the weighted space $H$ for the
W-scheme and of $L^2$ for the S-scheme. They are equivalent on
$L^2(0,1)$.

\begin{table}[H]
\centering\small
\caption{Theorems with and without mesh-uniform regularity (MUR).}\label{tab:intro-theory}
\setlength{\tabcolsep}{4pt}
\begin{tabular}{@{}>{\raggedright\arraybackslash}p{0.27\textwidth}
>{\raggedright\arraybackslash}p{0.33\textwidth}
>{\raggedright\arraybackslash}p{0.34\textwidth}@{}}
\toprule
 & with MUR (W-scheme) & without MUR (S-scheme)\\
\midrule
linear error, worst unit datum & $\le Ce^{-\nu t}\,h\,t_\wedge^{-3/2}$
 \cite[Thm~4.2]{Main} & $\ge e^{-\mathcal C_St}-o(1)$ (Cors~\ref{cor:principleS}, \ref{cor:lineargap}, \ref{cor:fixed})\\
linear error, fixed datum & rate $h$ for every datum & converges; every prescribed rate fails for some datum (Cor~\ref{cor:fixed})\\
full error, ball of radius $\varepsilon$ & $\le Ce^{-\nu t}\{a\varepsilon+h^{1/6}\varepsilon^2\}$ \cite[Thm~4.2]{Main} & relative error $\to e^{-\mathcal C_St}$ for data $\varepsilon_hf_h$, $\varepsilon_h=o(h)$ (Thm~\ref{thm:total})\\
ordering error, ball of radius $\varepsilon$ & $\le Ch^{1/6}\varepsilon^2$ for all $t\ge0$ \cite[Thm~4.3]{Main} & on $\varepsilon f_h$: $\varepsilon^{-2}\mathring\Cd_h(t)\to-E(t)\ne0$ as first $\varepsilon\downarrow0$, then $h\to0$; hence no bound $Ch^p\varepsilon^2$, $p>0$, uniform in $h$ and $\varepsilon$ (Cor~\ref{cor:noproduct})\\
\midrule
Crandall--Pazy constant $\sup_{t\le1}t^{3/2}\norm{A_hT_h(t)}$ & bounded \cite[Thm~3.5]{Main} & unbounded: $\|\mathring A_h\mathring T_h(t)\|\ge100h^{-3}e^{-1600t}$ (Cor~\ref{cor:failures})\\
Kato smoothing & uniform \cite[Section~3.2]{Main} & fails (Cor~\ref{cor:failures})\\
high-frequency spurious modes & none (Section~\ref{sec:mur}) & exist for $N\ge512$ (Thm~\ref{thm:top})\\
\bottomrule
\end{tabular}
\end{table}

\paragraph{The datum and the high-frequency spurious eigenpairs.}
The data $f_h$ are real mesh-dependent functions in $X_h$ with
$\norm{f_h}_2=1$, built from the high-frequency spurious eigenpairs of the S-scheme. For every uniform
mesh with $N\ge512$ elements, the discrete operator of the S-scheme has an exact
eigenpair whose eigenvalue $\zeta_h$ satisfies $|\Ima\zeta_h|\ge100h^{-3}$
and $0<-\Rea\zeta_h\le1600/L^3$ (Theorem~\ref{thm:top}(a)). It sits near a
frequency of zero group velocity of the grid waves. Such a mode oscillates
on the mesh scale and is damped only at a rate bounded independently of
the mesh. For the continuous problem and for the W-scheme, the damping of a
mode with eigenvalue $\zeta$ grows at least like $|\Ima\zeta|^{2/3}$
(Section~\ref{sec:mur}). We put $f_h=\sqrt2\Rea\varphi_h$, normalized,
where $\varphi_h$ is the eigenfunction. Since $f_h\in X_h$, both schemes
start from $f_h$ itself, without projection error.

\paragraph{Lean proofs.}
In substance, Theorem~\ref{thm:top} concerns a classical problem: the
asymptotic location of a sequence of eigenvalues. Its statement is easy
to read. Its proof is long and full of explicit constants.
An asymptotic expansion only indicates where a zero of the boundary
determinant should be. Existence needs a contraction argument on an
explicit disc, with explicit constants. For the threshold $N\ge512$ this
argument involves rational certificates, analytic root branches and
explicit Taylor remainders. The proof of Theorem~\ref{thm:top} is a
Lean~4 proof. Remark~\ref{rem:history} explains how it was obtained, by
successive iterations of the formalization, and what the author checked.

Given Theorem~\ref{thm:top} and the results derived from its proof, the
theoretical results of this paper can be proved and checked by hand.
Section~\ref{sec:spurious} states the linear approximation results with
sketches that suffice for a proof by hand. For the properties of the
eigenfunctions, it gives the idea and the key quantities.
Section~\ref{sec:gap} proves the results on the nonlinear problem by
hand. However, the author wanted to know whether the iterations of the
Lean formalization could reduce the premises to the essential one. In the
end, one Lean theorem collects these results, including the failure of
both error bounds for the S-scheme. Its only premise is the mesh-uniform
Crandall--Pazy estimate of the W-scheme, \cite[(3.15)]{Main} with
$\alpha=1$ on uniform meshes.
Remark~\ref{rem:leanscope} states what the Lean proofs cover. The Lean
statements and their translation are given in Appendix~\ref{app:lean}.

\paragraph{Main results.}
On the data $f_h$ the linear error of the S-scheme satisfies
$\norm{(\mathring T_h(t)-T(t))f_h}_2\to e^{-\mathcal C_St}$, uniformly on
compact subsets of $(0,\infty)$ (Corollary~\ref{cor:lineargap}). Here $\mathring T_h$ and $T$ are respectively the
discrete and the continuous semigroups. For data $\varepsilon_hf_h$ with
$\varepsilon_h=o(h)$, the full error divided by $\varepsilon_h$
has the same limit (Theorem~\ref{thm:total}(c)). In contrast, the linear
error of the W-scheme is at most $Ce^{-\nu t}ht_\wedge^{-3/2}$ for $t>0$, uniformly on the unit
ball~\cite[Theorem~4.2]{Main}. For a fixed datum the linear error of the S-scheme
converges to zero, but for every prescribed rate there is a datum on which
this rate fails (Corollary~\ref{cor:fixed}). The ordering error of the
W-scheme is at most $Ch^{1/6}\varepsilon^2$ for all $t\ge0$, uniformly on
an energy ball~\cite[Theorem~4.3]{Main}. For the S-scheme, on the data
$\varepsilon f_h$, the ordering error divided by $\varepsilon^2$ tends to
$-E(t)\ne0$ as first $\varepsilon\downarrow0$ and then $h\to0$
(Corollary~\ref{cor:noproduct}). The limit contains a boundary forcing
from the cell structure of the eigenfunction. Consequently the ordering
error of the S-scheme admits no bound $Ch^p\varepsilon^2$ with $p>0$ that
is uniform in $h$ and $\varepsilon$.

Behind the failure is a general principle
(Proposition~\ref{prop:necessary}). Suppose that a family of discrete
semigroups converges in operator norm, at some positive time, to an
eventually compact semigroup. Suppose also that the inverses of the
discrete generators converge in operator norm to an injective operator.
Then the damping of discrete eigenpairs tends to infinity whenever their
frequency does. For the S-scheme the inverses of the generators converge,
at the rate $h$. The spurious modes
have unbounded frequency, while their damping stays bounded
(Theorem~\ref{thm:top}). Hence the semigroups of the S-scheme do not
converge in operator norm at any positive time
(Corollary~\ref{cor:principleS}).

\paragraph{What the gap does not say.}
For a fixed datum the linear errors of both schemes converge. On smooth
fixed data the computations suggest that the S-scheme can even be more accurate than
the W-scheme (Appendix~\ref{app:num}). The gain of mesh-uniform regularity lies in
estimates that hold uniformly on balls of data of low regularity and in
the norm of the datum.

\paragraph{An application.}
A nonlinear correction computed on a coarse mesh can be added to a linear
solution on a fine mesh (Section~\ref{sec:twomesh}). The error of the sum
is the linear error on the fine mesh plus the ordering error on the coarse
mesh. Under the assumptions of Section~\ref{sec:twomesh}, the linear
solution on the fine mesh has an error of order $\varepsilon^2$, and the
nonlinear solution on the coarse mesh, of size $\hat h$, has an error of
order $\hat h^q\varepsilon$ for some $q>0$. A bound $C\hat h^p\varepsilon^2$
of the ordering error is smaller in order than both when $\hat h^p\ll1$
and $\varepsilon\ll\hat h^{q-p}$. For the W-scheme, the bound
of~\cite[Theorem~4.3]{Main} has this form with $p=1/6$. For the S-scheme no
such bound holds (Corollary~\ref{cor:noproduct}). The computed examples of
Section~\ref{sec:twomesh} with the W-scheme show that the two-mesh sum has
an accuracy comparable to that of the nonlinear solution on the fine mesh,
while taking about half as long to compute.

\paragraph{Related work.}
High-frequency numerical waves of vanishing group velocity are known for
finite-difference schemes~\cite{Trefethen}. For the wave equation they
escape boundary observation and damping~\cite{InfanteZuazua,Zuazua05}.
Uniform decay can be restored by numerical viscosity or other
modifications~\cite{TTZ,RTT,ErvedozaZuazua}, and in control problems
two-grid algorithms filter out these modes~\cite{Glowinski92,IgnatZuazua}.
For wave equations with boundary damping, Banks, Ito and
Wang~\cite{BanksItoWang} proved uniform exponential stability for
particular schemes and gave numerical evidence that several standard
schemes lose it. Here both discrete semigroups are contractive. The
difference that we prove lies in their regularity.

Two-grid methods~\cite{Xu94,Xu96,MarionXu,DWW} combine a nonlinear
computation on a coarse mesh with simpler computations on a fine mesh.
For evolution problems their error estimates hold on bounded time
intervals, or on $[t_0,\infty)$ with constants depending on $t_0$ and on
the datum~\cite{MarionXu,DWW}. The two-mesh sum of
Section~\ref{sec:twomesh} is built differently: its linear part on the
fine mesh does not use the coarse solution, so the parts can be computed
in parallel. Its error is the linear error on the fine mesh plus the
ordering error on the coarse mesh. For the W-scheme the latter is bounded
by $C\hat h^{1/6}\varepsilon^2$ for all $t\ge0$~\cite[Theorem~4.3]{Main}.

For parabolic problems, error estimates for nonsmooth initial data are
classical~\cite[Chapters~3 and~4]{Thomee}. For a selfadjoint problem, the
Galerkin semidiscretization inherits the analyticity of the semigroup
uniformly in $h$~\cite[Lemma~3.9]{Thomee}. The estimate of the W-scheme is
of the same type, with Crandall--Pazy regularity in place of analyticity.
The S-scheme does not have this mesh-uniform regularity.

\paragraph{Organization.}
Section~\ref{sec:setting} introduces the two schemes, the errors,
mesh-uniform regularity and high-frequency spurious modes.
Section~\ref{sec:spurious} contains the results for which Lean is used.
It states the existence of the high-frequency spurious eigenpairs, the
failure of mesh-uniform regularity and intermediate results of the Lean
proof. It then states, with sketches, the properties of the
eigenfunctions, the inverse of the discrete operator, the Galerkin limits
and the linear error gap. Section~\ref{sec:gap} proves the error gap for
the nonlinear problem by hand. Section~\ref{sec:twomesh} treats the
coarse-mesh correction. Appendix~\ref{app:lean}, which may be skipped,
gives the Lean statements of the theoretical results and of their only
premise, with their translation and tables of the corresponding Lean
declarations. It also outlines the Lean proof of Theorem~\ref{thm:top}
and describes the Lean files. Appendix~\ref{app:num} describes the
computations. Appendix~\ref{app:plain}
explains the computations for readers outside numerical analysis.
Table~\ref{tab:notation} collects the main notation. In Sections~\ref{sec:eigfun}--\ref{sec:gap}, $L=1$, $h=1/N$ and $N\ge512$.

\begin{table}[H]
\centering\small
\caption{Notation.}\label{tab:notation}
\begin{tabular}{@{}l>{\raggedright\arraybackslash}p{0.56\textwidth}l@{}}
\toprule
Symbol & Meaning & Defined in\\
\midrule
$L$; $N$, $h$ & length of the interval $(0,L)$; number of elements of the uniform mesh, $h=L/N$ & Section~\ref{sec:schemes}\\
$A$, $T(t)$ & generator $-\partial_x^3$ and semigroup of the linearized problem & Section~\ref{sec:schemes}\\
$X_h$; $H$ & trial space of $C^1$ piecewise cubics with the boundary conditions; weighted space $L^2((0,L),(L+x)\,dx)$ & Section~\ref{sec:schemes}\\
$(f,g)$; $(f,g)_2$, $\norm f_2$ & bilinear pairing $\int_0^Lfg\,dx$; inner product and norm of $L^2(0,L)$ & Section~\ref{sec:schemes}\\
$A_h$, $T_h(t)$, $\Pi_h$ & W-scheme: weighted Galerkin operator, its semigroup, $H$-orthogonal projection & Section~\ref{sec:schemes}\\
$\mathring A_h$, $\mathring T_h(t)$, $\mathring\Pi_h$ & S-scheme: Galerkin operator, its semigroup, $L^2$-orthogonal projection & \eqref{eq:Sscheme}\\
$u$, $v$; $u_h$, $v_h$; $\mathring u_h$, $\mathring v_h$ & nonlinear and linear solutions: continuous, W-scheme, S-scheme & Sections~\ref{sec:schemes}, \ref{sec:errors}\\
$\mathcal G_h$ & $(\mathcal G_h[f],\chi)=\tfrac12(f,\chi')$ for $\chi\in X_h$; $\mathcal G_h[\mathring u_h^2]$ is the nonlinear term of the S-scheme & Section~\ref{sec:schemes}\\
$L_h(t)$, $\mathring L_h(t)$ & linear error operators & Section~\ref{sec:errors}\\
$\Cd_h$, $\mathring\Cd_h$ & ordering errors & Section~\ref{sec:errors}\\
$\zeta_h$, $\varphi_h$; $\mu_h$ & high-frequency spurious eigenpair of the S-scheme, $\norm{\varphi_h}_2=1$; $\mu_h=\Rea\zeta_h$ & Section~\ref{sec:eigfun}\\
$\mathcal C_S$ & $\lim_{N\to\infty}\bigl(-L^3\Rea\zeta_h\bigr)$ & \eqref{eq:beta}, \eqref{eq:re}\\
$\tilde f_h$, $f_h$ & $\tilde f_h=\sqrt2\Rea\varphi_h$, $f_h=\tilde f_h/\norm{\tilde f_h}_2$ & \eqref{eq:fh}\\
$\mathfrak d$, $\mathfrak p$ & $2\sin^2\pi x$ and $x^2-x$ & \eqref{eq:fh}\\
$\mathsf m$, $\mathsf m_h$ & $\mathsf m(r)=r-3r^2+2r^3$; grid mode $\mathsf m_h(x)=h\,\mathsf m(x/h)$, $\mathsf m$ extended with period one & Section~\ref{sec:eigfun}\\
$H_0,\dots,H_3$ & $H_0=1-3r^2+2r^3$, $H_1=r-2r^2+r^3$, $H_2=3r^2-2r^3$, $H_3=-r^2+r^3$ & \eqref{eq:Hermite}\\
$\varpi$; $\xi_*$, $\varpi_*$ & optical branch of the dispersion relation \eqref{eq:dispersion}; $\xi_*\in(0,\pi)$ with $\varpi'(\xi_*)=0$, $\varpi_*=\varpi(\xi_*)$ & Section~\ref{sec:disp}\\
$z_*$, $\mathsf v_*$ & $z_*=e^{i\xi_*}$; $\mathsf v_*=\varpi_*(312+108\cos\xi_*)/(2-2\cos\xi_*-26\varpi_*\sin\xi_*)$ & Sections~\ref{sec:disp}, \ref{sec:eigfun}\\
$p$ & $p(r)=H_0(r)+i\mathsf v_*H_1(r)+z_*\bigl(H_2(r)+i\mathsf v_*H_3(r)\bigr)$, $0\le r\le1$ & \eqref{eq:bdef}\\
$\mathfrak g$, $\mathfrak b$ & $\mathfrak g=|p|^2/\int_0^1|p|^2$, extended with period one; $\mathfrak b=\tfrac12\int_0^1\mathfrak g\,\mathsf m'\,dr$ & \eqref{eq:bdef}\\
$z_\infty$ & $-\mathfrak d/(8\pi^2)+\mathfrak b\mathfrak p$ & \eqref{eq:zinf}\\
$\mathcal U$ & nonlinear solution as a Galerkin limit of the W-scheme & Proposition~\ref{prop:Ulimit}\\
$P_2$, $\mathring P_{2,h}$, $\mathring D_h$ & quadratic Taylor coefficients of the solution maps; discretization error of this coefficient for the S-scheme & Section~\ref{sec:quad}\\
$E(t)$ & limit of $\mathring P_{2,h}(\tilde f_h)(t)$ as $h\to0$ (Proposition~\ref{thm:quad}(a)) & \eqref{eq:Edef}\\
\bottomrule
\end{tabular}
\end{table}

\section{Setting}\label{sec:setting}

In this section, we introduce the two schemes, the errors, mesh-uniform
regularity and high-frequency spurious modes.

\subsection{The two schemes}\label{sec:schemes}
Let $L>0$, let $A=-\partial_x^3$ with
$D(A)=\{f\in H^3(0,L):f(0)=f(L)=f'(L)=0\}$, and let $T(t)=e^{tA}$. The
linearization of \eqref{eq:pde} is $v_t=Av$. For $f\in L^2(0,L)$, the
problem $-\varrho'''=f$, $\varrho(0)=\varrho(L)=\varrho'(L)=0$ has a unique
solution $\varrho\in H^3(0,L)$. It is given by variation of constants, and
$\norm\varrho_{H^3}\le C\norm f_2$. Hence $A$ is invertible. Since
$H^3(0,L)$ embeds compactly in $L^2(0,L)$, the generator $A$ has compact
resolvent. Both schemes use the space
$X_h$ of real $C^1$ piecewise cubic polynomials on a mesh of size $h$ that
satisfy the three boundary conditions. Eigenpairs are taken in its
complexification. Write $(f,g)=\int_0^Lfg\,dx$ for the bilinear
pairing and $(f,g)_2=\int_0^Lf\overline g\,dx$, $\norm f_2^2=(f,f)_2$.

\emph{The W-scheme} is the weighted scheme considered in~\cite{Main}. It works in
$H=L^2((0,L),(L+x)\,dx)$. Its weighted Galerkin operator $A_h$ on $X_h$
is the Galerkin form of $A$ tested in $H$; see~\cite[Section~3.1]{Main}.
The weighted energy identity of~\cite[Section~2.1]{Main} holds exactly on
$X_h$: $\Rea(-A_hv,v)_H=\tfrac32\norm{v'}_2^2+\tfrac L2|v'(0)|^2$. Its
interior term is coercive.
We write $T_h(t)=e^{tA_h}$, $\Pi_h$ for the $H$-orthogonal projection onto
$X_h$, and $u_h$, $v_h=T_h(\cdot)\Pi_hu_0$ for the nonlinear and the linear
discrete solutions.

\emph{The S-scheme} is the standard Galerkin method on $X_h$ in
the unweighted space $L^2(0,L)$. A ring accent marks quantities of the
S-scheme. For example, $\mathring T_h$ and $\mathring u_h$ correspond to
$T_h$ and $u_h$ of the W-scheme. The operator $\mathring A_h$ on $X_h$ is defined by
\begin{equation}\label{eq:Sscheme}
 (\mathring A_hv,\chi)_2=\int_0^Lv''\,\overline{\chi'}\,dx\qquad(v,\chi\in X_h),
\end{equation}
and integration by parts gives
\begin{equation}\label{eq:flux}
 2\Rea(\mathring A_hv,v)_2=|v'(L)|^2-|v'(0)|^2=-|v'(0)|^2 .
\end{equation}
Thus $\mathring T_h(t)=e^{t\mathring A_h}$ is a contraction on $L^2$, but the
dissipation acts only at the boundary. Every eigenpair
$\mathring A_h\varphi=\zeta\varphi$ satisfies
$|\varphi'(0)|^2=2|\Rea\zeta|\norm\varphi_2^2$ and $\Rea\zeta\le0$.
For $f\in L^1(0,L)$ let $\mathcal G_h[f]\in X_h$ be defined by
$(\mathcal G_h[f],\chi)=\tfrac12(f,\chi')$ for real $\chi\in X_h$, extended
complex linearly in $f$. The nonlinear S-scheme is
\[
 \mathring u_{h,t}=\mathring A_h\mathring u_h+\mathcal G_h[\mathring u_h^2],\qquad
 \mathring u_h(0)=\mathring\Pi_hu_0,
\]
where $\mathring\Pi_h$ is the $L^2$-orthogonal projection onto $X_h$. The
term $\mathcal G_h[\mathring u_h^2]$ is the Galerkin form of $-\mathring
u_h\mathring u_{h,x}$. For real data it has a global solution: it is a
polynomial system of ordinary differential equations, by uniqueness the
solution for real data is real, and
$(\mathcal G_h[v^2],v)=\tfrac12(v^2,v')=0$ for real $v$, together with \eqref{eq:flux},
bounds $\norm{\mathring u_h}_2$. We put $\mathring v_h=\mathring T_h(\cdot)\mathring\Pi_hu_0$.

For the S-scheme we use uniform meshes, $x_j=jh$, $h=L/N$. The continuous semigroup
$T$ is a contraction on $L^2$ by the unweighted energy identity, and on
$H$~\cite[Proposition~2.1]{Main}. The two norms are equivalent:
$\norm f_2\le L^{-1/2}\norm f_H\le\sqrt2\norm f_2$.

\paragraph{An inverse inequality.}
For any quadratic polynomial $p$ on $[0,1]$, expand in the orthonormal
basis $1$, $\sqrt3(2r-1)$, $\sqrt5(6r^2-6r+1)$, with coefficients
$a_0,a_1,a_2$. Direct integration gives
\[
 \|p'\|_{L^2(0,1)}^2=12|a_1|^2+60|a_2|^2
 \le60\sum_{k=0}^2|a_k|^2=60\|p\|_{L^2(0,1)}^2.
\]
Since $60\le8^2$, apply this to $p(r)=v'(x_j+hr)$ and sum over the
elements. This proves
\begin{equation}\label{eq:inverse}
 \|v''\|_2\le 8h^{-1}\|v'\|_2\qquad(v\in X_h).
\end{equation}

\subsection{Errors}\label{sec:errors}
Let $u$ solve \eqref{eq:pde} and $v=T(\cdot)u_0$. The \emph{linear errors}
are $v-v_h=-L_h(\cdot)u_0$ and $v-\mathring v_h=-\mathring L_h(\cdot)u_0$, with
\[
 L_h(t)=T_h(t)\Pi_h-T(t),\qquad \mathring L_h(t)=\mathring T_h(t)\mathring\Pi_h-T(t).
\]
The \emph{full errors} are $u-u_h$ and $u-\mathring u_h$. The
\emph{total errors} are $u-v_h$ and $u-\mathring v_h$. A total error
splits in two ways,
\[
 u-v_h=(u-v)+(v-v_h)=(u-u_h)+(u_h-v_h).
\]
The \emph{ordering error}, introduced in~\cite{Main}, is
\[
 \Cd_h=(u-v)-(u_h-v_h)=(u-u_h)-(v-v_h),
\]
and $\mathring\Cd_h$ is defined in the same way. It measures how much the
nonlinearity changes the discretization error.

\subsection{Mesh-uniform regularity and high-frequency spurious modes}\label{sec:mur}
Let a family of discrete semigroups $T_h$ with generators $A_h$
approximate a semigroup $T$. The family has \emph{mesh-uniform regularity}
if the semigroups $T_h$ inherit a regularity property of $T$, with bounds
independent of $h$. Here regularity refers to the regularity of semigroups
associated with linearized equations. This is one concept with many kinds. Examples are
analyticity, Gevrey regularity, Crandall--Pazy regularity and Kato
smoothing. Any one kind is a mesh-uniform regularity. Here two kinds
matter.
\begin{enumerate}
\item[(i)] Mesh-uniform Crandall--Pazy regularity of index $3/2$:
$\sup_h\sup_{0<t\le1}t^{3/2}\norm{A_hT_h(t)}<\infty$.
\item[(ii)] Mesh-uniform Kato smoothing:
$\int_0^\infty\norm{(T_h(t)v)'}_2^2\,dt\le C\norm v^2$ for $v\in X_h$,
with $C$ independent of $h$.
\end{enumerate}
The W-scheme has property (i)~\cite[Theorem~3.5]{Main} and property
(ii)~\cite[Section~3.2]{Main}. The S-scheme has neither
(Corollary~\ref{cor:failures}).

A family of eigenpairs of the discrete operators whose frequencies tend to
infinity is called a family of \emph{high-frequency spurious modes} if its
eigenvalues asymptotically deviate too much from the region where the
eigenvalues of the continuous operator lie. For the present problem this
region is parabolic. By~\cite[Theorem~2.4(i)]{Main},
$\norm{(i\tau-A)^{-1}}_H\le C\langle\tau\rangle^{-2/3}$ for $\tau\in\R$,
where $\langle\tau\rangle=(1+\tau^2)^{1/2}$. If $A\varphi=\zeta\varphi$
and $\tau=\Ima\zeta$, then $(i\tau-A)\varphi=-(\Rea\zeta)\varphi$, so
$|\Rea\zeta|\ge1/\norm{(i\tau-A)^{-1}}_H$. Hence, for any $\Omega>0$, the
eigenvalues of $A$ with $|\Ima\zeta|\ge\Omega$ lie in the region
\begin{equation}\label{eq:region}
 |\Rea\zeta|\ge c\,|\Ima\zeta|^{2/3},\qquad |\Ima\zeta|\ge\Omega,
\end{equation}
with $c=1/C$. The resolvent of the W-scheme satisfies the same bound,
with $C$ independent of $h$~\cite[Theorem~3.5]{Main}. The same argument
puts its eigenvalues in \eqref{eq:region}, so the W-scheme has no
high-frequency spurious modes. The S-scheme has them. For $L=1$ the
eigenvalues $\zeta_h$ given in Theorem~\ref{thm:top} satisfy
$|\Ima\zeta_h|\to\infty$ and
$|\Rea\zeta_h|\,|\Ima\zeta_h|^{-2/3}\le1600\,(100h^{-3})^{-2/3}\le75h^2$,
so no constants $c,\Omega>0$ make \eqref{eq:region} hold for them as $h\to0$.

\section{High-frequency spurious eigenpairs of the S-scheme}\label{sec:spurious}

In this section, we state the results on the eigenpairs of the S-scheme
and the linear approximation results. In
Sections~\ref{sec:disp}--\ref{sec:leanlemmas} we give the definitions and
the statements only. The proof of Theorem~\ref{thm:top} is a Lean~4 proof,
outlined in Appendix~\ref{app:lean}, which may be skipped.
Section~\ref{sec:leanlemmas} states intermediate results of this proof.
In Sections~\ref{sec:eigfun}--\ref{sec:lineargap} each statement is
followed by a sketch. For the linear approximation results the sketch
contains all the steps and can be completed by hand. For the properties
of the eigenfunctions it gives the idea and the key quantities.
Remark~\ref{rem:leanscope} describes the role of Lean.
The mesh is uniform, $x_j=jh$, $h=L/N$.

\subsection{Grid waves and constants}\label{sec:disp}
\paragraph{Element matrices.}
On an element write $r=(x-x_j)/h\in[0,1]$ and $d_j=hu'(x_j)$.
In the nodal order $(u_j,d_j,u_{j+1},d_{j+1})$, the four cubic Hermite
basis functions are
\begin{equation}\label{eq:Hermite}
 \begin{aligned}
 H_0(r)&=1-3r^2+2r^3,& H_1(r)&=r-2r^2+r^3,\\
 H_2(r)&=3r^2-2r^3,& H_3(r)&=-r^2+r^3.
 \end{aligned}
\end{equation}
Direct integration gives the element matrices (rows index test functions):
\begin{equation}\label{eq:element}
 M_e=\frac h{420}\underbrace{\begin{pmatrix}156&22&54&-13\\22&4&13&-3\\
 54&13&156&-22\\-13&-3&-22&4\end{pmatrix}}_{\widehat M},\qquad
 K_e=\frac1{h^2}\underbrace{\begin{pmatrix}0&1&0&-1\\-1&-\frac12&1&-\frac12\\
 0&-1&0&1\\1&\frac12&-1&\frac12\end{pmatrix}}_{\widehat K}.
\end{equation}

\paragraph{Dispersion relation.}
At an interior node, put $(u_j,d_j)=e^{ij\xi}z$, $c=\cos\xi$ and
$s=\sin\xi$. Assembly gives the dimensionless Bloch matrices
\begin{equation}\label{eq:bloch}
 K_0=\begin{pmatrix}0&2(1-c)\\-2(1-c)&-is\end{pmatrix},\qquad
 M_0=\begin{pmatrix}312+108c&-26is\\26is&8-6c\end{pmatrix}.
\end{equation}
The assembled rows are $h^{-2}K_0z$ and $(h/420)M_0z$. Put
\[
 a=455-252c+7c^2,\quad b=s(104+c),\quad d=(1-c)^2,\quad \Delta=b^2+4ad.
\]
Expansion of the determinant gives
\begin{equation}\label{eq:dispersion}
 \det(K_0-i\omega M_0)=-4P(\xi,\omega),\qquad
 P(\xi,\omega)=a\omega^2+b\omega-d.
\end{equation}
A root $\omega$ gives the interior spectral value $\zeta=420i\omega/h^3$.
The negative root
\[
 \varpi(\xi)=\frac{-b-\sqrt\Delta}{2a}
\]
is the optical branch. The group velocity of the corresponding grid wave
is proportional to $\varpi'(\xi)$.

\paragraph{The point of zero group velocity.}
Let
\[
 \tilde H(c)=-7644c^4+9790c^3+237409c^2+335938c-244743,\qquad
 R_c(c)=\frac{94-128c+34c^2}{-3809+6968c+131c^2-140c^3}.
\]
Let $c_*$ be the zero of $\tilde H$ in
$[c_{\rm lo},c_{\rm hi}]$, where
$c_{\rm lo}=528567235794090/10^{15}$ and
$c_{\rm hi}=528567235794092/10^{15}$. It is unique there, since
$\tilde H'>0$ on $[0,1]$ ($\tilde H''>0$ on $[0,1]$ and
$\tilde H'(0)>0$). Put
\[
 s_*=(1-c_*^2)^{1/2},\qquad \xi_*=\arccos c_*,\qquad w_*=R_c(c_*)s_*.
\]
Then
\begin{equation}\label{eq:zgv}
 0<\xi_*<\pi,\qquad \varpi(\xi_*)=w_*,\qquad \varpi'(\xi_*)=0,
 \qquad \varpi''(\xi_*)>0 .
\end{equation}
These are the Lean statements \lean{thetaStar_range},
\lean{optical_thetaStar}, \lean{optical_stationary} and
\lean{optical_curvature_pos} (Appendix~\ref{app:lean}). Thus $\xi_*$ is a
point of zero group velocity.

\paragraph{Constants.}
Set
\[
 \varpi_*=\varpi(\xi_*),\quad \varpi''=\varpi''(\xi_*),\quad
 \Omega_*=420|\varpi_*|,\quad z_*=e^{i\xi_*},\quad \eta_*=i\varpi_*.
\]
Numerically $\varpi_*\approx-0.27656$ and $\Omega_*\approx116.15$
(Appendix~\ref{app:num}). The interior characteristic polynomial, the
slope ratio and the left boundary functional are
\begin{align}
 q(z;\eta)&=(z-1)^4+\eta(z^4+208z^3-208z-1)+\eta^2(7z^4-504z^3+1834z^2-504z+7),\label{eq:q}\\
 \mathcal V(z,\eta)&=\frac{\eta\,(312+54\mathfrak s)}{\mathcal D(z,\eta)},\qquad \mathcal D=(2-\mathfrak s)+13\eta\mathfrak t,\qquad
 \mathfrak s=z+z^{-1},\ \mathfrak t=z-z^{-1},\label{eq:slope}\\
 B(z,\eta)&=z-\tfrac12(1+z)\mathcal V-\eta\,(4\mathcal V+13z-3z\mathcal V),\qquad \mathcal V=\mathcal V(z,\eta).\label{eq:B}
\end{align}
Let $z_3$ be the zero of $q(\cdot;\eta_*)$ in $\abs z<1/50$ and
$z_4=1/\overline{z_3}$. Lemma~\ref{lem:cert} in
Section~\ref{sec:leanlemmas} gives their existence and uniqueness. Put $\mathcal V_*=\mathcal V(z_*,\eta_*)$, $B_*=B(z_*,\eta_*)$,
$\mathcal V_4=\mathcal V(z_4,\eta_*)$, $B_3=B(z_3,\eta_*)$, and let
$\mathcal V'=iz_*\partial_z\mathcal V(z_*,\eta_*)$,
$B'=iz_*\partial_zB(z_*,\eta_*)$ be the $\xi$-derivatives at fixed
$\eta_*$. Define
\begin{equation}\label{eq:beta}
 \beta=-\frac{2B'}{B_*-B_3}+\frac{2\mathcal V'}{\mathcal V_*-\mathcal V_4},\qquad \mathcal C_S=210\pi^2\varpi''\Rea\beta .
\end{equation}

\subsection{The high-frequency spurious eigenpairs}
At $\xi_*$ the group velocity vanishes. Heuristically, a grid wave of
frequency $\Omega_*h^{-3}$ does not travel to the dissipative boundary.
The following theorem, the main result of this section, turns this into
an exact eigenpair.

\begin{theorem}[high-frequency spurious eigenpairs]\label{thm:top}
Fix $L>0$. Then the following hold.
\begin{enumerate}
\item[(a)] For every integer $N\ge512$ and $h=L/N$, the operator
$\mathring A_h$ has an eigenpair $(\zeta_h,\varphi_h)$, $\varphi_h\ne0$, with
\begin{equation}\label{eq:bounds}
 |\Ima\zeta_h|\ge100h^{-3},\qquad 0<-\Rea\zeta_h\le1600/L^3.
\end{equation}
\item[(b)] For the same $N$, the same eigenvalue $\zeta_h$ satisfies
\begin{align}
 \Bigl|h^3|\Ima\zeta_h|-\Bigl(\Omega_*-\frac{210\pi^2\varpi''}{N^2}
 +\frac{210\pi^2\varpi''\Ima\beta}{N^3}\Bigr)\Bigr|&\le\frac{1\,260\,000}{N^4},\label{eq:im}\\
 \bigl|L^3\Rea\zeta_h+\mathcal C_S\bigr|&\le\frac{1\,260\,000}{N},\label{eq:re}
\end{align}
where $\Omega_*$, $\varpi''$, $\beta$ and $\mathcal C_S$ are the constants
of Section~\ref{sec:disp}.
\end{enumerate}
\end{theorem}

The constants satisfy $\Omega_*>100$ and $0<\mathcal C_S<1600$
(Lemma~\ref{lem:cert}). Numerically, $\mathcal C_S=1575.383\ldots$
(Appendix~\ref{app:num}). The scheme \eqref{eq:Sscheme} has real
coefficients. Hence the conjugate $(\overline{\zeta_h},\overline{\varphi_h})$
is also an eigenpair. It satisfies the same bounds \eqref{eq:bounds}.

The constant $1\,260\,000$ in Theorem~\ref{thm:top}(b) is not optimized.
The computations in Appendix~\ref{app:num} suggest that on the meshes
computed there the errors are about $6000N^{-4}$, and about $700N^{-1}$
for $N=256,512$. Theorem~\ref{thm:top} concerns only the first optical
mode. A numerical scan of the
meshes $3\le N\le128$ and $N=256,512$ finds the bounds \eqref{eq:bounds}
satisfied from $N=7$ on (Appendix~\ref{app:num}).

The proof of Theorem~\ref{thm:top} is a Lean~4 proof. The following
remark explains how it was obtained and what the author checked.

\begin{remark}[how Theorem~\ref{thm:top} was obtained]\label{rem:history}
The first version of part~(a) had no explicit threshold: it held for all
sufficiently large $N$. Its proof was produced by Claude Opus~5.5 and
was computer-assisted. Rouch\'e's theorem and the argument principle
located the roots, and rational interval arithmetic in Python checked
the nondegeneracy conditions. The author tried to check this proof by
hand and found its complexity too large for a complete check. The
author therefore had the proof formalized in Lean~4. The formalization
was carried out by GPT-6~Astra. Mathlib has no Rouch\'e theorem. The
formal proof therefore uses contraction mappings and analytic inverse
functions. A contraction mapping comes with explicit constants. Hence
the threshold could be made explicit. In successive iterations, each ending in a
complete Lean theorem, it was lowered from $10^{200}$ to
$1.05\cdot10^{85}$, $2048$, $1024$ and finally $512$. Part~(b) was first
derived by hand from the Lean fixed point, with remainders
$O(N^{-4})$ and $O(N^{-1})$. Its formalization then gave the explicit
constants for $N\ge512$.

Thus the present statement of Theorem~\ref{thm:top} and its present
proof are both results of these iterations. The first proof is not used
in this paper. The present proof is the Lean proof of the Lean statements
\lean{SLean.Explicit.Sharp.Expansion.spectral_top_with_expansion} in the
ancillary files of arXiv:2611.12585. Lean checks it. GPT-6~Astra rewrote the proof along
the Lean route. Appendix~\ref{app:lean}, which may be skipped, outlines
this proof with the Lean name of each step.

The author does not check the Lean code directly. The Lean statements
were translated into mathematics by a fixed procedure. Claude Opus~5.5
and GPT-6~Astra each drafted a translation from the Lean source, without
seeing the other's draft. The translation covers the Lean statements and
the definitions on which these statements depend. Each then
checked the other's draft against the Lean source and reported in
writing. A single translation was written from the two drafts only
after both reports, and Claude Fable~5.1 checked this translation
against the Lean source. On this basis, the author made sure that the
translated Lean statements are the statements of Theorem~\ref{thm:top}.
Appendix~\ref{app:lean} reproduces the Lean statements with this
translation, so that a reader who knows Lean can check it directly.
\end{remark}

The eigenpairs rule out mesh-uniform regularity. The following corollary
has a short proof by hand and a Lean~4 proof (Appendix~\ref{sec:A5}).

\begin{corollary}[failure of mesh-uniform regularity]\label{cor:failures}
All operator norms below are $L^2$ norms.
\begin{enumerate}
\item[(a)] A mesh-uniform Crandall--Pazy bound
$\norm{\mathring A_he^{t\mathring A_h}}_2\le Ct^{-3/2}$ for $0<t\le1$ fails.
Indeed, for every fixed $t\in(0,1]$ and $N\ge512$,
\begin{equation}\label{eq:CPlower}
 \norm{\mathring A_he^{t\mathring A_h}}_2\ge100h^{-3}e^{-1600t/L^3}.
\end{equation}
\item[(b)] For every fixed $T>0$, a mesh-uniform finite-time Kato bound
fails, already for real data. For each $N\ge512$ there is a real $\chi_h\ne0$ such that
\begin{equation}\label{eq:katofinite}
 \int_0^T\norm{(e^{t\mathring A_h}\chi_h)'}_2^2\,dt
 \ge\frac{L^3}{256}(1-e^{-3200T/L^3})h^{-2}\norm{\chi_h}_2^2.
\end{equation}
\end{enumerate}
\end{corollary}

\begin{proof}
\emph{(a)} On the eigenvector,
$\mathring A_he^{t\mathring A_h}\varphi_h=\zeta_he^{t\zeta_h}\varphi_h$.
By \eqref{eq:bounds}, we obtain \eqref{eq:CPlower}. A bound $Ct^{-3/2}$
would be finite at any fixed $t\in(0,1]$, which contradicts
\eqref{eq:CPlower} as $h\to0$.

\emph{(b)} We first bound $\|\varphi_h'\|_2$ from below. The eigenvalue
equation and the inverse inequality \eqref{eq:inverse} yield
\[
 |\zeta_h|\|\varphi_h\|_2^2\le 8h^{-1}\|\varphi_h'\|_2^2,
 \qquad \|\varphi_h'\|_2^2\ge\frac{100}{8}h^{-2}\|\varphi_h\|_2^2.
\]
Let $\mu_h=\Rea\zeta_h$. By \eqref{eq:bounds}, $-1600/L^3\le\mu_h<0$.
Since $e^{t\mathring A_h}\varphi_h=e^{t\zeta_h}\varphi_h$ and
$e^{2\mu_ht}\ge e^{-3200t/L^3}$, integration over $[0,T]$ yields
\eqref{eq:katofinite} for $\varphi_h$. No lower bound on $|\mu_h|$ is needed.
We write $\varphi_h=p+iq$ with $p,q$ real. The initial squared norm and
the integrated derivative squared norm each split into the contributions
of $p$ and $q$, because the semigroup is real. It follows that at least
one nonzero real component has at least the same quotient. We take it as
$\chi_h$.
\end{proof}

By~\cite[Proposition~3.4]{Main}, part (a) is equivalent to the failure of
a mesh-uniform resolvent bound of order $|\tau|^{-2/3}$ for large $|\tau|$.
In fact, \eqref{eq:CPlower} shows that
$\sup_h\|\mathring A_h\mathring T_h(t)\|_2=\infty$ at every fixed
$t\in(0,1]$. It therefore excludes mesh-uniform Crandall--Pazy bounds of
every finite index. It also excludes the usual mesh-uniform analyticity
and Gevrey bounds. Part (b) excludes mesh-uniform Kato smoothing.
The finite-time Kato bound in (b) is weaker than the bound on
$(0,\infty)$ in the definition of mesh-uniform regularity, so (b) excludes
that bound as well. The damping $-\Rea\zeta_h$ of the eigenpairs stays
bounded, so the mode survives for times independent of the mesh. This is
what the error estimates of Section~\ref{sec:gap} require.

\subsection{Intermediate results of the Lean proof}\label{sec:leanlemmas}
The following results are steps of the Lean proof of
Theorem~\ref{thm:top}. They are stated without proof. Their Lean
names are listed in Table~\ref{tab:corr}. Sections~\ref{sec:eigfun}
and~\ref{sec:inverse} use them.

\begin{lemma}[certificates]\label{lem:cert}
The quantities below are defined in Section~\ref{sec:disp}. They satisfy
the following.
\begin{enumerate}
\item[(C1)] $a(\xi_*)>0$, $\Delta(\xi_*)>0$ and
\[
 -.277<w_*<-.276,\qquad
 \frac{372284}{10^6}\le\varpi''\le\frac{372286}{10^6},\qquad
 \frac{279038}{10^6}\le\frac{\sqrt{\Delta(\xi_*)}}{a(\xi_*)}\le\frac{279039}{10^6}.
\]
\item[(C2)] $z_*$ is a zero of $q(\cdot;\eta_*)$ of multiplicity exactly
two.
\item[(C3)] For $|\eta-\eta_*|\le1/100$, there is exactly one zero in
$|z|<1/50$, no zero on $|z|=1/50$, and exactly one zero in $|z|>50$.
The inner and outer zeros are simple.
\item[(C4)] At $\eta_*$, the small zero satisfies $|z_3|<9/1000$;
the large zero is $z_4=1/\overline{z_3}$.
\item[(C5)] $\mathcal D$ does not vanish at $(z_*,\eta_*)$,
$(z_3,\eta_*)$, $(z_4,\eta_*)$, and
\[
 |B_*-B_3|\ge411/50,\qquad |\mathcal V_*-\mathcal V_4|\ge1587/100,\qquad
 1276/625\le\Rea\beta\le10209/5000.
\]
\end{enumerate}
In particular, $\Omega_*>100$ and $0<\mathcal C_S<1600$.
\end{lemma}

\paragraph{Exact eigenvectors.}
Use the nodal values $u_j=u(x_j)$ and scaled slopes $d_j=hu'(x_j)$, and
put $\eta=\zeta h^3/420$.
\begin{lemma}[determinant criterion]\label{lem:det}
Let $N\ge3$ and let $z_1,\dots,z_4$ be four distinct nonzero zeros of
$q(\cdot;\eta)$ with $\mathcal D(z_k,\eta)\ne0$. Put
$\mathcal V_k=\mathcal V(z_k,\eta)$, $B_k=B(z_k,\eta)$. Let $\gamma_k\ne0$
and let $\D$ be the $4\times4$ matrix with columns
$\gamma_k\,(1,\ B_k,\ z_k^N,\ z_k^N\mathcal V_k)^T$. If $\D\alpha=0$ with
$\alpha\ne0$, then
\[
 (u_j,d_j)=\sum_k\alpha_k\gamma_kz_k^j\,(1,\mathcal V_k),\qquad 0\le j\le N,
\]
are the nodal data of an eigenfunction of $\mathring A_h$ with eigenvalue
$420\eta/h^3$.
\end{lemma}
The Lean statement \lean{Sdet} states the existence of an eigenfunction.
The displayed nodal formula and its nonvanishing are proved by hand; they
also have a Lean~4 proof (Appendix~\ref{sec:A5}). By
the mode identities (\lean{exact_mode_value}, \lean{exact_mode_slope},
\lean{mode_left}), each sequence $z_k^j(1,\mathcal V_k)$ solves the
interior equations, and $B_k$ is the left boundary functional on it. The
rows of $\D\alpha=0$ are the four boundary conditions. The nodal values
$u_j$ for $j=0,\dots,3$ are $\sum_k\alpha_k\gamma_kz_k^j$, and the
Vandermonde matrix of the distinct $z_k$ is invertible. Since $\alpha\ne0$
and $\gamma_k\ne0$, the sequence is nonzero.

\paragraph{Root branches.}
Let $\tilde\varpi(\varsigma)=\varpi(\xi_*+\varsigma)$, extended to complex
$\varsigma$ by the negative quadratic formula with the principal square
root. Write $\eta(\varsigma)=i\tilde\varpi(\varsigma)$,
$z_+(\varsigma)=e^{i(\xi_*+\varsigma)}$, and
\[
 a_2=\frac{\varpi''}{2},\qquad a_3=\frac{\varpi'''(\xi_*)}{6},
 \qquad c_2=\frac{a_3}{a_2}.
\]
\begin{lemma}[explicit domains]\label{lem:domain}
On $|\varsigma|<1/10$, $\tilde\varpi$ and four root branches $z_\pm$, $z_3$,
$z_4$ of $q(\cdot;\eta(\varsigma))$ are analytic, with $z_3(0)=z_3$ and
$z_4(0)=z_4$. The branch $z_-$ has the form
\[
 z_-=e^{i(\xi_*+\sigma(\varsigma))},
\]
where $\sigma$ is analytic on $|\varsigma|<1/10$ with $\sigma(0)=0$ and
$\sigma'(0)=-1$. On $|\varsigma|\le1/10$,
\[
 |\eta-\eta_*|\le3/1000,\qquad |\sigma|\le2|\varsigma|,\qquad
 |\varsigma+\sigma|\le(81/10)|\varsigma|^2 .
\]
For $0<|\varsigma|\le1/20$ the four roots are distinct and nonzero, their
$\mathcal D$ denominators are nonzero, and, with $B_k$, $\mathcal V_k$ as in
Lemma~\ref{lem:det},
\[
 |B_k|,|\mathcal V_k|\le10^8,\qquad |z_3|,|z_4^{-1}|\le1/50 .
\]
These bounds also hold at $\varsigma=0$. Near $\varsigma=0$,
$\tilde\varpi\circ\sigma=\tilde\varpi$, $B_-=B_+\circ\sigma$ and
$\mathcal V_-=\mathcal V_+\circ\sigma$.
\end{lemma}

\paragraph{Quantization.}
Let $\D_N(\varsigma)$ be the matrix of Lemma~\ref{lem:det} for the zeros
$z_+,z_-,z_3,z_4$ with $\gamma_+=\gamma_-=\gamma_3=1$, $\gamma_4=z_4^{-N}$:
\[
 \D_N=\begin{pmatrix}1&1&1&\delta_4\\ B_+&B_-&B_3&\delta_4B_4\\ Z_+&Z_-&\delta_3&1\\
 Z_+\mathcal V_+&Z_-\mathcal V_-&\delta_3\mathcal V_3&\mathcal V_4\end{pmatrix},
 \qquad Z_\pm=z_\pm^N,\ \delta_3=z_3^N,\ \delta_4=z_4^{-N}.
\]
Put $\mathcal R_+=(B_--B_3)(\mathcal V_+-\mathcal V_4)$,
$\mathcal R_-=(B_+-B_3)(\mathcal V_--\mathcal V_4)$, and let $\mathfrak r=\mathcal R_-/\mathcal R_+$.
\begin{lemma}[reflection]\label{lem:refl}
$\mathfrak r$ is analytic on $|\varsigma|<1/10$ and $\mathfrak r(0)=1$. Its
logarithm $\ell$ with $\ell(0)=0$ satisfies $\ell'(0)=-\beta$ and, for
$|\varsigma|\le1/10$, $|\ell(\varsigma)|\le15|\varsigma|$ and
$|\ell(\varsigma)+\beta\varsigma|\le32|\varsigma|^2$. Moreover,
$\mathcal R_+\ne0$ and $|\mathcal R_+^{-1}|\le1$ for $|\varsigma|\le1/20$.
\end{lemma}

For $N\ge512$ set $\varsigma=(\pi+y)/N$ and
\begin{equation}\label{eq:Phi}
 \Phi_N=\frac{\det\D_N}{\mathcal R_+z_-^N}
 =e^{iN(\varsigma-\sigma)}-\mathfrak r+R_N .
\end{equation}
\begin{lemma}[exponential tail]\label{lem:tail}
On $|y|\le2$, $|R_N|\le N^{-6}$, $|\mathfrak r^{-1}|\le2$, and
$J_N=\Log(1-R_N/\mathfrak r)$ satisfies $|J_N|\le4N^{-6}$.
\end{lemma}

\begin{lemma}[first-mode fixed point]\label{lem:fixed}
For every $N\ge512$ the map
\[
 \mathcal T_N(y)=\frac N2(\varsigma+\sigma(\varsigma))
       +\frac{\ell(\varsigma)+J_N(\varsigma)}{2i},\qquad \varsigma=(\pi+y)/N,
\]
has a unique fixed point $y_N$ in $|y|<1/2$. It satisfies
$|y_N|\le147/N$ and $\Phi_N(\varsigma_N)=0$ for
$\varsigma_N=(\pi+y_N)/N$.
\end{lemma}

\begin{lemma}[parameter box]\label{lem:box}
For $N\ge512$, $|y_N|\le8/N$ and $|\varsigma_N|\le16/(5N)$.
\end{lemma}

\paragraph{The eigenvalue.}
The eigenvalue of Theorem~\ref{thm:top} is
$\zeta_h=420i\tilde\varpi(\varsigma_N)/h^3$. The Lean statements give this
eigenvalue and assert the existence of an eigenfunction. They do not
single out one.

\subsection{The eigenfunctions}\label{sec:eigfun}
The results of this subsection are by-products of Theorem~\ref{thm:top}
and of its proof. They are used in the error estimates of
Section~\ref{sec:gap}. So are Lemmas~\ref{lem:slow} and~\ref{lem:fast} of
Section~\ref{sec:inverse}.

\paragraph{The eigenfunction and the data.}
From here to the end of Section~\ref{sec:gap}, $L=1$ and $h=1/N$ with $N\ge512$. Let $\zeta_h$ be the
eigenvalue of Theorem~\ref{thm:top}. Since $|y_N|\le8/N<\pi$,
Lemma~\ref{lem:box} gives $0<|\varsigma_N|\le16/(5N)<1/20$. By
Lemma~\ref{lem:domain}, the four roots $z_+,z_-,z_3,z_4$ at $\varsigma_N$
are distinct and nonzero, and $\mathcal D$ does not vanish at them. In
particular $\gamma_4=z_4^{-N}\ne0$. By Lemma~\ref{lem:refl},
$\mathcal R_+\ne0$. Hence the denominator in \eqref{eq:Phi} is nonzero. By
Lemma~\ref{lem:fixed}, $\Phi_N(\varsigma_N)=0$, so the matrix
$\D_N(\varsigma_N)$ is singular. For $\varphi_h$ we take the eigenfunction
that Lemma~\ref{lem:det} gives from a nonzero kernel vector of
$\D_N(\varsigma_N)$. Its eigenvalue is $420i\tilde\varpi(\varsigma_N)/h^3$,
which is $\zeta_h$ (Section~\ref{sec:leanlemmas}).
We normalize it by $\norm{\varphi_h}_2=1$ and fix its phase. Put
$\mu_h=\Rea\zeta_h$, so $\mu_h\to-\mathcal C_S$ by \eqref{eq:re}, and
\begin{equation}\label{eq:fh}
 \tilde f_h=\sqrt2\,\Rea\varphi_h,\qquad
 f_h=\tilde f_h/\norm{\tilde f_h}_2,\qquad
 \mathfrak d(x)=2\sin^2\pi x,\qquad \mathfrak p(x)=x^2-x .
\end{equation}
Here $\tilde f_h\ne0$. Otherwise $\Ima\varphi_h$ would be a real
eigenfunction of the real operator $\mathring A_h$ with the nonreal
eigenvalue $\zeta_h$, which is impossible. The data $f_h$ are real and lie
in $X_h$, and $\norm{f_h}_2=1$.

\paragraph{The grid mode.}
Let $H_0,\dots,H_3$ be the Hermite basis \eqref{eq:Hermite} and
$\mathsf m=H_1+H_3=r-3r^2+2r^3$. Let $\mathsf m_h(x)=h\,\mathsf m(x/h)$, with
$\mathsf m$ extended with period one. Thus $\mathsf m_h\in C^1[0,1]$ has zero
nodal values and unit nodal slopes. It is the constant-slope grid mode.
Since $|\mathsf m|\le1$ on $[0,1]$, also $|\mathsf m_h|\le h$.

\paragraph{Constants.}
The letter $C$ denotes a positive constant that may change from line to
line. It is independent of $N\ge512$, of $t$ in the indicated ranges,
and of the test functions $\vartheta$, $f$, $G$
where these appear. The notation $O(h)$ stands for a quantity bounded by
$Ch$ for all $N\ge512$. Limits are taken as $h\to0$, that is, as
$N\to\infty$.

\paragraph{The cell profile and the cell constant.}
At the stationary point, with $c=\cos\xi_*$ and $s=\sin\xi_*$, the slope
ratio \eqref{eq:slope} is purely imaginary, $\mathcal V_*=i\mathsf v_*$, with
\[
 \mathsf v_*=\frac{\varpi_*(312+108c)}{2-2c-26\varpi_*s}.
\]
Let
\[
 p(r)=H_0+i\mathsf v_*H_1+z_*(H_2+i\mathsf v_*H_3)
\]
be the cell profile. Since $p(0)=1$, we have $\int_0^1|p|^2>0$. Let
$\mathfrak g=|p|^2/\int_0^1|p|^2$, extended with period one, and put
\begin{equation}\label{eq:bdef}
 \mathfrak b=\frac12\int_0^1\mathfrak g(r)\mathsf m'(r)\,dr .
\end{equation}
We compute $\mathfrak b$ and show that it is positive. By
\eqref{eq:bdef},
$\mathfrak b=\tfrac12\int_0^1|p|^2\mathsf m'\big/\int_0^1|p|^2$, and
\[
 |p|^2=H_0^2+H_2^2+2cH_0H_2+\mathsf v_*^2(H_1^2+H_3^2+2cH_1H_3)+2\mathsf v_* s(H_1H_2-H_0H_3).
\]
We integrate $|p|^2$ against $1$ and against $\mathsf m'$. Against $1$, the
integrals are the entries of $\widehat M$ in \eqref{eq:element}, divided
by $420$. So
\[
 420\int_0^1|p|^2\,dr=312+108c+(8-6c)\mathsf v_*^2+52\mathsf v_*s .
\]
Against $\mathsf m'$, we put $\mathsf q=r(1-r)$. Then $\mathsf m'=1-6\mathsf q$ and
$\int_0^1\mathsf q^n\,dr=(n!)^2/(2n+1)!$. Since $H_0+H_2=1$ and $H_1+H_3=\mathsf m$,
\[
 H_0^2+H_2^2+2cH_0H_2=1-2(1-c)H_0H_2,\qquad
 H_1^2+H_3^2+2cH_1H_3=\mathsf m^2-2(1-c)H_1H_3 .
\]
Moreover $H_0H_2=\mathsf q^2(3+4\mathsf q)$, $H_1H_3=-\mathsf q^3$,
$H_1H_2-H_0H_3=\mathsf q^2(1+4\mathsf q)$,
$\int_0^1\mathsf m'=0$ and $\int_0^1\mathsf m^2\mathsf m'=0$. Direct
integration gives
\[
 \int_0^1H_0H_2\mathsf m'=-\frac4{105},\qquad
 \int_0^1H_1H_3\mathsf m'=\frac1{420},\qquad
 \int_0^1(H_1H_2-H_0H_3)\mathsf m'=-\frac2{105}.
\]
Hence
\[
 210\int_0^1|p|^2\mathsf m'\,dr=(1-c)(16-\mathsf v_*^2)-8\mathsf v_*s ,
\]
and therefore
\begin{equation}\label{eq:bclosed}
 \mathfrak b=\frac{(1-c)(16-\mathsf v_*^2)-8\mathsf v_*s}{312+108c+(8-6c)\mathsf v_*^2+52\mathsf v_*s}.
\end{equation}
By \eqref{eq:zgv}, $\varpi_*=w_*=R_c(c_*)s$. Since
$94-128c+34c^2=2(1-c)(47-17c)$, we obtain $\mathsf v_*=s\,\mathsf v(c)$ with
\[
 \mathsf v(c)=\frac{\mathsf N(c)}{\mathsf D(c)},\qquad \mathsf N(c)=14664-228c-1836c^2,\qquad
 \mathsf D(c)=-5031+7410c+1353c^2-582c^3 .
\]
We now show that
\begin{equation}\label{eq:bvalue}
 \mathfrak b>0 .
\end{equation}
Put $y=-\mathsf v(c)$. Since $s^2=(1-c)(1+c)$,
\[
 (1-c)(16-\mathsf v_*^2)-8\mathsf v_*s=(1-c)\bigl[16-(1+c)y\bigl((1-c)y-8\bigr)\bigr].
\]
The point $c=c_*$ lies in $[c_{\rm lo},c_{\rm hi}]\subset[0.528,0.529]$.
On $[0,1]$, $\mathsf N'(c)=-228-3672c<0$ and $\mathsf N(1)=12600$, so
$\mathsf N>0$. Moreover $\mathsf D'(c)=7410+2706c-1746c^2>0$ on $[0,1]$,
$\mathsf N(0.528)<14032$ and $\mathsf D(0.529)<-818$. Hence
$0<y<14032/818<17.16$. Thus $(1-c)y-8<0.472\cdot17.16-8<0.1$ and
$0<(1+c)y<1.529\cdot17.16<27$. Their product is less than $2.7$, also when
the first factor is negative. So the bracket exceeds $16-2.7>0$. Since
$\int_0^1|p|^2>0$ and $1-c>0$, this gives
$\mathfrak b>0$. Numerically, $\mathfrak b\approx0.0094219$.

\paragraph{The two-scale form of the eigenfunction.}
The next corollary follows from Theorem~\ref{thm:top} and from its proof,
through the results of Section~\ref{sec:leanlemmas}. It gives the
two-scale form of $|\varphi_h|^2$: up to $O(h)$, it is the product of
$\mathfrak d(x)$ and the one-periodic function $\mathfrak g(x/h)$.

\begin{corollary}[the eigenfunction]\label{lem:density}
Recall $\tilde f_h$, $f_h$, $\mathfrak d$, $\mathfrak g$ and $\mathfrak b$
from Table~\ref{tab:notation}. For $N\ge512$, the following
hold.
\begin{enumerate}
\item[(a)] $\sup_x\bigl||\varphi_h(x)|^2-\mathfrak d(x)\mathfrak g(x/h)\bigr|\le Ch$, $\norm{\varphi_h}_\infty\le C$ and $\norm{\varphi_h'}_\infty\le Ch^{-1}$.
\item[(b)] $(|\varphi_h|^2,G)=(\mathfrak d,G)+O(h)\norm G_{W^{1,\infty}}$ for $G\in W^{1,\infty}(0,1)$.
\item[(c)] $\tfrac12(|\varphi_h|^2,\mathsf m_h')=\mathfrak b+O(h)$.
\item[(d)] $\bigl|\int_0^1\varphi_h^2\,dx\bigr|\le16h^3$ and
$\bigl|\norm{\tilde f_h}_2^2-1\bigr|\le16h^3$. Moreover,
$\varphi_h,f_h\rightharpoonup0$ in $L^2$ as $h\to0$.
\end{enumerate}
\end{corollary}
\begin{proof}[Sketch]
Let $\varsigma=\varsigma_N$ and $\sigma=\sigma(\varsigma)$. By
Lemma~\ref{lem:det} and the choice of $\varphi_h$, the nodal data of
$\varphi_h$ are
\[
 (u_j,d_j)=\alpha_+z_+^j(1,\mathcal V_+)+\alpha_-z_-^j(1,\mathcal V_-)
 +\alpha_3z_3^j(1,\mathcal V_3)+\alpha_4z_4^{j-N}(1,\mathcal V_4),
\]
where $\alpha\ne0$ is the chosen kernel vector of $\D_N(\varsigma)$. The
key quantities are its four coefficients. We use $|y_N|\le8/N$
(Lemma~\ref{lem:box}) and the bounds and the analyticity of
Lemma~\ref{lem:domain}. They give $B_\pm=B_*+O(h)$,
$\mathcal V_\pm=\mathcal V_*+O(h)$ and $B_+-B_-=O(h)$.

Rows one and three of $\D_N\alpha=0$ are a $2\times2$ system for
$(\alpha_3,\alpha_4)$ with matrix
$\bigl(\begin{smallmatrix}1&\delta_4\\\delta_3&1\end{smallmatrix}\bigr)$,
where $|\delta_3|,|\delta_4|\le50^{-N}$ and $Z_\pm$ are bounded. Hence
$|\alpha_3|+|\alpha_4|\le C(|\alpha_+|+|\alpha_-|)$. Subtracting $B_3$
times row one from row two, we obtain
\[
 (B_+-B_3)\alpha_++(B_--B_3)\alpha_-=-\delta_4(B_4-B_3)\alpha_4 .
\]
The difference $B_--B_3$ is bounded away from zero. Indeed,
$|\varsigma_N|\le16/(5N)<1/20$ by Lemma~\ref{lem:box}, so
$|\mathcal R_+|\ge1$ by Lemma~\ref{lem:refl} and
$|\mathcal V_+-\mathcal V_4|\le2\cdot10^8$ by Lemma~\ref{lem:domain}.
Hence $|B_--B_3|\ge1/(2\cdot10^8)$. It follows that
$\alpha_+\ne0$ and $\alpha_-=-\alpha_+(1+O(h))$. Next,
\[
 Z_+/Z_-=e^{iN(\varsigma-\sigma)}=e^{i(2y_N-N(\varsigma+\sigma))}=1+O(h),
\]
by $|y_N|\le8/N$ and $|\varsigma+\sigma|\le(81/10)|\varsigma|^2$. Rows one
and three now give $\alpha_3=O(h|\alpha_+|)$ and
$\alpha_4=O(h|\alpha_+|)$.

Uniformly in $0\le j\le N$ we have $\varsigma j=\pi x_j+O(h)$ and
$\sigma j=-\pi x_j+O(h)$, so $z_\pm^j=z_*^je^{\pm i\pi x_j}(1+O(h))$.
With $c_h=2i\alpha_+$, the Hermite representation on $I_j$ gives the cell
form
\begin{equation}\label{eq:cell}
 \varphi_h(x_j+hr)=c_hz_*^j\sin(\pi x_j)\,p(r)+O(h|c_h|),
\end{equation}
uniformly in $j$ and $r$. Combining \eqref{eq:cell} with
$h\sum_{j=0}^{N-1}\sin^2\pi x_j=\tfrac12$ and the normalization, we
obtain $|c_h|^2\int_0^1|p|^2=2+O(h)$.

Part (a) follows from \eqref{eq:cell}. The slope bound follows from the
inverse inequality for cubics on each cell. For (b) and (c) we use an
averaging inequality.
For a bounded one-periodic $\mathsf c$ with mean $\bar{\mathsf c}$ and $G\in W^{1,1}(0,1)$,
\begin{equation}\label{eq:average}
 \Bigl|\int_0^1G(x)\mathsf c(x/h)\,dx-\bar{\mathsf c}\int_0^1G\Bigr|\le h\norm{\mathsf c-\bar{\mathsf c}}_\infty\norm{G'}_1 .
\end{equation}
To see this, we subtract $G(x_j)$ on each $I_j$.
Part (b) follows from (a) and \eqref{eq:average} with
$\mathsf c=\mathfrak g$, $\bar{\mathsf c}=1$. Part (c) follows from (a) and
\eqref{eq:average} with $\mathsf c=\mathfrak g\mathsf m'$ and
$G=\mathfrak d$, since $\int_0^1\mathfrak d=1$.

Part (d) does not use the cell form. We test \eqref{eq:Sscheme} with
$\chi=\overline{\varphi_h}$. Since $\varphi_h'(1)=0$,
\[
 \zeta_h\int_0^1\varphi_h^2\,dx=\int_0^1\varphi_h''\varphi_h'\,dx
 =-\tfrac12\varphi_h'(0)^2 .
\]
By \eqref{eq:flux}, $|\varphi_h'(0)|^2=2|\Rea\zeta_h|$. Hence, by
\eqref{eq:bounds},
$\bigl|\int_0^1\varphi_h^2\bigr|=|\Rea\zeta_h|/|\zeta_h|\le1600/(100h^{-3})=16h^3$.
The second bound follows from
$\norm{\tilde f_h}_2^2=1+\Rea\int_0^1\varphi_h^2$. The weak convergence
$\varphi_h\rightharpoonup0$ is Corollary~\ref{cor:principleS}(b) below,
whose sketch does not use the present corollary. Since
$(\Rea\varphi_h,\psi)=\Rea(\varphi_h,\psi)$ for real $\psi$, also
$\tilde f_h=\sqrt2\Rea\varphi_h\rightharpoonup0$. Since
$\norm{\tilde f_h}_2\to1$, finally
$f_h=\tilde f_h/\norm{\tilde f_h}_2\rightharpoonup0$.
\end{proof}

\subsection{The inverse of the discrete operator}\label{sec:inverse}
We now turn to the inverse of $\mathring A_h$.
Lemmas~\ref{lem:inverse} and~\ref{lem:twoscale} give it exactly and in a
two-scale form. They do not involve the eigenpairs.
Lemmas~\ref{lem:slow} and~\ref{lem:fast} apply the inverse to
$|\varphi_h|^2$ and $\varphi_h^2$.

\begin{lemma}[exact inverse]\label{lem:inverse}
Let $\vartheta\in L^2(0,1)$ be real. Let $\psi_0\in H^3$ solve
$\psi_0'''=\vartheta$, $\psi_0(0)=\psi_0'(0)=\psi_0(1)=0$, and put $g=\psi_0'$.
Let $\mathcal S_h\vartheta\in X_h$ have nodal values $\psi_0(x_j)$ and nodal slopes
\[
 p_j=g(x_j)-g(1)-\int_{x_j}^1\mathsf m_h\,\vartheta\,dx\qquad(0\le j\le N).
\]
Then $(v'',(\mathcal S_h\vartheta)')=(v,\vartheta)$ for every $v\in X_h$.
In particular $\mathring A_h$ is invertible.
\end{lemma}
\begin{proof}[Sketch]
Since $\psi_0(0)=\psi_0(1)=0$ and $p_N=0$, indeed $\psi:=\mathcal S_h\vartheta\in X_h$.
Integrating by parts twice and using $v(0)=v(1)=0$, $v'(1)=0$ and
$g(0)=0$ for the boundary terms, we obtain
\begin{equation}\label{eq:ibp}
 (v,\vartheta)=(v,\psi_0''')=(v'',g).
\end{equation}
It remains to show $(v'',\phi')=0$ for $\phi=\psi-\psi_0$.
On $I_j$ write $v''=m+n(r-\tfrac12)$. Since $\phi(x_j)=\phi(x_{j+1})=0$,
\[
 \int_{I_j}\phi'\,dx=0,\qquad
 \int_{I_j}(r-\tfrac12)\phi'\,dx=-\frac1h\int_{I_j}\phi\,dx .
\]
So it suffices that $\int_{I_j}\psi=\int_{I_j}\psi_0$ for each $j$. For the cubic $\psi$,
\[
 \int_{I_j}\psi\,dx=\frac h2\bigl(\psi_0(x_j)+\psi_0(x_{j+1})\bigr)-\frac{h^2}{12}(p_{j+1}-p_j).
\]
The Peano form of the trapezoidal rule gives
\[
 \int_{I_j}\psi_0\,dx=\frac h2\bigl(\psi_0(x_j)+\psi_0(x_{j+1})\bigr)
 -\frac{h^2}2\int_{I_j}r(1-r)g'\,dx .
\]
We need $p_{j+1}-p_j=\int_{I_j}6r(1-r)g'\,dx$. Since $6r(1-r)=1-\mathsf m'(r)$
and $\mathsf m_h'=\mathsf m'(r)$, integrating by parts and using that
$\mathsf m_h$ vanishes at the nodes, we obtain
\[
 \int_{I_j}6r(1-r)g'\,dx
 =g(x_{j+1})-g(x_j)-\int_{I_j}\mathsf m_h'g'\,dx
 =g(x_{j+1})-g(x_j)+\int_{I_j}\mathsf m_h\,\vartheta\,dx .
\]
This is $p_{j+1}-p_j$ by definition.

If $\mathring A_hv=0$, then $(v,\vartheta)=(\mathring A_hv,\mathcal S_h\vartheta)=0$ for all $\vartheta$, so $v=0$.
\end{proof}

\begin{lemma}[two-scale form]\label{lem:twoscale}
Let $\vartheta$, $g$ and $\mathcal S_h\vartheta$ be as in
Lemma~\ref{lem:inverse}. Then
\begin{align}
 \norm g_{W^{1,\infty}}&\le C\norm\vartheta_2,\label{eq:gbound}\\
 \norm{(\mathcal S_h\vartheta)'-g+g(1)\mathsf m_h'}_{L^\infty}&\le Ch\norm\vartheta_2.\label{eq:twoscale}
\end{align}
Consequently, for $f\in L^2(0,1)$,
\begin{equation}\label{eq:inverseload}
 (\mathring A_h^{-1}\mathcal G_h[f],\vartheta)=\tfrac12(f,g)-g(1)\,\tfrac12(f,\mathsf m_h')+R,
 \qquad |R|\le Ch\norm f_1\norm\vartheta_2,
\end{equation}
and $\norm{\mathring A_h^{-1}\mathcal G_h[f]}_2\le C\norm f_1$.
\end{lemma}
\begin{proof}[Sketch]
By variation of constants, $\psi_0=\int_0^x\tfrac12(x-y)^2\vartheta(y)dy+\tfrac12\alpha_0 x^2$,
with $\alpha_0=-\int_0^1(1-y)^2\vartheta\,dy$. This proves \eqref{eq:gbound}.
On $I_j$,
\[
 (\mathcal S_h\vartheta)'=\frac{\psi_0(x_{j+1})-\psi_0(x_j)}hH_2'(r)+p_jH_1'(r)+p_{j+1}H_3'(r).
\]
The difference quotient equals $g(x)+O(h\norm{g'}_\infty)$.
By $|\mathsf m_h|\le h$, also $p_j,p_{j+1}=g(x)-g(1)+O(h\norm\vartheta_2)$.
Since $H_1'+H_3'=\mathsf m'$ and $H_2'+\mathsf m'=1$, we obtain \eqref{eq:twoscale}.
For \eqref{eq:inverseload} let $z=\mathring A_h^{-1}\mathcal G_h[f]$. By Lemma~\ref{lem:inverse} with $v=z$,
\[
 (z,\vartheta)=(z'',(\mathcal S_h\vartheta)')=(\mathring A_hz,\mathcal S_h\vartheta)=\tfrac12(f,(\mathcal S_h\vartheta)').
\]
Inserting \eqref{eq:twoscale}, we obtain \eqref{eq:inverseload}. The norm bound follows from \eqref{eq:gbound} and $|\mathsf m_h'|\le1$.
\end{proof}
Both lemmas are linear in $\vartheta$. They extend to complex $\vartheta$,
with $\mathcal S_h\vartheta$ in the complexification of $X_h$ and with the
same bilinear pairings.

We now turn to the eigenfunctions.
Put
\begin{equation}\label{eq:zinf}
 z_\infty=-\frac{\mathfrak d}{8\pi^2}+\mathfrak b\mathfrak p .
\end{equation}

\begin{lemma}[slow part]\label{lem:slow}
Let $z_h=\mathring A_h^{-1}\mathcal G_h[|\varphi_h|^2]$. Then $\norm{z_h-z_\infty}_2\le Ch$.
\end{lemma}
\begin{proof}[Sketch]
Let $\vartheta\in L^2(0,1)$ be real. Let $\psi_0$ and $g$ be as in
Lemma~\ref{lem:inverse}. By \eqref{eq:inverseload},
$\norm{|\varphi_h|^2}_1=1$, Corollary~\ref{lem:density}(b,c) and
\eqref{eq:gbound},
\[
 (z_h,\vartheta)=\tfrac12(\mathfrak d,g)-\mathfrak b\,g(1)+O(h)\norm\vartheta_2 .
\]
It remains to identify $(z_\infty,\vartheta)$ with the first two terms.
Since $\mathfrak d$ and $\mathfrak d'$ vanish at both ends, $(\mathfrak d,\vartheta)=(\mathfrak d,\psi_0''')=(\mathfrak d'',g)$.
Further $\mathfrak d''=4\pi^2(1-\mathfrak d)$ and $\int_0^1g=0$. Thus
$-(8\pi^2)^{-1}(\mathfrak d,\vartheta)=\tfrac12(\mathfrak d,g)$.
Also $(\mathfrak p,\vartheta)=(\mathfrak p,\psi_0''')=-[\mathfrak p'g]_0^1+(\mathfrak p'',g)=-g(1)$,
since $g(0)=0$ and $\int_0^1g=0$. Hence
$|(z_h-z_\infty,\vartheta)|\le Ch\norm\vartheta_2$. Both $z_h$ and
$z_\infty$ are real, so the supremum over real $\vartheta$ gives the
claim.
\end{proof}

\begin{lemma}[oscillating part]\label{lem:fast}
Let
\[
 \mathcal O_h(t)=\int_0^t\mathring T_h(t-s)\mathcal G_h[\varphi_h^2]e^{2\zeta_hs}\,ds .
\]
Then, for $N\ge512$, $\norm{(2\zeta_h-\mathring A_h)^{-1}}_2\le h^3/15$ and
$\sup_{t\ge0}\norm{\mathcal O_h(t)}_2\le Ch^2$.
\end{lemma}
\begin{proof}[Sketch]
The idea is that $2\zeta_h$ lies outside the band of frequencies of the
grid waves. The key quantity is the bound $3/10$ for this band.

First, both roots of $P(\xi,\cdot)$ in \eqref{eq:dispersion} lie in $(-3/10,3/10)$
for every real $\xi$. Indeed $P(\xi,0)=-d\le0$, $a>0$, and
$P(\xi,\pm3/10)\ge U-\tfrac3{10}|b|$ with
\[
 U(c)=\tfrac9{100}a-d=\tfrac{799}{20}-\tfrac{517}{25}c-\tfrac{37}{100}c^2 .
\]
The function $U$ is concave and positive at $c=\pm1$, so $U>0$ on $[-1,1]$. It remains to show
$Q(c)=U(c)^2-\tfrac9{100}(1-c^2)(104+c)^2>0$ there. Direct expansion gives
\[
 10^4\bigl(Q(c)-119\bigr)=2269c^4+340232c^3+13714494c^2-16710520c+5035625 .
\]
The Sturm sequence of this polynomial has sign patterns $(+,-,-,+,+)$ at $c=-1$ and
$(+,+,+,-,+)$ at $c=1$. Both have two sign changes, so there is no zero in $[-1,1]$.
The value at $c=1$ is positive. Hence $Q>119$ on $[-1,1]$.

Second, let $v$ be in the complexification of $X_h$ with $\norm v_2=1$,
and put $A_0=\tfrac12|v'(0)|^2$, $I=\Ima\int_0^1v''\overline{v'}$. Let
$\tilde v$ be the piecewise cubic on $\mathbb R$ with the nodal data of
$v$ on $0\le j\le N$ and zero data elsewhere. Then $\tilde v=v$ on
$[0,1]$, and on $[-h,0]$ it equals $hv'(0)H_3$. So
$\norm{\tilde v}_{L^2(\mathbb R)}^2=1+\frac{h^3}{105}|v'(0)|^2$, and the
extra contribution to $\int\tilde v''\overline{\tilde v'}$ is real. Let
$\hat z(\xi)=\sum_j(u_j,d_j)^Te^{-ij\xi}$ be the Fourier series of the
nodal data of $\tilde v$. By Parseval's identity and the assembly rule of
Section~\ref{sec:disp}, $\norm{\tilde v}_{L^2(\mathbb R)}^2$ and
$\int_{\mathbb R}\tilde v''\overline{\tilde v'}$ are the means over
$\xi\in(-\pi,\pi)$ of the forms $(h/420)\,\hat z^*M_0\hat z$ and
$h^{-2}\,\hat z^*K_0\hat z$, with the matrices of \eqref{eq:bloch}. The
matrix $M_0$ is Hermitian and positive definite, and $K_0$ is
skew-Hermitian. By \eqref{eq:dispersion}, the eigenvalues of $-iK_0$
relative to $M_0$ are the roots of $P(\xi,\cdot)$. Since they lie in
$(-3/10,3/10)$, we have
$|\hat z^*K_0\hat z|\le\tfrac3{10}\,\hat z^*M_0\hat z$. Hence
\[
 |I|\le\frac{420\cdot3/10}{h^3}\norm{\tilde v}^2_{L^2(\mathbb R)}
 =\frac{126}{h^3}+\frac{12}5A_0 .
\]

Third, write $2\zeta_h=2\mu_h+i\omega$. By \eqref{eq:im},
$\Omega_*>420\cdot0.276$ (Lemma~\ref{lem:cert}, (C1)) and $|\beta|\le15$
(Lemma~\ref{lem:refl}), we have $|\omega|\ge230h^{-3}$ for $N\ge512$. By
\eqref{eq:flux},
\[
 \int_0^1\bigl((2\zeta_h-\mathring A_h)v\bigr)\overline v
 =2\mu_h+A_0+i(\omega-I).
\]
If $A_0\le30h^{-3}$, then $|I|\le198h^{-3}$, and the imaginary part has
absolute value $|\omega-I|\ge|\omega|-|I|\ge32h^{-3}$. Otherwise the real
part is at least $30h^{-3}-3200\ge15h^{-3}$, by \eqref{eq:bounds}. So
$\norm{(2\zeta_h-\mathring A_h)v}_2\ge15h^{-3}\norm v_2$. This is the first
claim.

Fourth, on the complexification of $X_h$ we have the inverse inequality
\begin{equation}\label{eq:inverse2}
 \norm{\chi'}_2\le C_Ih^{-1}\norm\chi_2 ,
\end{equation}
which follows as \eqref{eq:inverse} by expanding cubics in Legendre
polynomials on each cell. With Corollary~\ref{lem:density}(a) it gives
$\norm{\mathcal G_h[\varphi_h^2]}_2\le Ch^{-1}$. Thus
$q_h=(2\zeta_h-\mathring A_h)^{-1}\mathcal G_h[\varphi_h^2]$ has
$\norm{q_h}_2\le Ch^2$. Finally
$\mathcal O_h(t)=e^{2\zeta_ht}q_h-\mathring T_h(t)q_h$, and both terms are
bounded by $\norm{q_h}_2$.
\end{proof}

\subsection{Galerkin limits}\label{sec:galerkin}
The continuous objects $A^{-1}$ and $T$ are the limits of the W-scheme.
The semigroups of the S-scheme have the same strong limit. In this
subsection $L=1$ and the meshes are uniform. Recall the form of the
W-scheme from~\cite[(2.2) and (3.1)]{Main}: for $w,\psi\in X_h$,
\begin{equation}\label{eq:Wform}
 (A_hw,\psi)_H=\int_0^1w''\,\overline{\bigl((1+x)\psi\bigr)'}\,dx ,
\end{equation}
and $(Aw,\psi)_H$ is given by the same integral for $w\in D(A)$ and
$\psi\in H_0^1(0,1)$.

\begin{proposition}[$A^{-1}$ and $T$ as Galerkin limits]\label{prop:galerkin}
Let $L=1$ and let the meshes be uniform.
\begin{enumerate}
\item[(a)] $\norm{A_h^{-1}\Pi_h-A^{-1}}_{L^2\to L^2}\le Ch$, with $C$
independent of $h$.
\item[(b)] For every $f\in L^2(0,1)$ and $T_0>0$,
$\sup_{0\le t\le T_0}\norm{T_h(t)\Pi_hf-T(t)f}_2\to0$ as $h\to0$.
\item[(c)] For every $f\in L^2(0,1)$, $\varrho=A^{-1}f$ is the unique
$\varrho\in H^3(0,1)$ with $-\varrho'''=f$ and
$\varrho(0)=\varrho(1)=\varrho'(1)=0$. Moreover $T$ is a
$C_0$-semigroup on $L^2(0,1)$, and
$A^{-1}\bigl(T(t)f-f\bigr)=\int_0^tT(s)f\,ds$ for $f\in L^2(0,1)$ and
$t\ge0$.
\end{enumerate}
\end{proposition}
\begin{proof}[Sketch]
\emph{(a)} Let $f\in L^2(0,1)$, $w=A^{-1}f$ and $w_h=A_h^{-1}\Pi_hf$. For
$\psi\in X_h$, $(A_hw_h,\psi)_H=(\Pi_hf,\psi)_H=(f,\psi)_H=(Aw,\psi)_H$.
By \eqref{eq:Wform},
\[
 \int_0^1(w_h-w)''\,\overline{\bigl((1+x)\psi\bigr)'}\,dx=0\qquad(\psi\in X_h).
\]
Thus $w_h$ is the Ritz projection of $w$. Let $I_hw\in X_h$ be the cubic
Hermite interpolant of $w$ and put $\phi=w_h-I_hw$. By the weighted
energy identity of Section~\ref{sec:setting} and the orthogonality with
$\psi=\phi$,
\[
 \tfrac32\norm{\phi'}_2^2\le\Rea(-A_h\phi,\phi)_H
 =\Rea\int_0^1(I_hw-w)''\,\overline{\bigl((1+x)\phi\bigr)'}\,dx
 \le3\norm{(w-I_hw)''}_2\norm{\phi'}_2 ,
\]
because $\norm{((1+x)\phi)'}_2\le\norm\phi_2+2\norm{\phi'}_2\le3\norm{\phi'}_2$.
On each cell $I_j$, $w-I_hw$ and its derivative vanish at both ends.
Two integrations by parts show that $(w-I_hw)''$ is orthogonal to the
linear functions on $I_j$. So $(I_hw)''$ is the $L^2(I_j)$ projection of
$w''$ onto the linear functions, and
\[
 \norm{(w-I_hw)''}_{L^2(I_j)}\le\norm{w''-w''(x_j)}_{L^2(I_j)}
 \le h\norm{w'''}_{L^2(I_j)}=h\norm f_{L^2(I_j)} .
\]
Hence $\norm{\phi'}_2\le2h\norm f_2$. Since $w-I_hw$ and its derivative
vanish at $x_j$, also
$\norm{w-I_hw}_{L^2(I_j)}\le h^2\norm{(w-I_hw)''}_{L^2(I_j)}$. Together
with $\norm\phi_2\le\norm{\phi'}_2$, this gives
$\norm{w_h-w}_2\le\norm\phi_2+\norm{I_hw-w}_2\le3h\norm f_2$.

\emph{(b)} Put $J_h=A_h^{-1}\Pi_h$. Then $A_hJ_h=\Pi_h$ and
$\Pi_hJ_h=J_h$, so $\norm{J_h-\Pi_hA^{-1}}=\norm{\Pi_h(J_h-A^{-1})}\le Ch$
by (a). We also use $A^{-1}T(s)=T(s)A^{-1}$ and
$\frac d{ds}T(s)A^{-1}x=T(s)x$ for $x\in L^2(0,1)$. Hence, for
$0\le s\le t$,
\begin{equation}\label{eq:TKidentity}
 \frac d{ds}\Bigl[T_h(t-s)J_hT(s)A^{-1}x\Bigr]
 =T_h(t-s)\bigl(J_h-\Pi_hA^{-1}\bigr)T(s)x .
\end{equation}
Both semigroups are contractions on $H$. Integrating
\eqref{eq:TKidentity} over $[0,t]$, we obtain
\[
 \norm{J_hT(t)A^{-1}x-T_h(t)J_hA^{-1}x}_2\le Cth\norm x_2 .
\]
In the first term we replace $J_h$ by $A^{-1}$, and in the second by
$\Pi_hA^{-1}$. Each replacement costs at most $Ch\norm x_2$. Thus
\begin{equation}\label{eq:TKrate}
 \norm{\bigl(T_h(t)\Pi_h-T(t)\bigr)A^{-2}x}_2\le C(1+t)h\norm x_2 .
\end{equation}
The range of $A^{-2}$ is $D(A^2)$, which is dense in $L^2(0,1)$, and the
operators $T_h(t)\Pi_h-T(t)$ are bounded uniformly in $h$ and $t$. This
gives (b).

\emph{(c)} The first statement is the definition of $A$ in
Section~\ref{sec:setting}. The integral identity holds for every
$C_0$-semigroup and its generator.
\end{proof}
The sketch uses the Ritz projection and one comparison function on each
cell. It does not use the mesh-uniform regularity of the W-scheme. With
that regularity, the convergence in (b) holds in operator norm at
positive times:
$\norm{T_h(t)\Pi_h-T(t)}_{H\to H}\le Ce^{-\nu t}ht_\wedge^{-3/2}$
\cite[Proposition~4.6]{Main}.

The Lean proof goes from the scheme to the continuous objects. It does
not presuppose $A^{-1}$ and $T$. It defines them as the limits in (a) and
(b), proves that these limits exist, and proves (c) for the limits. Part
(c) identifies them with the objects of Section~\ref{sec:setting}; in
particular, the generator of the limit semigroup is $A$.

\paragraph{Compactness.}
Recall from Section~\ref{sec:setting} that $A$ has compact resolvent. In
particular, $A^{-1}$ is a compact injective operator on $L^2(0,1)$. On
$L^2(0,1)$, the semigroup $T$ is also compact and has Crandall--Pazy
regularity of index $3/2$. Indeed, by~\cite[Theorem~2.4(ii)]{Main} and the
equivalence of the norms of $H$ and $L^2$,
$\norm{AT(t)}_{L^2\to L^2}\le Ct_\wedge^{-3/2}$ for $t>0$. Hence $T(t)$
maps $L^2(0,1)$ into $D(A)$, and $T(t)=A^{-1}AT(t)$ is compact for every
$t>0$. These two facts also follow from the scheme. By~\cite[(3.15)]{Main}
with $\alpha=1$, $\sup_N\norm{A_hT_h(t)}<\infty$ for each $t>0$. Now
$T_h(t)\Pi_hf=A_h^{-1}\Pi_h\,g_h$ with $g_h=A_hT_h(t)\Pi_hf$, and the
family $(g_h)$ is bounded in $L^2$. By
Proposition~\ref{prop:galerkin}, every weak limit point $g$ of $(g_h)$
satisfies $A^{-1}g=T(t)f$. So $T(t)f\in D(A)$, and
$\norm{AT(t)f}_2\le C_t\norm f_2$. In the Lean proofs this is the only
place where the mesh-uniform regularity of the W-scheme is used; see
Remark~\ref{rem:leanscope}.

\begin{lemma}[strong convergence]\label{lem:TK}
For each fixed $f\in L^2$ and $T_0>0$,
$\sup_{0\le t\le T_0}\norm{\mathring T_h(t)\mathring\Pi_hf-T(t)f}_2\to0$.
\end{lemma}
\begin{proof}[Sketch]
The proof of Proposition~\ref{prop:galerkin}(b) applies, with
$\mathring J_h=\mathring A_h^{-1}\mathring\Pi_h$, $\mathring T_h$ and
$\mathring\Pi_h$ in place of $J_h$, $T_h$ and $\Pi_h$. It uses
$\mathring A_h\mathring J_h=\mathring\Pi_h$,
$\mathring\Pi_h\mathring J_h=\mathring J_h$, the contraction property of
$\mathring T_h$ and $T$ on $L^2$, and
$\norm{\mathring J_h-A^{-1}}_{L^2\to L^2}\le Ch$. The last bound is
Corollary~\ref{cor:principleS}(a) below. Its sketch uses only
Lemmas~\ref{lem:inverse} and~\ref{lem:twoscale}. We obtain
\[
 \norm{\bigl(\mathring T_h(t)\mathring\Pi_h-T(t)\bigr)A^{-2}x}_2\le C(1+t)h\norm x_2
 \qquad(x\in L^2(0,1)),
\]
and the lemma follows, since $D(A^2)$ is dense and both families are
contractions on $L^2$.
\end{proof}
The rate holds for fixed $f\in D(A^2)$ only. The lemma gives strong
convergence, not convergence in operator norm. By
Corollary~\ref{cor:principleS}(c), convergence in operator norm is false.

\subsection{The linear error gap}\label{sec:lineargap}
Behind the failure of the S-scheme is a general principle. We state it
for abstract families and then apply it to the S-scheme.

\begin{proposition}[uniform high-frequency damping is necessary]\label{prop:necessary}
Let $\mathcal X$ be a complex Hilbert space. Let $X_n\subset\mathcal X$ be
closed subspaces, with orthogonal projections $\Pi_n$, and let $A_n$ be
bounded operators on $X_n$. Let $(\zeta_n,x_n)$ be eigenpairs of $A_n$
with $\norm{x_n}=1$.
\begin{enumerate}
\item[(a)] Suppose that each $A_n$ is invertible and that
$\norm{A_n^{-1}\Pi_n-K}\to0$ for an injective bounded operator $K$ on
$\mathcal X$. If $|\zeta_n|\to\infty$, then $x_n\rightharpoonup0$.
\item[(b)] Let $S$ be an eventually compact $C_0$-semigroup on $\mathcal X$,
that is, $S(t_1)$ is compact for some $t_1>0$. Suppose that
$x_n\rightharpoonup0$. Then, for every $t\ge t_1$,
\[
 \liminf_{n\to\infty}\norm{e^{tA_n}\Pi_n-S(t)}\ge\liminf_{n\to\infty}e^{t\Rea\zeta_n}.
\]
If, in addition, $\norm{e^{t_0A_n}\Pi_n-S(t_0)}\to0$ for some $t_0>0$,
then $\Rea\zeta_n\to-\infty$.
\end{enumerate}
\end{proposition}
\begin{proof}[Sketch]
Write $\langle\cdot,\cdot\rangle$ for the inner product of $\mathcal X$.

\emph{(a)} Since $A_n$ is invertible, $\zeta_n\ne0$. We have
$\Pi_nx_n=x_n$ and $A_n^{-1}x_n=\zeta_n^{-1}x_n$. Hence
\[
 \norm{Kx_n}\le\norm{(K-A_n^{-1}\Pi_n)x_n}+\norm{A_n^{-1}x_n}
 \le\norm{K-A_n^{-1}\Pi_n}+|\zeta_n|^{-1}\to0 .
\]
For $y\in\mathcal X$ it follows that
$\langle x_n,K^*y\rangle=\langle Kx_n,y\rangle\to0$. Since $K$ is
injective, the range of $K^*$ is dense in $\mathcal X$. Since
$\norm{x_n}=1$, we obtain $x_n\rightharpoonup0$.

\emph{(b)} Let $t\ge t_1$. Then $S(t)=S(t-t_1)S(t_1)$ is compact. We
have $x_n=\Pi_nx_n$ and $e^{tA_n}x_n=e^{t\zeta_n}x_n$, so
\begin{equation}\label{eq:necessary}
 e^{t\Rea\zeta_n}=\norm{e^{tA_n}x_n}\le\norm{e^{tA_n}\Pi_n-S(t)}+\norm{S(t)x_n}.
\end{equation}
Since $S(t)$ is compact and $x_n\rightharpoonup0$, we have
$\norm{S(t)x_n}\to0$. Taking the lower limit in \eqref{eq:necessary}, we
obtain the inequality.

Now suppose that $\norm{P_n-S(t_0)}\to0$, where $P_n=e^{t_0A_n}\Pi_n$.
The range of $P_n$ lies in $X_n$, so $P_n^k=e^{kt_0A_n}\Pi_n$ for every
integer $k\ge1$. The norms $\norm{P_n}$ are bounded, because $P_n$
converges. Hence
$\norm{e^{kt_0A_n}\Pi_n-S(kt_0)}=\norm{P_n^k-S(t_0)^k}\to0$. We choose
$k$ with $kt_0\ge t_1$. Then \eqref{eq:necessary} with $t=kt_0$ gives
$e^{kt_0\Rea\zeta_n}\to0$, that is, $\Rea\zeta_n\to-\infty$.
\end{proof}
Suppose that the inverses converge as in (a), and that the semigroups
converge as in (b) at some positive time. Then $|\zeta_n|\to\infty$
implies $\Rea\zeta_n\to-\infty$. Eigenvalues of large modulus must be
strongly damped.

With Lemmas~\ref{lem:inverse} and~\ref{lem:twoscale}, the principle
applies to the S-scheme.

\begin{corollary}[the principle for the S-scheme]\label{cor:principleS}
Recall $(\zeta_h,\varphi_h)$ and $\mathcal C_S$ from
Table~\ref{tab:notation}. Then the
following hold.
\begin{enumerate}
\item[(a)] $\norm{\mathring A_h^{-1}\mathring\Pi_h-A^{-1}}_{L^2\to L^2}\le Ch$,
with $C$ independent of $h$.
\item[(b)] $\varphi_h\rightharpoonup0$ in $L^2$.
\item[(c)] For every $t>0$,
$\liminf_{h\to0}\norm{\mathring T_h(t)\mathring\Pi_h-T(t)}_{L^2\to L^2}\ge e^{-\mathcal C_St}$.
\end{enumerate}
\end{corollary}
\begin{proof}[Sketch]
\emph{(a)} Let $f\in L^2(0,1)$. Put
$\varrho_h=\mathring A_h^{-1}\mathring\Pi_hf$ and $\varrho=A^{-1}f$. Let
$\vartheta\in L^2(0,1)$ be complex, and let $\psi_0$ and $g$ be as in
Lemma~\ref{lem:inverse}. In bilinear form, \eqref{eq:Sscheme} reads
$(\mathring A_hv,\chi)=(v'',\chi')$, and $(\mathring\Pi_hf,\chi)=(f,\chi)$
for $\chi$ in the complexification of $X_h$. By Lemma~\ref{lem:inverse},
\[
 (\varrho_h,\vartheta)=(\varrho_h'',(\mathcal S_h\vartheta)')
 =(\mathring A_h\varrho_h,\mathcal S_h\vartheta)
 =(\mathring\Pi_hf,\mathcal S_h\vartheta)=(f,\mathcal S_h\vartheta).
\]
Integrating by parts three times, we obtain
$(\varrho,\vartheta)=(\varrho,\psi_0''')=-(\varrho''',\psi_0)=(f,\psi_0)$.
The boundary terms vanish, since $\varrho(0)=\varrho(1)=\varrho'(1)=0$ and
$\psi_0(0)=\psi_0'(0)=\psi_0(1)=0$. Hence
$(\varrho_h-\varrho,\vartheta)=(f,\mathcal S_h\vartheta-\psi_0)$.

We integrate \eqref{eq:twoscale} from the left end. There
$\mathcal S_h\vartheta$, $\psi_0$ and $\mathsf m_h$ vanish. Moreover
$|\mathsf m_h|\le h$, and $|g(1)|\le C\norm\vartheta_2$ by \eqref{eq:gbound}.
Thus $\norm{\mathcal S_h\vartheta-\psi_0}_\infty\le Ch\norm\vartheta_2$, and
$|(\varrho_h-\varrho,\vartheta)|\le Ch\norm f_2\norm\vartheta_2$. We take
$\vartheta=\overline{\varrho_h-\varrho}$ and obtain
$\norm{\varrho_h-\varrho}_2\le Ch\norm f_2$.

\emph{(b)} By Lemma~\ref{lem:inverse}, $\mathring A_h$ is invertible. We
apply Proposition~\ref{prop:necessary}(a), with the index $n$ replaced by
$N\ge512$, $\mathcal X=L^2(0,1)$, $X_n$ the complexification of $X_h$,
$\Pi_n=\mathring\Pi_h$, $A_n=\mathring A_h$ and $K=A^{-1}$, which is
injective. The assumptions hold by (a), by $\norm{\varphi_h}_2=1$, and by
$|\zeta_h|\ge100h^{-3}$ from \eqref{eq:bounds}. Hence
$\varphi_h\rightharpoonup0$.

\emph{(c)} By Section~\ref{sec:galerkin}, $T(t)$ is compact on $L^2$ for
every $t>0$. In particular, $T$ is eventually compact. By (b),
Proposition~\ref{prop:necessary}(b) applies with $S=T$ and with any
$t_1\in(0,t]$. Since $\Rea\zeta_h\to-\mathcal C_S$ by \eqref{eq:re}, it
gives
$\liminf_{h\to0}\norm{\mathring T_h(t)\mathring\Pi_h-T(t)}\ge e^{-\mathcal C_St}$.
\end{proof}
Thus the inverses of the generators of the S-scheme converge in operator norm, at the rate
$h$. Its semigroups do not converge in operator norm at any positive
time. The obstruction comes from the eigenpairs $(\zeta_h,\varphi_h)$.
Their frequency tends to infinity. Their damping stays bounded.

On the data $f_h$ the linear error of the S-scheme does not tend to
zero.

\begin{corollary}[linear error gap]\label{cor:lineargap}
Recall $f_h$ and $\mathcal C_S$ from Table~\ref{tab:notation}.
For every $0<t_0<T_0$,
\[
 \sup_{t_0\le t\le T_0}\Bigl|\norm{(T(t)-\mathring T_h(t))f_h}_2-e^{-\mathcal C_St}\Bigr|\to0
 \qquad\text{as }h\to0 .
\]
\end{corollary}
\begin{proof}[Sketch]
We have $\mathring T_h(t)\tilde f_h=\sqrt2\,\Rea w$ with
$w=e^{\zeta_ht}\varphi_h$. Put $m_h=\int_0^1\varphi_h^2$, so
$|m_h|\le16h^3$ by Corollary~\ref{lem:density}(d). From
$2(\Rea w)^2=|w|^2+\Rea w^2$ and
$\norm{\tilde f_h}_2^2=1+\Rea m_h$,
\[
 \norm{\mathring T_h(t)f_h}_2^2
 =\frac{e^{2\mu_ht}+\Rea(e^{2\zeta_ht}m_h)}{1+\Rea m_h},\qquad
 \Bigl|\norm{\mathring T_h(t)f_h}_2^2-e^{2\mu_ht}\Bigr|\le64h^3 ,
\]
uniformly in $t\ge0$. Since $f_h\rightharpoonup0$ and
$T(t_0)$ is compact, $\norm{T(t)f_h}_2\le\norm{T(t_0)f_h}_2\to0$ for
$t\ge t_0$. The claim follows from the triangle inequality and
$\mu_h\to-\mathcal C_S$.
\end{proof}
The condition $t_0>0$ cannot be dropped: at $t=0$ the difference
vanishes.

For every fixed datum the linear error of the S-scheme converges to zero.
For each prescribed rate, there is a fixed datum for which that rate fails.

\begin{corollary}[fixed data]\label{cor:fixed}
Recall $\mathcal C_S$ and $\mathring L_h(t)$ from Table~\ref{tab:notation}. Fix $t>0$.
\begin{enumerate}
\item[(a)] $\liminf_{h\to0}\norm{\mathring L_h(t)}_{L^2\to L^2}\ge e^{-\mathcal C_St}$, while
 $\mathring L_h(t)f\to0$ in $L^2$ for every $f\in L^2$.
\item[(b)] For every sequence $h_k\to0$ and every sequence $r_k\to0$ of
 positive numbers there is a real $f\in L^2(0,1)$ with
 $\sup_k\norm{\mathring L_{h_k}(t)f}_2/r_k=\infty$.
\end{enumerate}
\end{corollary}
\begin{proof}[Sketch]
\emph{(a)} Since $f_h\in X_h$, $\mathring\Pi_hf_h=f_h$, and
$\norm{\mathring L_h(t)}\ge\norm{(\mathring T_h(t)-T(t))f_h}_2$. The lower
bound follows from Corollary~\ref{cor:lineargap}. The strong convergence
is Lemma~\ref{lem:TK}.

\emph{(b)} Let $h_k\to0$ and $r_k\to0$, and suppose that
$\sup_k\norm{\mathring L_{h_k}(t)f}_2/r_k<\infty$ for every real
$f\in L^2(0,1)$. The real functions form a real Banach space. By the
uniform boundedness principle on this space, there is $C$ with
$\norm{\mathring L_{h_k}(t)f}_2\le Cr_k\norm f_2$ for all $k$ and all real
$f$. The functions $f_{h_k}$ are real, with $\norm{f_{h_k}}_2=1$. Hence
$\norm{\mathring L_{h_k}(t)f_{h_k}}_2\le Cr_k\to0$. This contradicts
Corollary~\ref{cor:lineargap}.
\end{proof}
The lower bound in (a) also follows from
Corollary~\ref{cor:principleS}(c), without Corollary~\ref{cor:lineargap}.
For the W-scheme, $\norm{L_h(t)}_{H\to H}\le Ce^{-\nu t}ht_\wedge^{-3/2}$ for $h\le t_\wedge^{3/2}$
\cite[Theorem~4.2]{Main}, so the rate $h$ holds for every datum.

\begin{remark}[the role of Lean]\label{rem:leanscope}
For the existence of the high-frequency spurious eigenpairs,
Theorem~\ref{thm:top}, the Lean proof is essential. The properties of the
eigenpairs, Corollary~\ref{lem:density} and Lemmas~\ref{lem:slow}
and~\ref{lem:fast}, can be derived by hand, but the derivations are long.
The linear approximation results, Lemmas~\ref{lem:inverse},
\ref{lem:twoscale} and~\ref{lem:TK}, Propositions~\ref{prop:galerkin}
and~\ref{prop:necessary} and
Corollaries~\ref{cor:principleS}--\ref{cor:fixed}, are simple, and the
sketches above suffice for a proof by hand.

We have covered in Lean as much as was practical. All results of this
section have Lean~4 proofs, checked by the Lean kernel. So have the
results of Section~\ref{sec:gap}, except
Proposition~\ref{thm:quad}(b) and the last assertion of
Proposition~\ref{prop:Ulimit}; they also have complete proofs by hand. A
single premise remains: the mesh-uniform Crandall--Pazy estimate of the
W-scheme, \cite[(3.15)]{Main} with $\alpha=1$ on uniform meshes.

This premise is the main discrete estimate of~\cite{Main}. Its place in
the Lean proof shows what separates the two schemes: the W-scheme has
mesh-uniform regularity, and the S-scheme does not. The Lean statements
and their translation are given in Appendix~\ref{app:lean}.
\end{remark}

\section{The error gap}\label{sec:gap}

In this section, we prove the error gap for the nonlinear problem. The
proofs are by hand and complete. They use the results stated in
Section~\ref{sec:spurious} as statements, in the same way as
Theorem~\ref{thm:top}. These are Corollary~\ref{lem:density}, the
positivity \eqref{eq:bvalue} of the cell constant,
Lemmas~\ref{lem:slow}, \ref{lem:fast} and~\ref{lem:TK},
Corollary~\ref{cor:lineargap}, and the compactness of $T(t)$ from
Section~\ref{sec:galerkin}. They do not use steps of the Lean proofs.

The route is as follows. Section~\ref{sec:Ulimit} shows that the
solution of \eqref{eq:pde} for small data is the limit of the W-scheme.
Section~\ref{sec:quad} computes the limit $E$ of the quadratic Taylor
coefficient of the S-scheme on the data $\tilde f_h$, and shows that
$E(t)\ne0$ for $t>0$. Section~\ref{sec:gapthm} proves the error gap,
Theorem~\ref{thm:total}, and shows that the ordering error of the
S-scheme admits no product bound.

Throughout, $L=1$ and $h=1/N$ with $N\ge512$. The eigenpairs
$(\zeta_h,\varphi_h)$, the data $\tilde f_h$ and $f_h$, the functions
$\mathfrak d$ and $\mathfrak p$, and the conventions on constants are
those of Section~\ref{sec:eigfun}. Recall that
$\mu_h=\Rea\zeta_h\to-\mathcal C_S$ by \eqref{eq:re}.

\subsection{The nonlinear solution as a Galerkin limit}\label{sec:Ulimit}
The error gap compares the S-scheme with the solution $u$ of
\eqref{eq:pde}. In~\cite{Main} this solution is the sum of a power series
in the datum~\cite[Theorem~2.2]{Main}. By
Proposition~\ref{prop:galerkin}, $A^{-1}$ and $T$ are limits of the
W-scheme. The next proposition shows the same for $u$. Thus the gap is a
statement about two Galerkin schemes on the same trial space.

We use the notation of~\cite{Main} for $L=1$ and uniform meshes. Let
$F(f,g)=-\tfrac12(fg)'$ and $F_h=\Pi_hF$. The W-scheme is
\begin{equation}\label{eq:Wscheme}
 u_{h,t}=A_hu_h+F_h(u_h,u_h),\qquad u_h(0)=\Pi_hu_0 .
\end{equation}
Fix $0<\nu<3\pi^2/4$, write $|f|_1=\norm{f'}_2$, and put
\[
 \mathcal E_\nu(w)=\sup_{t\ge0}e^{\nu t}\norm{w(t)}_H
 +\Bigl(\int_0^\infty e^{2\nu t}|w(t)|_1^2\,dt\Bigr)^{1/2}.
\]
Let $\mathcal X_0$ be the energy trajectory space
of~\cite[(2.11)]{Main}. It consists of the $w\in C([0,\infty);H)$ with
$w(t)\in H_0^1(0,1)$ for almost every $t$, $\mathcal E_\nu(w)<\infty$ and
$e^{\nu t}\norm{w(t)}_H\to0$ as $t\to\infty$, and $\mathcal E_\nu$ is its
norm.

\begin{proposition}[the nonlinear solution as a Galerkin limit]\label{prop:Ulimit}
Let $L=1$, let the meshes be uniform, and fix $0<\nu<3\pi^2/4$. There are
$r_0,R_0>0$ such that the following hold for every real $f\in L^2(0,1)$
with $\norm f_H<r_0$.
\begin{enumerate}
\item[(a)] For every $N$ the W-scheme \eqref{eq:Wscheme} with $u_0=f$ has
a unique global solution $u_h$, $u_h\in\mathcal X_0$ and $\mathcal E_\nu(u_h)\le R_0$.
\item[(b)] There is $\mathcal U(f)\in\mathcal X_0$ such that
$\mathcal E_\nu(u_h-\mathcal U(f))\to0$ as $h\to0$. In particular,
$u_h(t)\to\mathcal U(f)(t)$ in $L^2$, uniformly in $t\ge0$.
\item[(c)] The function $w=\mathcal U(f)$ is real,
$\mathcal E_\nu(w)\le R_0$, and
\begin{equation}\label{eq:mild}
 w(t)=T(t)f-\int_0^tT(t-s)\bigl[w(s)\,\partial_xw(s)\bigr]\,ds,\qquad t\ge0,
\end{equation}
where the integral converges in $L^2(0,1)$. Moreover, $w$ is the only
solution of \eqref{eq:mild} in
$\{w\in\mathcal X_0:\mathcal E_\nu(w)\le R_0\}$.
\end{enumerate}
Consequently $\mathcal U(f)$ is the solution of \eqref{eq:pde} given
by~\cite[Theorem~2.2]{Main} with $\theta=0$ and $u_0=f$.
\end{proposition}
\begin{proof}
For $w,\tilde w\in\mathcal X_0$ we have
$\norm{F(w(s),\tilde w(s))}_H\le C|w(s)|_1|\tilde w(s)|_1$
by~\cite[(3.12)]{Main}. Hence $s\mapsto F(w(s),\tilde w(s))$ is
integrable on $(0,\infty)$ with values in $H$, also with the weight
$e^{\nu s}$. We write $\mathcal L_hf=T_h(\cdot)\Pi_hf$,
$\mathcal Lf=T(\cdot)f$, and
\[
 \mathcal Q_h(w,\tilde w)(t)=\int_0^tT_h(t-s)\Pi_hF(w(s),\tilde w(s))\,ds,\qquad
 \mathcal Q(w,\tilde w)(t)=\int_0^tT(t-s)F(w(s),\tilde w(s))\,ds .
\]

\emph{Step 1: uniform bounds.} Let $V_1^*$ be the dual of $H_0^1(0,1)$
with respect to the pairing of $H$. By~\cite[(3.12)]{Main},
$\norm{F(w,\tilde w)}_{V_1^*}\le C|w|_1\norm{\tilde w}_H$ on $H_0^1(0,1)$.
This uses only $\norm w_\infty\le|w|_1$. Agmon's inequality is not
needed. Hence, for $w,\tilde w\in\mathcal X_0$,
\[
 \Bigl(\int_0^\infty e^{2\nu t}\norm{F(w,\tilde w)}_{V_1^*}^2\,dt\Bigr)^{1/2}
 \le C\Bigl(\int_0^\infty e^{2\nu t}|w|_1^2\,dt\Bigr)^{1/2}\sup_{t\ge0}\norm{\tilde w(t)}_H
 \le C\,\mathcal E_\nu(w)\,\mathcal E_\nu(\tilde w).
\]
By the forced energy estimates~\cite[(2.19) and (3.11)]{Main}, there are
constants $c_1,c_2>0$, independent of $h$, such that
\begin{equation}\label{eq:LQbounds}
 \mathcal E_\nu(\mathcal L_hf),\ \mathcal E_\nu(\mathcal Lf)\le c_1\norm f_H,\qquad
 \mathcal E_\nu(\mathcal Q_h(w,\tilde w)),\ \mathcal E_\nu(\mathcal Q(w,\tilde w))
 \le c_2\,\mathcal E_\nu(w)\,\mathcal E_\nu(\tilde w).
\end{equation}
These trajectories belong to $\mathcal X_0$ (the decay
$e^{\nu t}\norm{\cdot}_H\to0$ follows from the same estimates on $[0,t]$
with a weight $\nu'\in(\nu,3\pi^2/4)$ and dominated convergence).
Put $r_0=(8c_1c_2)^{-1}$ and $R_0=(4c_2)^{-1}$.

\emph{Step 2: the discrete solutions.} Let $P_{1,h}=\mathcal L_hf$ and
$P_{n,h}=\sum_{j+k=n}\mathcal Q_h(P_{j,h},P_{k,h})$ for $n\ge2$. By
induction and \eqref{eq:LQbounds}, with the Catalan numbers bounded by
$4^{n-1}$ as in~\cite[Section~2.4]{Main},
\begin{equation}\label{eq:Pnbound}
 \mathcal E_\nu(P_{n,h})\le4^{n-1}c_2^{n-1}c_1^n\norm f_H^n
 =R_0\Bigl(\frac{\norm f_H}{2r_0}\Bigr)^n .
\end{equation}
For $\norm f_H<r_0$ the series $u_h=\sum_{n\ge1}P_{n,h}$ converges in
$\mathcal X_0$, and $\mathcal E_\nu(u_h)\le R_0$. Rearranging the
absolutely convergent double series
$\sum_{j,k}\mathcal Q_h(P_{j,h},P_{k,h})$, we obtain
$u_h=\mathcal L_hf+\mathcal Q_h(u_h,u_h)$. This is the mild form of
\eqref{eq:Wscheme} in the finite-dimensional space $X_h$. Hence $u_h$ is
continuously differentiable and solves \eqref{eq:Wscheme} for all
$t\ge0$. The right side of \eqref{eq:Wscheme} is a polynomial, so the
solution is unique. It is real, since the scheme has real coefficients.
This proves (a). It is the argument of~\cite[Proposition~3.6]{Main}. The
global existence comes from the series. The a priori bound
of~\cite[Proposition~3.2]{Main} is not needed.

\emph{Step 3: convergence of the linear maps.}
By~\cite[Lemma~4.8]{Main},
$\mathcal E_\nu(\mathcal L_hx-\mathcal Lx)\le Ch^2\norm{Ax}_H$ for
$x\in D(A)$. Since $D(A)$ is dense in $H$, \eqref{eq:LQbounds} gives
$\mathcal E_\nu(\mathcal L_hx-\mathcal Lx)\to0$ for every $x\in H$. Now
let $w,\tilde w\in\mathcal X_0$ and $g(s)=F(w(s),\tilde w(s))$. Then
\[
 \bigl(\mathcal Q_h(w,\tilde w)-\mathcal Q(w,\tilde w)\bigr)(t)
 =\int_0^t\bigl(T_h(t-s)\Pi_h-T(t-s)\bigr)g(s)\,ds .
\]
A time shift by $s$ multiplies the weighted norm by $e^{\nu s}$. Hence,
by Minkowski's inequality,
\[
 \mathcal E_\nu\bigl(\mathcal Q_h(w,\tilde w)-\mathcal Q(w,\tilde w)\bigr)
 \le\int_0^\infty e^{\nu s}\,\mathcal E_\nu\bigl((\mathcal L_h-\mathcal L)g(s)\bigr)\,ds .
\]
The integrand tends to zero for almost every $s$. It is at most
$2c_1e^{\nu s}\norm{g(s)}_H$, which is integrable. By dominated
convergence,
$\mathcal E_\nu(\mathcal Q_h(w,\tilde w)-\mathcal Q(w,\tilde w))\to0$ for
fixed $w,\tilde w\in\mathcal X_0$.

\emph{Step 4: the limit.} Let $P_1=\mathcal Lf$ and
$P_n=\sum_{j+k=n}\mathcal Q(P_j,P_k)$. These are the Taylor coefficients
of~\cite[Theorem~2.2]{Main} with $\theta=0$. They satisfy
\eqref{eq:Pnbound} as well. We show
$\mathcal E_\nu(P_{n,h}-P_n)\to0$ by induction on $n$. For $n=1$ this is
Step~3. For $j+k=n$,
\[
 \mathcal Q_h(P_{j,h},P_{k,h})-\mathcal Q(P_j,P_k)
 =\mathcal Q_h(P_{j,h}-P_j,P_{k,h})+\mathcal Q_h(P_j,P_{k,h}-P_k)
 +\bigl(\mathcal Q_h-\mathcal Q\bigr)(P_j,P_k).
\]
The first two terms tend to zero by the induction hypothesis and
\eqref{eq:LQbounds}, the third by Step~3. By \eqref{eq:Pnbound},
$\mathcal E_\nu(P_{n,h}-P_n)\le2R_0(\norm f_H/(2r_0))^n$. This bound is
summable and independent of $h$. So we may pass to the limit term by
term. Hence $\mathcal E_\nu(u_h-\mathcal U(f))\to0$ with
$\mathcal U(f)=\sum_{n\ge1}P_n$. This proves (b), since
$\norm\cdot_2\le\norm\cdot_H$.

\emph{Step 5: the equation and uniqueness.} As in Step~2,
$w=\mathcal U(f)$ satisfies $\mathcal E_\nu(w)\le R_0$ and
$w=\mathcal Lf+\mathcal Q(w,w)$. It is real, as a limit of real
functions. Since $F(w,w)=-w\,\partial_xw$ and
$\norm{w\,\partial_xw}_2\le|w|_1^2$, the equation
$w=\mathcal Lf+\mathcal Q(w,w)$ is \eqref{eq:mild}, with an integral that
converges in $L^2(0,1)$. Let $\tilde w$ be another solution with
$\mathcal E_\nu(\tilde w)\le R_0$. Then, by \eqref{eq:LQbounds},
\[
 \mathcal E_\nu(w-\tilde w)
 =\mathcal E_\nu\bigl(\mathcal Q(w-\tilde w,w)+\mathcal Q(\tilde w,w-\tilde w)\bigr)
 \le2c_2R_0\,\mathcal E_\nu(w-\tilde w)=\tfrac12\mathcal E_\nu(w-\tilde w),
\]
so $\tilde w=w$. This proves (c). The series $\sum_{n\ge1}P_n$ is the
series of~\cite[Theorem~2.2]{Main}. This proves the last assertion.
\end{proof}

\subsection{The quadratic coefficient}\label{sec:quad}
For real $f$ let
\[
 P_2(f)(t)=-\tfrac12\int_0^tT(t-s)\partial_x\bigl[(T(s)f)^2\bigr]ds,\qquad
 \mathring P_{2,h}(f)(t)=\int_0^t\mathring T_h(t-s)\mathcal G_h\bigl[(\mathring T_h(s)f)^2\bigr]ds
\]
be the quadratic Taylor coefficients of the solution maps, the second for
$f\in X_h$, and put $\mathring D_h(f)=P_2(f)-\mathring P_{2,h}(\mathring\Pi_hf)$.
For the W-scheme, the discretization error of the quadratic Taylor coefficient
$D_h(u_0)=P_2(u_0)-P_{2,h}(\Pi_hu_0)$ is defined
in~\cite[Theorem~4.3]{Main}, which gives
$\sup_{t\ge0}e^{\nu t}\norm{\Cd_h(t)-D_h(u_0)(t)}_H\le Ch^{1/6}\varepsilon^3$ for
$\varepsilon=\norm{u_0}_H$ small. Thus $D_h$ gives the ordering error to
leading order.

The next proposition identifies the limit of the quadratic Taylor
coefficient of the S-scheme on the data $\tilde f_h$.

\begin{proposition}[quadratic coefficient]\label{thm:quad}
Recall $\tilde f_h$, $\mathfrak d$, $\mathfrak p$, $\mathfrak b$, $z_\infty$
and $\mathcal C_S$ from Table~\ref{tab:notation}. Put
\begin{equation}\label{eq:Edef}
 E(t)=T(t)z_\infty-e^{-2\mathcal C_St}z_\infty
 -2\mathcal C_S\int_0^tT(t-s)z_\infty e^{-2\mathcal C_Ss}\,ds .
\end{equation}
Then the following hold.
\begin{enumerate}
\item[(a)] For every $T_0>0$,
 $\sup_{0\le t\le T_0}\norm{\mathring P_{2,h}(\tilde f_h)(t)-E(t)}_2\to0$ as $h\to0$.
\item[(b)] For $t>0$, $E(t)\in H^3(0,1)$, $E(0,t)=E(1,t)=0$ and
 $E_x(1,t)=-\mathfrak be^{-2\mathcal C_St}$. Moreover $E(0)=0$ in $L^2$, and
 $E_t+E_{xxx}=-\tfrac12\mathfrak d'(x)e^{-2\mathcal C_St}$ in the sense of distributions.
\item[(c)] $E(t)\ne0$ for every $t>0$.
\end{enumerate}
\end{proposition}
\begin{proof}
\emph{(a)} We first split the quadratic coefficient into a slow part and
an oscillating part. Write $\mu=\mu_h$, $\zeta=\zeta_h$, $w=e^{\zeta s}\varphi_h$.
Then $\mathring T_h(s)\tilde f_h=\sqrt2\,\Rea w$ and $2(\Rea w)^2=|w|^2+\Rea w^2$.
Let $z_h$ be as in Lemma~\ref{lem:slow} and $\mathcal O_h$ as in
Lemma~\ref{lem:fast}. Since $\mathcal G_h$ and $\mathring T_h(t)$ commute
with complex conjugation,
\[
 \mathring P_{2,h}(\tilde f_h)(t)=\Sigma_h(t)+\Rea\mathcal O_h(t),\qquad
 \Sigma_h(t)=\int_0^t\mathring T_h(t-s)\mathring A_hz_h\,e^{2\mu s}ds .
\]
Differentiating $s\mapsto\mathring T_h(t-s)z_he^{2\mu s}$ and integrating, we obtain
\begin{equation}\label{eq:SN}
 \Sigma_h(t)=\mathring T_h(t)z_h-e^{2\mu t}z_h+2\mu\int_0^t\mathring T_h(t-s)z_he^{2\mu s}ds .
\end{equation}
By Lemma~\ref{lem:slow}, $\norm{z_h-\mathring\Pi_hz_\infty}_2\le\norm{z_h-z_\infty}_2\le Ch$. So contraction lets us replace
$\mathring T_h(\cdot)z_h$ by $\mathring T_h(\cdot)\mathring\Pi_hz_\infty$ with error $O(h)$.
Combining this with Lemma~\ref{lem:TK} and $\mu_h\to-\mathcal C_S$, we conclude that
$\Sigma_h\to E$ uniformly on $[0,T_0]$.
Lemma~\ref{lem:fast} removes $\mathcal O_h$.

\emph{(b)} We now prove the regularity and the equation of $E$. Put
\[
 \mathcal H(t)=\int_0^tT(t-s)z_\infty e^{-2\mathcal C_Ss}\,ds .
\]
We claim that
\begin{equation}\label{eq:HE}
 \mathcal H(t)=A^{-1}E(t),\qquad t\ge0 .
\end{equation}
Indeed, $\frac d{dr}A^{-1}T(r)z_\infty=T(r)z_\infty$, so
\[
 \frac d{ds}\Bigl[e^{-2\mathcal C_Ss}A^{-1}T(t-s)z_\infty\Bigr]
 =-e^{-2\mathcal C_Ss}\bigl(2\mathcal C_SA^{-1}+I\bigr)T(t-s)z_\infty .
\]
Integrating over $[0,t]$, we obtain
\[
 e^{-2\mathcal C_St}A^{-1}z_\infty-A^{-1}T(t)z_\infty
 =-2\mathcal C_SA^{-1}\mathcal H(t)-\mathcal H(t).
\]
By \eqref{eq:Edef}, this is \eqref{eq:HE}. In particular
$\mathcal H(t)\in D(A)$.

Put $Y=-\mathfrak d/(8\pi^2)$. Then $Y\in D(A)$ and $AY=-\tfrac12\mathfrak d'$.
Let $\mathcal Y(t)=E(t)+\mathfrak b e^{-2\mathcal C_St}\mathfrak p$. By \eqref{eq:Edef} and $z_\infty=Y+\mathfrak b\mathfrak p$,
\begin{equation}\label{eq:Y}
 \mathcal Y(t)=T(t)z_\infty-e^{-2\mathcal C_St}Y-2\mathcal C_S\mathcal H(t).
\end{equation}
For $t>0$, $T(t)z_\infty\in D(A)$ by \cite[Theorem~2.4(ii)]{Main} with $\alpha=1$. Hence $\mathcal Y(t)\in D(A)$.
This yields the boundary values, since $\mathfrak p(0)=\mathfrak p(1)=0$ and $\mathfrak p'(1)=1$.
The relation $E(0)=0$ is immediate. No compatibility of the boundary trace at $t=0$ is claimed.
For $\vartheta\in C_c^\infty(0,1)$, it follows from \eqref{eq:Edef} that
$\frac d{dt}(E,\vartheta)=(E,\vartheta''')+e^{-2\mathcal C_St}(z_\infty,\vartheta''')$. Moreover
$(z_\infty,\vartheta''')=-(z_\infty''',\vartheta)=-\tfrac12(\mathfrak d',\vartheta)$. This is the equation.

\emph{(c)} Suppose $E(t_0)=0$ for some $t_0>0$. By \eqref{eq:HE},
$\mathcal H(t_0)=0$. Then \eqref{eq:Edef} gives
$T(t_0)z_\infty=e^{-2\mathcal C_St_0}z_\infty$. The left side lies in
$D(A)$, as noted in (b). Hence $z_\infty\in D(A)$ and
$\mathfrak b\mathfrak p=z_\infty-Y\in D(A)$. This is false, because
$\mathfrak b\ne0$ by \eqref{eq:bvalue} and $\mathfrak p'(1)=1$.
\end{proof}

The proposition has the following consequence for the discretization
error $\mathring D_h$.

\begin{corollary}\label{cor:Dh}
Recall $f_h$ and $E$ from Table~\ref{tab:notation}. Then
$\norm{f_h}_2=1$ and $f_h\rightharpoonup0$ in $L^2$. For every $T_0>0$,
\[
 \sup_{0\le t\le T_0}\norm{P_2(f_h)(t)}_2\to0,\qquad
 \sup_{0\le t\le T_0}\norm{\mathring D_h(f_h)(t)+E(t)}_2\to0 .
\]
Hence $\lim_{h\to0}\norm{\mathring D_h(f_h)(t)}_2=\norm{E(t)}_2>0$ for every $t>0$.
\end{corollary}
\begin{proof}
By Corollary~\ref{lem:density}(d), $\norm{\tilde f_h}_2\to1$ and $f_h\rightharpoonup0$.
Let $v=T(\cdot)f_h$. For $s>0$, $T(s)$ is compact
(Section~\ref{sec:galerkin}). So $v(s)\to0$ in $L^2$.
By dominated convergence, $\norm v_{L^2(0,T_0;L^2)}\to0$.
Kato smoothing \cite[Proposition~2.1]{Main} bounds $\norm{v_x}_{L^2(0,T_0;L^2)}$.
Now $\norm v_\infty^2\le2\norm v_2\norm{v_x}_2$, and $T$ is a contraction on $L^2$. So
\[
 \norm{P_2(f_h)(t)}_2\le C\int_0^{T_0}\norm v_2^{1/2}\norm{v_x}_2^{3/2}ds
 \le C\norm v_{L^2(0,T_0;L^2)}^{1/2}\norm{v_x}_{L^2(0,T_0;L^2)}^{3/2}\to0 .
\]
Since $\mathring P_{2,h}$ is homogeneous of degree two,
$\mathring P_{2,h}(f_h)=\mathring P_{2,h}(\tilde f_h)/\norm{\tilde f_h}_2^2$.
The rest follows from Proposition~\ref{thm:quad}.
\end{proof}
For the W-scheme, \cite[Theorem~4.3]{Main} and the homogeneity of the quadratic
difference give $\norm{D_h(f)(t)}_H\le Ch^{1/2}t_\wedge^{-5/2}e^{-\nu t}\norm f_H^2$
for real $f$. The discretization error of the quadratic Taylor
coefficient therefore tends to zero uniformly on the unit ball for the W-scheme, for each
fixed $t>0$ and uniformly for $t\ge t_0>0$, and does not tend to zero for the S-scheme.

\begin{remark}[the trace of the cell structure]\label{rem:cell}
The term $\mathfrak b\mathfrak p$ in $z_\infty$ is the only trace of the
cell structure of the eigenfunctions. The weak limit
$|\varphi_h|^2\rightharpoonup\mathfrak d$ alone would suggest the limit
$-\tfrac12\int_0^tT(t-s)\mathfrak d'e^{-2\mathcal C_Ss}\,ds$ in
Proposition~\ref{thm:quad}(a), without this term. The term arises because
the load is tested against derivatives that oscillate on the mesh scale,
namely the derivative of the grid mode $\mathsf m_h$
(Lemma~\ref{lem:twoscale}). It comes from
$-g(1)\tfrac12(f,\mathsf m_h')$ in \eqref{eq:inverseload}. In
Proposition~\ref{thm:quad}(b) it appears as the boundary flux
$E_x(1,t)=-\mathfrak b\,e^{-2\mathcal C_St}$. The proof of
Proposition~\ref{thm:quad}(c) rests on it, since $\mathfrak p\notin D(A)$
and $\mathfrak b>0$. The numerical results in Appendix~\ref{app:num} are
consistent with Proposition~\ref{thm:quad}(a): the computed differences
decrease by a factor of about two per halving of $h$.
\end{remark}

\subsection{The error gap}\label{sec:gapthm}
\begin{theorem}[error gap]\label{thm:total}
Recall $f_h$ and $\mathcal C_S$ from Table~\ref{tab:notation}.
Part (a) concerns the linear error. Part (b) concerns the total error.
Part (c) concerns the full error.
\begin{enumerate}
\item[(a)] For every $0<t_0<T_0$,
 $\sup_{t_0\le t\le T_0}\bigl|\norm{(T(t)-\mathring T_h(t))f_h}_2-e^{-\mathcal C_St}\bigr|\to0$.
 This is Corollary~\ref{cor:lineargap}.
\item[(b)] There are $\varepsilon_0,C>0$ such that for $0<\varepsilon\le\varepsilon_0$, $u_0=\varepsilon f_h$
 and every $t>0$,
 \[
  \limsup_{h\to0}\bigl|\norm{(u-\mathring v_h)(t)}_2-\varepsilon e^{-\mathcal C_St}\bigr|\le C\varepsilon^3 .
 \]
\item[(c)] Let $\varepsilon_h>0$ with $\varepsilon_h/h\to0$, and $u_0=\varepsilon_hf_h$. Then for every $0<t_0<T_0$,
 \[
  \sup_{t_0\le t\le T_0}\Bigl|\varepsilon_h^{-1}\norm{(u-\mathring u_h)(t)}_2-e^{-\mathcal C_St}\Bigr|\to0 .
 \]
\end{enumerate}
\end{theorem}
\begin{proof}
\emph{(a)} This is Corollary~\ref{cor:lineargap}.

\emph{(b)} Let $u_0=\varepsilon f_h$. We fix $0<\nu<3\pi^2/4$ and apply
\cite[Theorem~2.2]{Main} with $\theta=0$ and $m=1,2$. Let $r$ be its
radius and $r_0$ the radius of Proposition~\ref{prop:Ulimit}. Take
$\varepsilon_0=\min\{r,r_0\}/(2\sqrt2)$, so that
$\norm{u_0}_H\le\min\{r,r_0\}/2$. By Proposition~\ref{prop:Ulimit}, the
solution $u$ is also the limit of the W-scheme. Moreover,
\begin{equation}\label{eq:taylor}
 \norm{u(t)-\varepsilon T(t)f_h}_2\le C\varepsilon^2,\qquad
 \norm{u(t)-\varepsilon T(t)f_h-\varepsilon^2P_2(f_h)(t)}_2\le C\varepsilon^3,
\end{equation}
with $C$ independent of $N$ and $t$. Since
$\mathring v_h=\varepsilon\mathring T_h(\cdot)f_h$, the second estimate
gives
\[
 \Bigl|\norm{(u-\mathring v_h)(t)}_2
 -\varepsilon\norm{(T(t)-\mathring T_h(t))f_h}_2\Bigr|
 \le\varepsilon^2\norm{P_2(f_h)(t)}_2+C\varepsilon^3 .
\]
By Corollary~\ref{cor:Dh}, $\norm{P_2(f_h)(t)}_2\to0$. Together with
(a), this gives (b).

\emph{(c)} It remains to treat the nonlinear scheme. Let
$\varepsilon=\varepsilon_h$, $u_0=\varepsilon f_h$, $v=\mathring v_h$ and
$w=\mathring u_h-\mathring v_h$. Subscripts $t,x$ denote time and space
derivatives. Then $w(0)=0$ and
$w_t=\mathring A_hw+\mathcal G_h[(v+w)^2]$. We test with $w$. All
functions are real, $C^1$ in $x$ and vanish at both ends, so
\[
 \tfrac12(w^2,w_x)=0,\qquad (vw,w_x)=-\tfrac12(v_x,w^2),\qquad \tfrac12(v^2,w_x)=-(vv_x,w).
\]
Together with \eqref{eq:flux}, this yields
\[
 \tfrac12\frac d{dt}\norm w_2^2=-\tfrac12|w_x(0,t)|^2-\tfrac12(v_x,w^2)-(vv_x,w).
\]
We have
$v=\sqrt2\,\varepsilon\Rea(e^{\zeta_ht}\varphi_h)/\norm{\tilde f_h}_2$.
By Corollary~\ref{lem:density}(a,d), $\norm v_\infty\le\varkappa_0$ and
$\norm{v_x}_\infty\le\varkappa_1$ for all $t\ge0$, with
$\varkappa_0=C\varepsilon$ and $\varkappa_1=C\varepsilon h^{-1}$. Since
$2\varkappa_0\norm w_2\le\varkappa_0^2+\norm w_2^2$,
\[
 \frac d{dt}\norm w_2^2\le\varkappa_1\norm w_2^2+2\varkappa_0\varkappa_1\norm w_2
 \le\varkappa_0^2\varkappa_1+2\varkappa_1\norm w_2^2 .
\]
By Gronwall's inequality,
$\norm{w(t)}_2^2\le\varkappa_0^2\varkappa_1\,t\,e^{2\varkappa_1t}$, that
is,
\[
 \norm{w(t)}_2\le C\varepsilon\Bigl(\frac{\varepsilon t}h\Bigr)^{1/2}e^{C\varepsilon t/h}.
\]
So $\sup_{t\le T_0}\norm{w(t)}_2=o(\varepsilon)$, because
$\varepsilon/h\to0$. For $N$ large, $\varepsilon\le\varepsilon_0$. Now
\[
 u-\mathring u_h=\bigl(u-\varepsilon T(\cdot)f_h\bigr)
 +\varepsilon\bigl(T(\cdot)-\mathring T_h(\cdot)\bigr)f_h-w .
\]
By the first estimate in \eqref{eq:taylor}, the first term is at most
$C\varepsilon^2$ in $L^2$. Hence
\[
 \sup_{0\le t\le T_0}\Bigl|\varepsilon^{-1}\norm{(u-\mathring u_h)(t)}_2
 -\norm{(T(t)-\mathring T_h(t))f_h}_2\Bigr|\to0 .
\]
The claim follows from (a).
\end{proof}

For the W-scheme, with $\varepsilon=\norm{u_0}_H$ small and $0<\nu<3\pi^2/4$,
\cite[Theorem~4.2]{Main} gives
$\norm{(v-v_h)(t)}_H\le Ce^{-\nu t}a(h,t)\varepsilon$ and
$\norm{(u-u_h)(t)}_H\le Ce^{-\nu t}\{a(h,t)\varepsilon+h^{1/6}\varepsilon^2\}$,
with $a(h,t)=\min\{1,ht_\wedge^{-3/2}\}$.
There the error divided by $\norm{u_0}_H$ tends to zero at each fixed $t>0$,
uniformly on a ball. The S-scheme has no such bound. Fix $t_*>0$ and $R>0$.
By (c) with $\varepsilon_h=h^2$, no constant $C$ independent of $h$ and $u_0$
makes the bound for $u-\mathring u_h$ hold for all real
$u_0\in X_h$ with $0<\norm{u_0}_H\le R$. The same argument excludes the
analogue of the uniform $L^4$ estimate in~\cite[Theorem~4.1(i)]{Main}. The data in (c) tend to zero with $h$,
so (c) concerns estimates linear in the datum.

By Corollary~\ref{cor:Dh}, on the data $f_h$ the discretization error of
the quadratic Taylor coefficient of the S-scheme does not tend to zero. This excludes a
product bound for the ordering error itself, uniform in the mesh and in
the size of the datum.

\begin{corollary}[the ordering error of the S-scheme]\label{cor:noproduct}
Recall $f_h$, $E$ and $\mathring D_h$ from Table~\ref{tab:notation}. Fix $t>0$.
Let $\mathring\Cd_h(t)$ be the ordering error of the S-scheme for the datum
$u_0=\varepsilon f_h$.
\begin{enumerate}
\item[(a)] For every $N\ge512$, $\varepsilon^{-2}\mathring\Cd_h(t)\to\mathring D_h(f_h)(t)$
in $L^2$ as $\varepsilon\downarrow0$. Hence
\[
 \lim_{h\to0}\lim_{\varepsilon\downarrow0}\varepsilon^{-2}\mathring\Cd_h(t)=-E(t)\ne0
 \quad\text{in }L^2 .
\]
\item[(b)] Let $p>0$ and $\varepsilon_1>0$, and let $\varepsilon_0$ be as
in Theorem~\ref{thm:total}(b). There is no constant $C$ such that
$\norm{\mathring\Cd_h(t)}_2\le Ch^p\varepsilon^2$ for all $N\ge512$ and all
$0<\varepsilon\le\min\{\varepsilon_1,\varepsilon_0\}$.
\end{enumerate}
\end{corollary}
\begin{proof}
\emph{(a)} We fix $N\ge512$ and write $f=f_h$. The argument compares the Taylor
expansions in $\varepsilon$ of the two solutions. Recall that $\mathring\Pi_hf=f$. By
the estimates \eqref{eq:taylor} in the proof of Theorem~\ref{thm:total},
which use only \cite[Theorem~2.2]{Main}, the
continuous solution with datum $\varepsilon f$ satisfies
$u(t)=\varepsilon T(t)f+\varepsilon^2P_2(f)(t)+O(\varepsilon^3)$ in $L^2$ for
$0<\varepsilon\le\varepsilon_0$. For fixed $h$ the scheme is a quadratic
system of ordinary differential equations on the finite-dimensional space
$X_h$, so $\norm{\mathcal G_h[w^2]}_2\le C_h\norm w_2^2$ for $w\in X_h$. On
$C([0,t];X_h)$ with the supremum norm, the solution is therefore the sum
of the series in $\varepsilon$ generated by the recursion
of~\cite[Proposition~3.6]{Main}, for $\varepsilon$ small depending on $h$
and $t$. Its first two terms are $\varepsilon\mathring T_h(\cdot)f$ and
$\varepsilon^2\mathring P_{2,h}(f)$, and the remainder after two terms is
bounded by $C_{h,t}\varepsilon^3$. By uniqueness for the system of
ordinary differential equations, this sum is $\mathring u_h$. Thus
$\mathring u_h(t)=\varepsilon\mathring T_h(t)f+\varepsilon^2\mathring P_{2,h}(f)(t)+O_h(\varepsilon^3)$,
where the constant in $O_h$ may depend on $h$. Since
$\mathring v_h=\varepsilon\mathring T_h(\cdot)f$,
\[
 \mathring\Cd_h(t)=(u-\varepsilon T(t)f)-(\mathring u_h-\varepsilon\mathring T_h(t)f)
 =\varepsilon^2\mathring D_h(f)(t)+O_h(\varepsilon^3).
\]
Dividing by $\varepsilon^2$ and letting $\varepsilon\downarrow0$ at fixed $h$, we
obtain the first claim. The iterated limit then follows from
Corollary~\ref{cor:Dh}, which also gives $\norm{E(t)}_2>0$.

\emph{(b)} Suppose $\norm{\mathring\Cd_h(t)}_2\le Ch^p\varepsilon^2$ for all $N\ge512$ and
$0<\varepsilon\le\min\{\varepsilon_1,\varepsilon_0\}$. By (a), letting $\varepsilon\downarrow0$ at fixed
$h$, we obtain $\norm{\mathring D_h(f_h)(t)}_2\le Ch^p$ for all $N\ge512$.
As $h\to0$, this contradicts (a), since $\norm{E(t)}_2>0$.
The order of the limits matters: the argument gives no information on
$\mathring\Cd_h$ as $h\to0$ at a fixed $\varepsilon$.
\end{proof}
The corollary concerns bounds uniform in $h$ and $\varepsilon$. It says
nothing about the limit $h\to0$ of the ordering error at a fixed
amplitude. Table~\ref{tab:intro-numbers} reports computed values.

\section{A nonlinear correction on a coarse mesh}\label{sec:twomesh}

In this section, we add a nonlinear correction computed on a coarse mesh
to a linear simulation on a fine mesh. The ordering error has a direct
computational use. Since a nonlinear simulation costs more than a linear
one, the two can be run on different meshes. Let $h\le\hat h$ be the sizes
of two meshes as in~\cite[Section~3.1]{Main}, with the same constant
$\kappa$. They need not be nested. The constants below depend on $\kappa$.
Put
$\tilde u=v_h+(u_{\hat h}-v_{\hat h})$. The
definition of the ordering error on the mesh of size $\hat h$ gives
\[
 u-\tilde u=(v-v_h)+\Cd_{\hat h}.
\]
For the W-scheme, \cite[Corollary~4.5]{Main} gives, for $t>0$ and real data
with $\varepsilon=\norm{u_0}_H\le r$,
\[
 \norm{(u-\tilde u)(t)}_H\le Ce^{-\nu t}\bigl\{a(h,t)\varepsilon+\hat h^{1/6}\varepsilon^2\bigr\}.
\]
This has the form of the estimate for $u-u_h$, with $\hat h$ in place of
$h$ in the second term. The nonlinear mesh enters only at second order in
the datum.

\paragraph{When does the correction gain?}
We compare error bounds with error scales. Suppose that the linearization
share and the discretization errors on the coarse mesh are of sizes
$\norm{u-v}\approx\varepsilon^2$ and
$\norm{v-v_{\hat h}}\approx\norm{u-u_{\hat h}}\approx\hat h^q\varepsilon$.
Then the errors of $v_h$ and of $u_{\hat h}$ are of the orders
$\varepsilon^2$ and $\hat h^q\varepsilon$. Suppose also that the fine mesh
is so fine that $\norm{v-v_h}$ is much smaller than
$\min\{\varepsilon^2,\hat h^q\varepsilon\}$. Bounding the two
linearization shares $u-v$ and $u_{\hat h}-v_{\hat h}$ separately gives
$\norm{\Cd_{\hat h}}=O(\varepsilon^2)$. Bounding the full error
$u-u_{\hat h}$ and the linear error $v-v_{\hat h}$ separately gives
$\norm{\Cd_{\hat h}}=O(\hat h^q\varepsilon)$. Together they give only the
smaller of the two scales. So the bound for $u-\tilde u$ is not smaller in
order than both errors. A product bound $C\hat h^p\varepsilon^2$ is of
smaller order than $\varepsilon^2$ when $\hat h^p\ll1$, independently of
$\varepsilon$. It is of smaller order than $\hat h^q\varepsilon$ when
$\varepsilon\ll\hat h^{q-p}$. The bound of~\cite[Theorem~4.3]{Main} is of
this form with $p=1/6$, uniformly on a ball and for all time. In this
sense the estimate of the ordering error under mesh-uniform regularity is
what guarantees the gain in the error bound. For the S-scheme no product
bound holds (Corollary~\ref{cor:noproduct}). These comparisons concern
error scales. Comparing actual errors also requires control of their
sizes and of possible cancellation. The examples below compare actual
errors for the W-scheme.

The ordering error is also the gap in an error budget. Attributing the
total error $u-v_h$ to linearization and to discretization by the shares
$(u_h-v_h)$ and $(v-v_h)$ leaves
\[
 (u-v_h)-\bigl[(u_h-v_h)+(v-v_h)\bigr]=\Cd_h
\]
unexplained. The attribution is unambiguous only when both shares are much
larger than $\Cd_h$. This is the familiar residual of a two-factor
decomposition.

\paragraph{Examples.}
A nonlinear simulation costs more than a linear one on the same mesh,
since every time step requires an iteration. We take the W-scheme with $L=1$,
$\varepsilon=1/8$, the fine mesh $h=1/1024$ ($2047$ unknowns), the coarse
mesh $\hat h=1/64$ ($127$ unknowns), and the data $u_0=\varepsilon f$,
where $f$ is D1 (smooth), D2 (jump) or D3 (kink) of
Appendix~\ref{app:S2}, normalized by $\norm f_H=1$. Table~\ref{tab:twomesh-main} gives the errors in
$H$ at $t=0.01$, computed from the stored trajectories of
Appendix~\ref{app:S2}. Its columns give the errors of the linear solution
$v_h$ for the datum $u_0$, the nonlinear solution $u_{\hat h}$ on the
coarse mesh, the combined solution $\tilde u$ and the nonlinear solution
$u_h$ on the fine mesh. Appendix~\ref{app:num} writes the same linear
solution as $\varepsilon v_h$ under its unit-datum convention.

\begin{table}[H]
\centering\small
\caption{Errors in $H$ at $t=0.01$.}\label{tab:twomesh-main}
\begin{tabular}{@{}lcccc@{}}
\toprule
datum & $v_h$ & $u_{\hat h}$ & $\tilde u$ & $u_h$\\
\midrule
D1 (smooth) & $8.19\times10^{-5}$ & $6.37\times10^{-7}$ & $1.74\times10^{-9}$ & $2.17\times10^{-9}$\\
D2 (jump) & $1.21\times10^{-4}$ & $2.76\times10^{-6}$ & $6.29\times10^{-9}$ & $4.18\times10^{-9}$\\
D3 (kink) & $9.86\times10^{-5}$ & $8.64\times10^{-7}$ & $2.28\times10^{-9}$ & $3.65\times10^{-9}$\\
\bottomrule
\end{tabular}
\end{table}

The reference solution has an uncertainty of about $10^{-9}$. Smaller
differences are not distinguishable. Over all stored cases (D1--D3, three
amplitudes, four times), the ratio of the error of $\tilde u$ to that of
$u_h$ lies between $0.6$ and $4.9$ for $\hat h=1/64$. It reaches $14$ for
$\hat h=1/32$ and $73$ for $\hat h=1/16$.

For D2 we also measured running times (Appendix~\ref{app:twomesh}). The
four solutions were computed with the same Radau IIA method and with time
grids generated by the same rule. The numbers of steps differ by $0.3\%$.
Up to $t=0.1$, computing $u_h$ took $207$ seconds, $4.9$ times as long as
$v_h$. Each time step needed about two fixed-point iterations. The three
parts of $\tilde u$ took $99$ seconds when run one after another. They do
not depend on each other. Run simultaneously on three cores, they took
$47$ seconds of wall-clock time, about the time of $v_h$ alone. The errors
of $\tilde u$ and $u_h$ at $t=0.1$ are $1.0\times10^{-10}$ and
$6.6\times10^{-11}$.

\section*{Declarations}

\paragraph{Declaration of generative AI and AI-assisted technologies in the
manuscript preparation process.}
During the preparation of this work, the author used Claude Opus~5.5,
Claude Fable~5.1, Claude Fable~5.2 and GPT-6~Astra to check the
mathematical arguments in the initial draft, to search for relevant literature, to evaluate the
appropriateness of the English, to assist with mathematical
derivations, and to draft the reorganization of
Sections~\ref{sec:spurious} and~\ref{sec:gap} in the present version.
Their role in the proof of Theorem~\ref{thm:top} and in its
Lean formalization is described in Remark~\ref{rem:history}. They also
drafted the Lean proofs described in Remark~\ref{rem:leanscope}. Their role
in the numerical experiments, and the independent assessment of these
experiments, are described in Remark~\ref{rem:numerics}. The author did
not check the programs line by line. Appendix~\ref{app:num} states what
each computation computed and which checks it passed. The programs, their
outputs and logs are in the ancillary files.
The author reviewed and edited the output as needed and takes full
responsibility for the content of the published article, with two
exceptions. The proof of Theorem~\ref{thm:top} is checked by Lean
(Remark~\ref{rem:history}). The programs behind the numerical results
were not checked line by line by the author.

\paragraph{Data availability.}
The Lean files, the programs and their outputs are available as ancillary
files of arXiv:2611.12585.


\appendix
\input{appA_lean}

\input{appB_numerics}
\input{appZ_plain}

\end{document}

%% file: appA_lean.tex
\section{The Lean formalization}\label{app:lean}

This appendix may be skipped. The proof of Theorem~\ref{thm:top} is the
Lean proof, and Lean checks it. The handwritten arguments use
Theorem~\ref{thm:top}, \eqref{eq:zgv} and the intermediate results stated
in Section~\ref{sec:leanlemmas}. Sections~\ref{sec:A1}, \ref{sec:A3}
and~\ref{sec:A4} outline the formal proof of Theorem~\ref{thm:top},
translate its main statements and describe the files.
Section~\ref{sec:A5} states the Lean theorem that collects the
theoretical results. Its only premise is the mesh-uniform Crandall--Pazy
estimate of the W-scheme. Table~\ref{A5:corr} lists their Lean
declarations. These sections are meant for readers who want to examine
the Lean proofs. Notation is that of
Section~\ref{sec:disp}. Lean names without a prefix are in
\lean{SLean.Cert}, \lean{SLean.Algebra} or \lean{SLean.Explicit.Sharp};
Table~\ref{tab:corr} gives the full names.

\subsection{Outline of the proof}\label{sec:A1}
The outline follows the route that the Lean formalization took.
\begin{enumerate}
\item \emph{Stationary point.} Elimination in $c=\cos\xi$ reduces
$\varpi'(\xi)=0$ to $\tilde H(c)=0$. The quartic $\tilde H$ changes sign
on $[c_{\rm lo},c_{\rm hi}]$, an interval of length $2\cdot10^{-15}$.
Hence $c_*$ exists by the intermediate value theorem. Two polynomial
identities yield $\varpi(\xi_*)=w_*$ and $\varpi'(\xi_*)=0$.
Differentiating the explicit formula for $\varpi$, we obtain
$\varpi''(\xi_*)$ (\lean{exists_critical_c}, \lean{optical_thetaStar},
\lean{optical_stationary}, \lean{optical_second_deriv}).
\item \emph{Rational certificates.} The certificates of
Lemma~\ref{lem:cert} follow from exact rational arithmetic. The decaying roots $z_3$, $z_4$ are
located by a contraction mapping on an explicit disc. Rouch\'e's
theorem is not used, since Mathlib has none
(\lean{certificate_C1}--\lean{certificate_C5}).
\item \emph{Exact eigenvectors.} Exponential solutions solve the interior
recurrences. A zero of a $4\times4$ boundary determinant yields an
eigenpair of $\mathring A_h$ (Lemma~\ref{lem:det}; \lean{Sexp},
\lean{Sdet}).
\item \emph{Normal form and root branches.} Since $\varpi'(\xi_*)=0$ and
$\varpi''(\xi_*)>0$, the branch has a quadratic normal form near $\xi_*$.
It defines the companion map $\sigma$, which pairs two grid waves of equal
frequency. All four root branches are constructed on explicit discs
(Lemma~\ref{lem:domain}; \lean{exists_normal_coordinate},
\lean{exists_companion_germ}, \lean{roots_work}).
\item \emph{Quantization by contraction.} For $N\ge512$ the normalized
determinant leads to a fixed-point equation for the first mode. The
map is a contraction on an explicit disc (Lemmas~\ref{lem:tail}
and~\ref{lem:fixed}; \lean{analytic_fixed_point},
\lean{quantized_zero_exists}).
\item \emph{Explicit bounds.} A bootstrap yields $|y_N|\le8/N$
(Lemma~\ref{lem:box}). From a Horner certificate for an auxiliary
polynomial on a parameter box, we obtain $0<-L^3\Rea\zeta_h\le1600$. The
bound $|\Ima\zeta_h|\ge100h^{-3}$ follows from $w_*<-0.276$ and
$|\eta-\eta_*|\le3/1000$
(\lean{horner_upper_512}, \lean{damping_bound_512},
\lean{spectral_top_explicit_512}).
\item \emph{Part (b).} From the fixed-point equation, we obtain
$y_N=(i\beta\pi-c_2\pi^2)/(2N)$ up to an explicit $O(N^{-2})$ term. In the
expansion of $\tilde\varpi$ the terms with $c_2$ cancel, because
$a_3=a_2c_2$. Scaling by $420/h^3$ yields \eqref{eq:im}--\eqref{eq:re}
(\lean{fixed_point_leading_error}, \lean{optical_model_eq},
\lean{physical_expansions}).
\end{enumerate}

\subsection{Intermediate results}\label{sec:A2}
The intermediate results of the Lean proof that the paper uses are
stated in Section~\ref{sec:leanlemmas}, as
Lemmas~\ref{lem:cert}--\ref{lem:box}. Table~\ref{tab:corr} gives their
Lean names. There $\varsigma_N=(\pi+y_N)/N$ is defined. The eigenvalue of
Theorem~\ref{thm:top} is $\zeta_h=420i\tilde\varpi(\varsigma_N)/h^3$.

The translation procedure of Remark~\ref{rem:history} covered the
statements of Theorem~\ref{thm:top} and \eqref{eq:zgv}. It did not cover
the lemmas of Section~\ref{sec:leanlemmas}. Claude Opus~5.5 wrote their
statements from the Lean source. Claude Fable~5.1 compared the parts used
later in the paper with the Lean declarations. Apart from the hand steps
stated after Lemma~\ref{lem:det}, each part is the same as, or weaker
than, the corresponding Lean statement.

\subsection{The Lean statements and their translation}\label{sec:A3}
The translation below was produced by the procedure described in
Remark~\ref{rem:history}.

\paragraph{Theorem~\ref{thm:top}.}
The Lean statements are
\lean{SLean.Explicit.Sharp.Expansion.spectral_top_with_expansion}:
{\small
\begin{alltt}
theorem spectral_top_with_expansion \{L : \(\mathbb{R}\)\} (hL : 0 < L) : \(\forall\) N : \(\mathbb{N}\), 512 \(\le\) N \(\to\)
    \(\exists\) hN : 0 < N, \(\exists\) f : Statement.X_h N (Statement.h L N), f \(\ne\) 0 \(\wedge\)
      Statement.A_h hN (Statement.h L N) (Main.mesh_pos hL hN) f = Statement.\(\zeta\) L N \(\bullet\) f \(\wedge\)
      100/(Statement.h L N)^3 \(\le\) |(Statement.\(\zeta\) L N).im| \(\wedge\)
      0 < -(Statement.\(\zeta\) L N).re \(\wedge\) -(Statement.\(\zeta\) L N).re \(\le\) 1600/L^3 \(\wedge\)
      |(Statement.h L N)^3*|(Statement.\(\zeta\) L N).im|-
        (420*|Cert.wStar|-210*Real.pi^2*deriv (deriv Cert.optical) Cert.thetaStar/(N:\(\mathbb{R}\))^2+
          210*Real.pi^2*deriv (deriv Cert.optical) Cert.thetaStar*Cert.betaStar.im/(N:\(\mathbb{R}\))^3)|
          \(\le\) 1260000/(N:\(\mathbb{R}\))^4 \(\wedge\)
      |L^3*(Statement.\(\zeta\) L N).re+
        210*Real.pi^2*deriv (deriv Cert.optical) Cert.thetaStar*Cert.betaStar.re|
          \(\le\) 1260000/(N:\(\mathbb{R}\))
\end{alltt}}
\noindent
\emph{Translation.} For every $L>0$ and every integer $N\ge512$, with
$h=L/N$, there is $\varphi\in X_h$, $\varphi\ne0$, such that
$\mathring A_h\varphi=\zeta_h\varphi$ and \eqref{eq:bounds}--\eqref{eq:re}
hold. All conclusions hold for the same
$\zeta_h$ and the same $\varphi$. The eigenvalue $\zeta_h$ is the explicit
number \texttt{Statement.$\zeta$ L N}. The eigenfunction $\varphi$ is only
asserted to exist.

The names in the statement have the following meaning
(Table~\ref{tab:names}). Each line is fixed by a definition or a proved
Lean statement, not by the name.

{\small
\begin{longtable}{@{}>{\raggedright\arraybackslash}p{.08\textwidth}>{\raggedright\arraybackslash}p{.4\textwidth}>{\raggedright\arraybackslash}p{.46\textwidth}@{}}
\caption{The names in the Lean statement of Theorem~\ref{thm:top}.}\label{tab:names}\\
\toprule
Paper & Lean & Fixed by\\\midrule
$h$, $X_h$, $\mathring A_h$ & \texttt{Statement.h L N}, \texttt{Statement.X\_h}, \texttt{Statement.A\_h} & $h=L/N$; \lean{SLean.Statement.space_characterization} and \lean{SLean.Statement.weak_form} identify $X_h$ and \eqref{eq:Sscheme}\\
$\zeta_h$ & \texttt{Statement.$\zeta$ L N} & $\zeta_h=420i\tilde\varpi(\varsigma_N)/h^3$ with $\varsigma_N=(\pi+y_N)/N$ (Section~\ref{sec:A2})\\
$\varpi$ & \lean{Cert.optical} & the formula $(-b-\sqrt\Delta)/(2a)$ of Section~\ref{sec:disp}\\
$\xi_*$ & \lean{Cert.thetaStar} & $\arccos c_*$ with $\tilde H(c_*)=0$, $c_*\in[c_{\rm lo},c_{\rm hi}]$\\
$\Omega_*$ & \texttt{420*|Cert.wStar|} & \lean{optical_thetaStar}: $w_*=\varpi(\xi_*)$\\
$\varpi''$ & \texttt{deriv (deriv Cert.optical) Cert.thetaStar} & $\varpi$ is analytic near $\xi_*$, so this is $\varpi''(\xi_*)$\\
$\beta$ & \lean{Cert.betaStar} & \lean{betaStar_derivative_formula}: the definition \eqref{eq:beta}\\
$\mathcal C_S$ & written out as $210\pi^2\varpi''\Rea\beta$ & the definition \eqref{eq:beta}\\
\bottomrule
\end{longtable}}

\paragraph{Zero group velocity.}
The four statements in \eqref{eq:zgv} are
{\small
\begin{alltt}
theorem thetaStar_range : 0 < thetaStar \(\wedge\) thetaStar < Real.pi
theorem optical_thetaStar : optical thetaStar = wStar
theorem optical_stationary : deriv optical thetaStar = 0
theorem optical_curvature_pos : 0 < deriv (deriv optical) thetaStar
\end{alltt}}
\noindent
in \lean{SLean.Cert}. \emph{Translation:} $0<\xi_*<\pi$,
$\varpi(\xi_*)=w_*$, $\varpi'(\xi_*)=0$ and $\varpi''(\xi_*)>0$. Since
$\varpi$ is analytic near $\xi_*$, these derivatives are the true
derivatives.

\subsection{The files}\label{sec:A4}
\paragraph{Scope.}
This appendix concerns the Lean proof of Theorem~\ref{thm:top}. The
results of this proof that the paper uses are Theorem~\ref{thm:top},
\eqref{eq:zgv} and the lemmas of Section~\ref{sec:leanlemmas}.
Table~\ref{tab:corr} gives their Lean names. The conjugate eigenpair
stated after Theorem~\ref{thm:top} is derived by hand from the real
coefficients. Its coordinate form is \lean{SLean.Bridge.generator_star}.
The Lean proofs of the other results of Sections~\ref{sec:spurious}
and~\ref{sec:gap} (Remark~\ref{rem:leanscope}) are not described in
Sections~\ref{sec:A1}--\ref{sec:A4}; see Section~\ref{sec:A5}.
The Lean
development also contains statements on the failure of resolvent,
Crandall--Pazy and Kato bounds, in \lean{SLean.Statement} and
\lean{SLean.CoreStatement}. The Crandall--Pazy and Kato arguments are
used for (H1) in Section~\ref{sec:A5}. The five files in \lean{SLean.Explicit.Sharp.Expansion} are
included in the build target \texttt{SLean} and in the ancillary source
package.

\paragraph{Toolchain and build.}
The development uses Lean~4.35.0 and the mathlib commit
\begin{quote}\small\ttfamily
56a005855ea17e0b0958f73ec801ae140b0a1679
\end{quote}
From the project root, run \texttt{lake exe cache get}, followed by
\texttt{lake build SLean}. This target includes all five files in
\lean{SLean.Explicit.Sharp.Expansion}. The Lean statements of
Theorem~\ref{thm:top} are
\lean{SLean.Explicit.Sharp.Expansion.spectral_top_with_expansion}. On
4~November 2026, the incremental build completed in $21.2$ seconds. The
separate command \texttt{lake env lean referee/CheckExpansion.lean}
checked the axiom dependencies of this statement, of the original
spectral theorem and of the four statements in \eqref{eq:zgv}, in
$19.4$ seconds. All six use only \texttt{propext},
\texttt{Classical.choice} and \texttt{Quot.sound}. A clean build
measured on 1~November 2026, before the Expansion files were added, took
$49$ minutes $46$ seconds. Retrieving the dependencies and preparing the
cache took a further $8$ minutes $6$ seconds, on Windows with an Intel
Core i7-1165G7 processor. This earlier measurement concerns the
preceding $169$-module target. A clean-build time for the enlarged
target has not been measured. The ancillary files of arXiv:2611.12585
contain the sources,
the pinned configuration, the README, the counting script and the
validation logs.

\paragraph{Size.}
Counted over the full transitive import closure of the default entry
\lean{SLean}, including the modules on semigroup consequences and
excluding mathlib and other dependencies, the development contains
$174$ project modules, $983$ theorem and lemma declarations and
$35{,}912$ physical source lines. The \lean{SLean.Explicit} portion
accounts for $88$ modules, $388$ declarations and $23{,}701$ lines.
Within it, the five Expansion files contribute $17$ declarations and
$423$ lines. Module counts include entry and audit modules. Declaration
counts include private theorems and lemmas, but exclude definitions and
audit commands. Line counts include comments and blank lines. These are
source-size counts, not counts of independent mathematical results. The
counting script \texttt{referee/count\_sources.py}, the exact module
lists and the per-file counts are included in the ancillary files.

\paragraph{Iterations.}
Three technical changes occurred during the iterations of
Remark~\ref{rem:history}. Elimination in $c=\cos\xi$ replaced a
high-precision bisection for $\xi_*$. Explicit discs for the root branches
replaced Rouch\'e's theorem. The cubic bound for the reflection logarithm
was sharpened from $2200$ to $340$ for the damping certificate.

\paragraph{Notation.}
The Lean code keeps the names of an earlier draft. The main differences
are as follows.
\begin{itemize}
\item $\varsigma$ is written $\varepsilon$ in comments, and \texttt{z} or
\texttt{xi} in code.
\item The slope ratio $\mathcal V$ is \texttt{W}. It is not the weighted
scheme.
\item The column scalings $\gamma_k$ are \texttt{lam}.
\item $\mathfrak r$ is \texttt{reflection}, and $\ell$ is
\texttt{logReflection}.
\item $\mathcal T_N$ is \texttt{quantMap}, and $\tilde H$ is
\texttt{criticalPolynomial}.
\item The eigenpair is indexed by $N$.
\end{itemize}

\paragraph{Correspondence.}
Table~\ref{tab:corr} lists the Lean declarations for each step of
Section~\ref{sec:A1} and each lemma of Section~\ref{sec:leanlemmas}.

\input{correspondence_table}

\subsection{The Lean statements of the other results}\label{sec:A5}
\input{appA5_statements}

%% file: correspondence_table.tex
\begingroup\small
Prefixes: S = \nolinkurl{SLean.Explicit.Sharp}, Ex = \nolinkurl{SLean.Explicit.Sharp.Expansion},
C = \nolinkurl{SLean.Cert}, A = \nolinkurl{SLean.Algebra}, B = \nolinkurl{SLean.Bridge},
DR = \nolinkurl{SLean.Pilots.DoubleRoot}, H = \nolinkurl{SLean.Hermite},
E = \nolinkurl{SLean.Explicit}, An = \nolinkurl{SLean.Analytic}.
\begin{longtable}{@{}>{\raggedright\arraybackslash}p{.25\textwidth}>{\raggedright\arraybackslash}p{.71\textwidth}@{}}
\caption{Lean declarations for Theorem~\ref{thm:top}.}\label{tab:corr}\\
\toprule Statement & Lean declarations\\\midrule\endfirsthead
\toprule Statement & Lean declarations (continued)\\\midrule\endhead
\bottomrule\endfoot
Theorem~\ref{thm:top} & \nolinkurl{Ex.spectral_top_with_expansion}; \linebreak[1] \nolinkurl{SLean.Statement.spectral_top}; \linebreak[1] \nolinkurl{S.spectral_top_explicit_512}\\[3pt]
\eqref{eq:zgv} & \nolinkurl{C.thetaStar_range}; \linebreak[1] \nolinkurl{C.optical_thetaStar}; \linebreak[1] \nolinkurl{C.optical_stationary}; \linebreak[1] \nolinkurl{C.optical_curvature_pos}\\[3pt]
Step 1, $\xi_*$ & \nolinkurl{C.critical_endpoint_signs}; \linebreak[1] \nolinkurl{C.exists_critical_c}; \linebreak[1] \nolinkurl{C.stationary_elimination}; \linebreak[1] \nolinkurl{C.stationary_relation}; \linebreak[1] \nolinkurl{C.optical_second_deriv}\\[3pt]
Lemma~\ref{lem:cert} & \nolinkurl{C.certificate_C1}; \linebreak[1] \nolinkurl{C.certificate_C2}; \linebreak[1] \nolinkurl{C.certificate_C3}; \linebreak[1] \nolinkurl{C.certificate_C4}; \linebreak[1] \nolinkurl{C.certificate_C5}; \linebreak[1] \nolinkurl{C.wStar_bounds}; \linebreak[1] \nolinkurl{C.optical_curvature_bounds}; \linebreak[1] \nolinkurl{C.boundary_data_bounds}; \linebreak[1] \nolinkurl{C.t1_constant_margins}\\[3pt]
$\beta$ in \eqref{eq:beta} & \nolinkurl{C.betaStar_derivative_formula}\\[3pt]
Lemma~\ref{lem:det} & \nolinkurl{A.symbol_kernel}; \linebreak[1] \nolinkurl{A.exact_mode_value}; \linebreak[1] \nolinkurl{A.exact_mode_slope}; \linebreak[1] \nolinkurl{A.mode_left}; \linebreak[1] \nolinkurl{A.Sexp}; \linebreak[1] \nolinkurl{A.Sdet}; \linebreak[1] \nolinkurl{A.Sdecomp}; \linebreak[1] \nolinkurl{B.matrix_eigenpair_function_equivalence}\\[3pt]
Step 4, normal form & \nolinkurl{DR.exists_normal_coordinate}; \linebreak[1] \nolinkurl{DR.exists_companion_germ}\\[3pt]
Lemma~\ref{lem:domain} & \nolinkurl{S.F_analytic_wide}; \linebreak[1] \nolinkurl{S.eta_parameter_small}; \linebreak[1] \nolinkurl{S.roots_analytic_disk}; \linebreak[1] \nolinkurl{S.roots_factor}; \linebreak[1] \nolinkurl{S.sigma_bounds}; \linebreak[1] \nolinkurl{S.sigma_disk}; \linebreak[1] \nolinkurl{E.sigma_agrees}; \linebreak[1] \nolinkurl{E.sigma_zero}; \linebreak[1] \nolinkurl{E.sigma_deriv}; \linebreak[1] \nolinkurl{An.sigma_spec}; \linebreak[1] \nolinkurl{An.bMinus_eq_eventually}; \linebreak[1] \nolinkurl{An.wMinus_eq_eventually}; \linebreak[1] \nolinkurl{E.roots_agree}; \linebreak[1] \nolinkurl{E.roots_zero}; \linebreak[1] \nolinkurl{An.zInnerBranch_zero}; \linebreak[1] \nolinkurl{An.zOuterBranch_zero}; \linebreak[1] \nolinkurl{S.minus_eq_zPlus_sigma}; \linebreak[1] \nolinkurl{S.roots_injective}; \linebreak[1] \nolinkurl{S.roots_work}; \linebreak[1] \nolinkurl{S.boundary_work_bounds}; \linebreak[1] \nolinkurl{S.tails_work}; \linebreak[1] \nolinkurl{An.analytic_implicit_root}; \linebreak[1] \nolinkurl{E.inner_analytic}; \linebreak[1] \nolinkurl{E.outer_analytic}; \linebreak[1] \nolinkurl{E.minus_agrees}; \linebreak[1] \nolinkurl{S.optical_domain}; \linebreak[1] \nolinkurl{S.productRatio_bound}; \linebreak[1] \nolinkurl{S.phaseCorrection_bound}\\[3pt]
Lemma~\ref{lem:refl} & \nolinkurl{E.reflection_zero}; \linebreak[1] \nolinkurl{S.reflection_eq_exp}; \linebreak[1] \nolinkurl{S.logReflection_analytic}; \linebreak[1] \nolinkurl{S.logReflection_zero}; \linebreak[1] \nolinkurl{E.reflection_deriv}; \linebreak[1] \nolinkurl{S.logReflection_deriv}; \linebreak[1] \nolinkurl{S.logReflection_bound}; \linebreak[1] \nolinkurl{S.logReflection_second_tight}; \linebreak[1] \nolinkurl{S.logReflection_third_tight}; \linebreak[1] \nolinkurl{S.varpi_work}\\[3pt]
Lemma~\ref{lem:tail}, \eqref{eq:Phi} & \nolinkurl{S.normalized_decomposition}; \linebreak[1] \nolinkurl{S.quant_norm}; \linebreak[1] \nolinkurl{S.quant_in_work}; \linebreak[1] \nolinkurl{S.quant_powers}; \linebreak[1] \nolinkurl{S.normalizedError_bound}; \linebreak[1] \nolinkurl{S.geometric_sixth}; \linebreak[1] \nolinkurl{S.quant_error}; \linebreak[1] \nolinkurl{S.quant_reflection_inv}; \linebreak[1] \nolinkurl{S.errorLog_bound}\\[3pt]
Lemma~\ref{lem:fixed} & \nolinkurl{S.analytic_fixed_point}; \linebreak[1] \nolinkurl{S.quantMap_analytic}; \linebreak[1] \nolinkurl{S.quantMap_bound}; \linebreak[1] \nolinkurl{S.quantMap_fixed_zero}; \linebreak[1] \nolinkurl{S.quantized_zero_exists}; \linebreak[1] \nolinkurl{S.ySequence_spec}\\[3pt]
Lemma~\ref{lem:box} & \nolinkurl{S.fixed_scaled_bound}; \linebreak[1] \nolinkurl{S.fixed_bootstrap}; \linebreak[1] \nolinkurl{S.scaled_X_box_512}; \linebreak[1] \nolinkurl{S.scaled_V_bootstrap}; \linebreak[1] \nolinkurl{S.fixed_real_upper}\\[3pt]
Taylor remainders (Steps 6, 7) & \nolinkurl{S.sigma_fourth_quotient}; \linebreak[1] \nolinkurl{S.sigma_fourth_im}; \linebreak[1] \nolinkurl{S.F_fourth_quotient_tight}; \linebreak[1] \nolinkurl{S.F_fourth_im_tight}; \linebreak[1] \nolinkurl{S.power_remainder_im}; \linebreak[1] \nolinkurl{S.companion_bound}; \linebreak[1] \nolinkurl{S.companion_negative}; \linebreak[1] \nolinkurl{S.beta_box}; \linebreak[1] \nolinkurl{S.beta_bound_tight}\\[3pt]
Step 6, bounds \eqref{eq:bounds} & \nolinkurl{S.scaled_quant_identity}; \linebreak[1] \nolinkurl{S.scaled_optical_identity}; \linebreak[1] \nolinkurl{S.scaled_quant_box_512}; \linebreak[1] \nolinkurl{S.scaledP_box}; \linebreak[1] \nolinkurl{S.scaled_quant_error_tight}; \linebreak[1] \nolinkurl{S.scaled_optical_error_tight}; \linebreak[1] \nolinkurl{S.horner_upper_512}; \linebreak[1] \nolinkurl{S.scaled_upper_imp}; \linebreak[1] \nolinkurl{S.damping_bound_512}; \linebreak[1] \nolinkurl{S.frequency_bound}; \linebreak[1] \nolinkurl{S.parameter_bounds_512}\\[3pt]
Step 7, \eqref{eq:im}--\eqref{eq:re} & \nolinkurl{Ex.fixed_point_leading_error}; \linebreak[1] \nolinkurl{Ex.optical_model_eq}; \linebreak[1] \nolinkurl{Ex.sequence_optical_expansion}; \linebreak[1] \nolinkurl{Ex.physical_expansions}; \linebreak[1] \nolinkurl{S.optical_third_factor}\\[3pt]
\eqref{eq:Hermite}, \eqref{eq:element} & \nolinkurl{H.nodal_values}; \linebreak[1] \nolinkurl{H.reference_mass}; \linebreak[1] \nolinkurl{H.physical_mass}; \linebreak[1] \nolinkurl{H.physical_stiffness}; \linebreak[1] \nolinkurl{H.physicalM_eq}; \linebreak[1] \nolinkurl{H.physicalK_eq}\\[3pt]
\eqref{eq:bloch}, \eqref{eq:dispersion} & \nolinkurl{A.symbolK_exp}; \linebreak[1] \nolinkurl{A.symbolM_exp}; \linebreak[1] \nolinkurl{A.bloch_dispersion}; \linebreak[1] \nolinkurl{A.q_exp_dispersion}\\[3pt]
\eqref{eq:Sscheme}, \eqref{eq:inverse}, $X_h$ & \nolinkurl{SLean.Statement.space_characterization}; \linebreak[1] \nolinkurl{SLean.Statement.norm_identification}; \linebreak[1] \nolinkurl{SLean.Statement.weak_form}; \linebreak[1] \nolinkurl{B.inverse_inequality_squared}\\[3pt]
\end{longtable}\endgroup

%% file: appA5_statements.tex

Theorem~\ref{thm:top}, Corollary~\ref{cor:failures} and the main
conclusions of Sections~\ref{sec:eigfun}--\ref{sec:lineargap} and
Section~\ref{sec:gap} are collected in one Lean theorem,
\lean{CompanionFinal.main_theorem}. It has a Lean~4 proof, checked by the
Lean kernel. Its only premise is (P1). Proposition~\ref{prop:necessary} is
a separate Lean theorem without premise, \lean{CompanionFinal.prop318}.
Further results are internal Lean theorems outside the main theorem, and
some parts are not in Lean. Table~\ref{A5:corr} lists, for each numbered
result and remark of Sections~\ref{sec:eigfun}--\ref{sec:lineargap} and
Section~\ref{sec:gap}, where it stands. Unprefixed Lean names are in
\lean{CompanionPilot}.

The translation below was produced by the procedure of
Remark~\ref{rem:history}, with Claude Fable~5.2 in place of Claude
Fable~5.1. On this basis, the author made sure that the translated
premise is \cite[(3.15)]{Main} for uniform meshes, and that each
translated tail is the corresponding result of the paper or a stronger
one, with $\mathcal U$ in place of $u$ in Theorem~\ref{thm:total}(b)(c)
and Corollary~\ref{cor:noproduct}.

\paragraph{The Lean statements.}
{\small
\begin{alltt}
def P1 : Prop :=
  \(\exists\) C : \(\mathbb{R}\), 0 < C \(\wedge\) \(\forall\) (N : \(\mathbb{N}\)) (hN : 0 < N) (t : \(\mathbb{R}\)), 0 < t \(\to\) \(\forall\) x \(\in\) Vh hN,
    normH (AW hN (TW hN t x)) \(\le\)
      C * (min 1 t) ^ (-(3/2 : \(\mathbb{R}\))) * Real.exp (-(3 * Real.pi ^ 2 / 4) * t) * normH x
def Tails : Prop :=
  H0 \(\wedge\) CompanionA.H1 \(\wedge\) H2 \(\wedge\) CompanionE.H3a \(\wedge\) CompanionE.H3c \(\wedge\) CompanionE.H4 \(\wedge\)
    H6i \(\wedge\) H6ii \(\wedge\) CompanionB.H6iiiV \(\wedge\) H3bH5V \(\wedge\) FamilyProps
theorem main_theorem : P1 \(\to\) Tails
theorem prop318 : CompanionH.Abstract.Prop41.\{u\}
\end{alltt}}
\noindent
In the displayed Lean text, names without a prefix are in
\lean{CompanionPilot} or \lean{CompanionFinal}.
The eleven conjuncts of \lean{CompanionFinal.Tails} are called the tails.
We translate them below as (H0), (H1), (H2), (H3a), (H3c), (H4), (H6i),
(H6ii), (H6iii), (B) and (F), in this order.
\emph{Translation of the main theorem:} if (P1) holds, then all eleven
tails hold. The tails (H0), (H1), (H6i), (H6ii), (H6iii) and (F) are
proved without (P1) (\lean{CompanionFinal.tails_without_P1}).

\paragraph{The premise.}
\begin{enumerate}
\item[(P1)] There is $C>0$ such that, for every integer $N\ge1$, $h=1/N$,
every $t>0$ and every $x\in X_h^{\C}$,
\[
 \norm{A_hT_h(t)x}_H\le C\,t_\wedge^{-3/2}e^{-3\pi^2t/4}\norm x_H,
 \qquad t_\wedge=\min\{1,t\}.
\]
\end{enumerate}
Here $X_h^{\C}$, $A_h$, $T_h$ and $\norm\cdot_H$ are those of the
W-scheme below. One constant $C$ serves all $N$, $t$ and $x$. This is
\cite[(3.15)]{Main} with $\alpha=1$, $L=1$, uniform meshes, $B_h=-A_h$
and $\lambda=3\pi^2/4$.

\paragraph{Conventions.}
\begin{itemize}
\item $L^2=L^2(0,1)$ is the complex space \lean{L2}, for Lebesgue measure
$\mu$ on $(0,1]$. Its elements are classes of functions modulo null sets.
We write $\norm f_2$ and $(f,g)_2=\int_0^1f\overline g\,dx$.
$\psi_n\rightharpoonup0$ means $(\psi_n,g)_2\to0$ for every $g\in L^2$
(\lean{WeakNull}; Mathlib writes $\langle g,\psi_n\rangle$).
\item $(f,g)_H=\int_0^1f\overline g\,(1+x)\,dx$ (\lean{ipH}) and
$\norm f_H=(f,f)_H^{1/2}$ (\lean{normH}). Thus
$\norm f_2\le\norm f_H\le\sqrt2\norm f_2$.
\item $f\in L^2$ is real if $\Ima f=0$ almost everywhere
(\lean{IsRealL2}). The bounded operators below are complex linear. An
operator norm $\norm\cdot$ is that of a bounded operator on $L^2$ with
$\norm\cdot_2$, unless stated otherwise.
\item $N\in\mathbb N$ and $h=1/N$ (\lean{hm}), except in (H0) and (H1),
where $L>0$ and $h=L/N$. ``For all large $N$'' means ``there is $N_0$
such that for all $N\ge N_0$''. Limits as $h\to0$ are limits as
$N\to\infty$. A clause ``for all $t\in[a,b]$'' is vacuous if $b<a$.
\item Integrals in time, $\int_a^b$ with $a\le b$, are over $(a,b]$, for
Lebesgue measure. For $L^2$-valued functions they are Bochner integrals.
In Lean the integral of a non-integrable function is $0$. All integrals
in the tails are of integrable functions; where this is not obvious, we
say why.
\item Some objects are defined by a choice (\lean{Classical.choose}), with
the value $0$ if nothing exists to choose. Lean's supremum of a set of
reals that is not bounded above is $0$, and $0^{-1}=0$. For each such
default we name the statement that excludes it.
\item Lean's derivative \lean{deriv} is $0$ where a function is not
differentiable. Below it is applied to differentiable functions, except
for $w''$ of a spline at the nodes, a null set, and for $\varpi''(\xi_*)$
in (H0) and $\mathcal C_S$, where $\varpi$ is analytic
(Table~\ref{tab:names}).
\end{itemize}

\paragraph{The trial space and the S-scheme.}
\begin{itemize}
\item A coefficient vector $v$ (\lean{SLean.Vector}) lists the values
$u_1,\dots,u_{N-1}$ and the scaled slopes $d_0,\dots,d_{N-1}$; put
$u_0=u_N=d_N=0$. Its spline (\lean{SLean.Bridge.spline}) is the $C^1$
function that on $[x_j,x_{j+1}]$, $x_j=jh$, equals the cubic of
\eqref{eq:Hermite} with the nodal data $(u_j,d_j,u_{j+1},d_{j+1})$. Thus
$d_j=h\,v'(x_j)$.
\item $X_h^{\C}\subset L^2$ (\lean{Vh}) is the set of classes of these
splines. By \lean{SLean.Statement.space_characterization}, it is the set
of $C^1$ functions on $[0,1]$ that are cubic on each cell and satisfy
$\chi(0)=\chi(1)=\chi'(1)=0$. For $\chi\in X_h^{\C}$ we use the
continuous representative, and $\chi'$ is its derivative. $X_h$ is the set
of real elements of $X_h^{\C}$.
\item $\mathring A_h$ on $X_h^{\C}$ (\lean{AS}) is
\lean{SLean.Statement.A_h}, transported to $X_h^{\C}$. It satisfies
\eqref{eq:Sscheme} for $v,\chi\in X_h^{\C}$
(\lean{SLean.Statement.weak_form}). $\mathring T_h(t)=e^{t\mathring A_h}$
(\lean{TS}) is the operator exponential, in Hermite coordinates.
$\mathring\Pi_h$ is the $L^2$-orthogonal projection onto $X_h^{\C}$
(\lean{PiS}). $\mathring T_h(t)\mathring\Pi_h$ (\lean{TSPi}) and
$\mathring A_h\mathring\Pi_h$ (\lean{CompanionC.ASL}) are bounded
operators on $L^2$.
\item For $q\in L^2$, $\mathcal G_h[q]\in X_h^{\C}$ is defined by
$(\mathcal G_h[q],\chi)_2=\tfrac12\int_0^1q\,\overline{\chi'}$ for every
$\chi\in X_h^{\C}$ (\lean{CompanionC.Gh}). This is $\mathcal G_h$ of
Section~\ref{sec:setting}, for $q\in L^2$.
\item $\mathcal C_S$ (\lean{C_S}) is written with the Lean terms of
\eqref{eq:re} in (H0), so it is the constant $\mathcal C_S$ of
\eqref{eq:beta} (Table~\ref{tab:names}).
\end{itemize}

\paragraph{The W-scheme.}
\begin{itemize}
\item $b(w,\psi)=-\int_0^1w''\,\overline{\psi+(1+x)\psi'}\,dx$
(\lean{bForm}). For a spline, $w''$ exists except at the nodes, a null set.
\item $A_h$ on $X_h^{\C}$ (\lean{AHlin}) is defined through an
$H$-orthonormal basis of $X_h^{\C}$, so that $(A_hw,\psi)_H=-b(w,\psi)$
for $w,\psi\in X_h^{\C}$ (\lean{inner_AHlin}). This is \eqref{eq:Wform}.
$A_h$ is injective (\lean{probe_AW_injective}).
\item $\Pi_h$ is the $H$-orthogonal projection onto $X_h^{\C}$, as an
operator on $L^2$ (\lean{PiW}). $A_h\Pi_h$ (\lean{AW}) is a bounded
operator on $L^2$. $T_h(t)=\exp(tA_h\Pi_h)$ on $L^2$ (\lean{TW}); on
$X_h^{\C}$ it is $e^{tA_h}$, and on the $H$-orthogonal complement of
$X_h^{\C}$ it is the identity. $T_h(t)\Pi_h$ is \lean{UW}, and
$A_h^{-1}\Pi_h$ is \lean{JW}. As sequences in $N$, \lean{JWt} and
\lean{UWt} take these values for $N\ge1$ and the value $0$ at $N=0$.
\end{itemize}

\paragraph{The limit operators.}
\begin{itemize}
\item $A^{-1}$ (\lean{Ainv}) is the limit of $A_h^{-1}\Pi_h$ in operator
norm if it exists, and $0$ otherwise. By (H6i),
$\norm{A_h^{-1}\Pi_h-A^{-1}}\le2h$, so the limit exists.
\item For $t\in\R$, $T(t)$ (\lean{Tsg}) is a bounded operator $B$ with
$T_h(t)\Pi_hf\to Bf$ in $L^2$ for every $f$, if such $B$ exists, and $0$
otherwise. Such $B$ is unique. (H6ii) asserts that it exists for $t\ge0$.
The tails use $T(t)$ only for $t\ge0$.
\item The identification of these limits with $A^{-1}$ and $T$ of
Section~\ref{sec:setting} is explained after
Proposition~\ref{prop:galerkin}.
\item $\mathring L_h(t)=\mathring T_h(t)\mathring\Pi_h-T(t)$
(\lean{CompanionE.LS}).
\item For $f\in L^2$, a triple $(\varrho,\varrho_1,\varrho_2)$ of
functions $\R\to\C$ solves the boundary value problem
(\lean{SolvesBVP}) if, on $[0,1]$, $\varrho$ has the derivative
$\varrho_1$ and $\varrho_1$ has the derivative $\varrho_2$ (one-sided at
the ends), $\varrho_2(x)=\varrho_2(0)-\int_0^xf\,d\mu$ for $x\in[0,1]$, and
$\varrho(0)=\varrho(1)=\varrho_1(1)=0$. Thus $\varrho\in C^2[0,1]$,
$\varrho''$ is absolutely continuous and $-\varrho'''=f$ almost
everywhere. Nothing is required outside $[0,1]$.
\end{itemize}

\paragraph{Trajectories.}
\begin{itemize}
\item $w\in L^2$ lies in $H_0^1$ (\lean{CompanionB.IsH01}) if there is
$g\in L^2$ with $\int_0^1g=0$ and $w(x)=\int_0^xg$ for almost every $x$.
Then $g$ is unique, and $\partial_xw=g$ (\lean{CompanionB.dx}). For
$w\notin H_0^1$, $\partial_xw=0$ by default. $|w|_1=\norm{\partial_xw}_2$.
\item $w\,\partial_xw$ (\lean{CompanionB.wdxw}) is the class of
$x\mapsto(\int_0^x\partial_xw)\,\partial_xw(x)$. For $w\in H_0^1$ it
equals $w\cdot\partial_xw$ almost everywhere
(\lean{CompanionB.coeFn_wdxw}). For $w\notin H_0^1$ it is $0$.
$X_h^{\C}\subset H_0^1$, with $\partial_x\chi=\chi'$
(\lean{CompanionB.isH01_of_mem}, \lean{CompanionB.dx_of_mem}).
\item A trajectory is a map $w:\R\to L^2$, used for $t\ge0$. The time
measure is Lebesgue measure on $(0,\infty)$.
\item $w\in\mathcal X_0(\nu)$ (\lean{CompanionB.InX0}) if $w$ is
continuous on $[0,\infty)$ into $L^2$; $e^{\nu t}\norm{w(t)}_H\to0$ as
$t\to\infty$; $w(t)\in H_0^1$ for almost every $t>0$;
$t\mapsto\partial_xw(t)$ is almost everywhere strongly measurable on
$(0,\infty)$; and $t\mapsto e^{2\nu t}|w(t)|_1^2$ is integrable on
$(0,\infty)$.
\item $\mathcal E_\nu(w)=\sup_{t\ge0}e^{\nu t}\norm{w(t)}_H
+\bigl(\int_0^\infty e^{2\nu t}|w(t)|_1^2\,dt\bigr)^{1/2}$
(\lean{CompanionB.energy}). For $w\in\mathcal X_0(\nu)$ neither default
occurs: the supremum is over a bounded set (\lean{CompanionB.bdd_of_InX0}),
and the integrand is integrable. $\mathcal X_0(\nu)$ is closed under
linear combinations (\lean{CompanionB.InX0_add},
\lean{CompanionB.InX0_smul}). The tails apply $\mathcal E_\nu$ only to
elements of $\mathcal X_0(\nu)$ asserted in the same statement and to
differences of two of them.
\end{itemize}

\paragraph{The nonlinear schemes.}
\begin{itemize}
\item $u$ solves the W-scheme with datum $f$ (\lean{CompanionB.IsWSol})
if $u(0)=\Pi_hf$ and, for every $t\ge0$, $u(t)\in X_h^{\C}$ and
$u'(t)=A_hu(t)-\Pi_h\bigl(u(t)\,\partial_xu(t)\bigr)$. The derivative is
in $L^2$, within $[0,\infty)$, hence one-sided at $t=0$. Nothing is
required for $t<0$. This is \eqref{eq:Wscheme} with
$F(u,u)=-u\,\partial_xu$. $u_h[f]$ (\lean{CompanionB.uW}) is a chosen
solution, $0$ by default; for $N=0$, $u_h[f]=0$ (\lean{CompanionB.uWt}).
\item $\mathcal U(f)(t)$ (\lean{CompanionB.Ucal}) is the limit of
$u_h[f](t)$ in $L^2$ as $N\to\infty$ if it exists, and $0$ otherwise. It
does not depend on $\nu$.
\item $u$ solves the S-scheme with datum $u_0$ (\lean{CompanionE.IsSSol})
if $u(0)=\mathring\Pi_hu_0$ and, for every $t\ge0$, $u(t)\in X_h^{\C}$
and $u'(t)=\mathring A_hu(t)-\mathring\Pi_h\bigl(u(t)\,\partial_xu(t)\bigr)$,
in the same sense. For $v\in X_h^{\C}$,
$-\mathring\Pi_h(v\,\partial_xv)=\mathcal G_h[v^2]$
(\lean{CompanionFinal.unified_objects}); so this is the S-scheme of
Section~\ref{sec:setting}. $\mathring u_h[u_0]$ (\lean{CompanionE.uS}) is a
chosen solution, $0$ by default.
\item For $\varepsilon\in\R$ and $t\ge0$,
$\mathring\Cd_h(\varepsilon,t)=\bigl(\mathcal U(\varepsilon f_h)(t)-\varepsilon T(t)f_h\bigr)
-\bigl(\mathring u_h[\varepsilon f_h](t)-\varepsilon\mathring T_h(t)\mathring\Pi_hf_h\bigr)$
(\lean{CompanionH.CS}), with $f_h$ as follows.
\end{itemize}

\paragraph{The eigenfunction family.}
Let $N\ge512$ and $h=1/N$. Lean fixes $y_N$ by a choice
(\lean{SLean.Explicit.Sharp.ySequence}) from
\lean{SLean.Explicit.Sharp.quantized_zero_exists}, with $|y_N|<1/2$,
$\mathcal T_N(y_N)=y_N$, $\Phi_N(\varsigma_N)=0$ and $|y_N|\le147/N$,
where $\varsigma_N=(\pi+y_N)/N$ (Lemma~\ref{lem:fixed}). Put
$\eta_N=i\tilde\varpi(\varsigma_N)$. Then
$\zeta_h=420\eta_N/h^3$ is \texttt{Statement.$\zeta$ 1 N}, the eigenvalue
of (H0) with $L=1$. Let $(z_1,z_2,z_3,z_4)=(z_+,z_-,z_3,z_4)$ be the four
roots of Lemma~\ref{lem:domain} at $\varsigma_N$
(\lean{SLean.Explicit.roots}), and
$(\gamma_1,\gamma_2,\gamma_3,\gamma_4)=(1,1,1,z_4^{-N})$. The $k$th column
of $\D_N=\D_N(\varsigma_N)$ is
$\gamma_k\bigl(1,B(z_k,\eta_N),z_k^N,z_k^N\mathcal V(z_k,\eta_N)\bigr)^{\mathsf T}$,
with $\mathcal V$ and $B$ of \eqref{eq:slope} and \eqref{eq:B}
(\lean{SLean.boundaryMatrix}). Lean chooses a nonzero
$\alpha_N\in\ker\D_N$ (\lean{CompanionA.alphaN}). Let $v_N\in X_h^{\C}$
have the nodal values and scaled slopes
\[
 v_N(x_j)=\sum_{k=1}^4\alpha_{N,k}\gamma_kz_k^j,\qquad
 h\,v_N'(x_j)=\sum_{k=1}^4\alpha_{N,k}\gamma_kz_k^j\,\mathcal V(z_k,\eta_N)
 \qquad(0\le j\le N)
\]
(\lean{CompanionA.xN}, \lean{CompanionA.xN_spec}). The general form is
\lean{CompanionA.nodal_eigenvector}, the nodal formula of
Lemma~\ref{lem:det} with its nonvanishing: for $N\ge3$, distinct nonzero
zeros $z_k$ of $q(\cdot;\eta)$ with $\mathcal D(z_k,\eta)\ne0$,
$\gamma_k\ne0$ and $0\ne\alpha\in\ker\D$, these nodal data give a nonzero
coefficient vector $x$ with $Kx=\eta Mx$ for the assembled matrices of
\eqref{eq:element}. Then
$\varphi_h=v_N/\norm{v_N}_2$ (\lean{CompanionA.phiA}),
$\tilde f_h=\sqrt2\,\Rea\varphi_h$ (\lean{CompanionA.ftilA}) and
$f_h=\norm{\tilde f_h}_2^{-1}\tilde f_h$ (\lean{CompanionA.fA}). The
denominators are nonzero by \lean{CompanionA.kernel_exists},
\lean{CompanionA.zN_data}, \lean{CompanionA.xN_spec} and
\lean{CompanionA.ftilA_ne_zero}. For $N<512$ all three are $0$. So
$\norm{f_h}_2\le1$ for every $N$, with equality for $N\ge512$. The tails
concern this one family. Its phase is fixed by the choices of $y_N$ and
$\alpha_N$.

\paragraph{The quadratic coefficients.}
These objects enter (H5). $\mathsf m(r)=r-3r^2+2r^3$ and
$\mathsf m'(r)=1-6r+6r^2$, as in Section~\ref{sec:eigfun}
(\lean{Companion.m}, \lean{Companion.mp}).
$p(r)=H_0(r)+\mathcal V_*H_1(r)+z_*\bigl(H_2(r)+\mathcal V_*H_3(r)\bigr)$
(\lean{CompanionD.pr}), with $H_i$ of \eqref{eq:Hermite},
$z_*$ = \lean{SLean.Cert.zStar} and $\mathcal V_*=\mathcal V(z_*,\eta_*)$
(\lean{CompanionD.Vs}). Since $\mathcal V_*=i\mathsf v_*$
(\lean{CompanionD.Vs_eq}), this is the cell profile $p$ of
Section~\ref{sec:eigfun}. $\mathfrak g(y)=|p(\{y\})|^2/\int_0^1|p|^2$, with
$\{\cdot\}$ the fractional part, as in Section~\ref{sec:eigfun}
(\lean{CompanionD.gfrak}). $\mathfrak d(x)=2\sin^2\pi x$ and
$\mathfrak p(x)=x^2-x$, as in \eqref{eq:fh} (\lean{CompanionD.dfrak},
\lean{CompanionG.pfrak}); $\mathfrak b=\tfrac12\int_0^1\mathfrak g\,\mathsf m'$,
as in \eqref{eq:bdef} (\lean{CompanionG.bfrak});
$z_\infty=-\mathfrak d/(8\pi^2)+\mathfrak b\mathfrak p$, as in
\eqref{eq:zinf} (\lean{CompanionG.zinf}). $E(t)$, for $t\ge0$, is as in
\eqref{eq:Edef} (\lean{CompanionG.Efun}). For $f\in L^2$ and $t\ge0$,
$P_2(f)(t)=-\int_0^tT(t-s)\bigl[T(s)f\,\partial_xT(s)f\bigr]\,ds$
(\lean{CompanionB.P2}; the integrand is integrable by
\lean{CompanionB.Qc_tr}). Since $\partial_x(v^2)=2v\,\partial_xv$ for
$v\in H_0^1$, this is $P_2(f)(t)$ of Section~\ref{sec:quad}. Further
$\mathring P_{2,h}(f)(t)=\int_0^t\mathring T_h(t-s)\mathcal G_h\bigl[(\mathring T_h(s)\mathring\Pi_hf)^2\bigr]\,ds$
(\lean{CompanionF.P2h}); for $f\in X_h^{\C}$ this is
$\mathring P_{2,h}(f)(t)$ of Section~\ref{sec:quad}. Finally
$\mathring D_h(f)(t)=P_2(f)(t)-\mathring P_{2,h}(\mathring\Pi_hf)(t)$, as in
Section~\ref{sec:quad} (\lean{CompanionH.Dh}).

\paragraph{The tails.}
\begin{enumerate}
\item[(H0)] \lean{CompanionFinal.H0} is the conjunction of the Lean
statements of Section~\ref{sec:A3}, with the implicit $L$ written as
$\forall L>0$. The text is the same as in the SLean source and in
Section~\ref{sec:A3}. The translation of Section~\ref{sec:A3} applies.
This is Theorem~\ref{thm:top} and \eqref{eq:zgv}.

\item[(H1)] (\lean{CompanionA.H1}) This is Corollary~\ref{cor:failures}.
Let $L>0$ and $h=L/N$. Here
$L^2(0,L)$ is taken for Lebesgue measure on $(0,L]$, and $X_h^{\C}$,
$\mathring A_h$, $e^{t\mathring A_h}$ and $\chi'$ are defined as above on
$(0,L)$ (\lean{CompanionA.VhL}, \lean{CompanionA.ASL},
\lean{CompanionA.TSL}, \lean{CompanionA.dVL}). The norm of
$\mathring A_he^{t\mathring A_h}$ (\lean{CompanionA.CPop}) is the operator
norm on $X_h^{\C}$ with the norm of $L^2(0,L)$.
 \begin{enumerate}
 \item[(a)] For every $N\ge512$ and $0<t\le1$,
 $\norm{\mathring A_he^{t\mathring A_h}}\ge100h^{-3}e^{-1600t/L^3}$.
 There is no $C$ with
 $\norm{\mathring A_he^{t\mathring A_h}}\le Ct^{-3/2}$ for all $N\ge512$
 and all $0<t\le1$.
 \item[(b)] For every $T>0$ and every $N\ge512$ there is a real
 $\chi\in X_h^{\C}$, $\chi\ne0$, with
 \[
 \int_0^T\norm{(e^{t\mathring A_h}\chi)'}_{L^2(0,L)}^2\,dt
 \ge\frac{L^3}{256}\bigl(1-e^{-3200T/L^3}\bigr)h^{-2}\norm\chi_{L^2(0,L)}^2 .
 \]
 \end{enumerate}

\item[(H2)] (\lean{H2}) For every $t>0$ and every $c<e^{-\mathcal C_St}$,
$\norm{\mathring T_h(t)\mathring\Pi_h-T(t)}\ge c$ for all large $N=1/h$.
That is, $\liminf_{h\to0}\norm{\mathring T_h(t)\mathring\Pi_h-T(t)}\ge e^{-\mathcal C_St}$.

\item[(H3a)] (\lean{CompanionE.H3a}) For every $0<t_0<T_0$,
\[
 \sup_{t_0\le t\le T_0}\Bigl|\norm{\bigl(T(t)-\mathring T_h(t)\mathring\Pi_h\bigr)f_h}_2
 -e^{-\mathcal C_St}\Bigr|\to0\qquad\text{as }h=1/N\to0 .
\]
Since $f_h\in X_h$, $\mathring\Pi_hf_h=f_h$; so the projection can be
dropped here and in (H3b), as in Corollary~\ref{cor:lineargap} and
Theorem~\ref{thm:total}.

\item[(H3c)] (\lean{CompanionE.H3c})
 \begin{enumerate}
 \item[(i)] (\lean{CompanionE.SWell}) For every $N\ge1$, $h=1/N$, and
 every real $u_0\in L^2$, $\mathring u_h[u_0]$ solves the S-scheme with
 datum $u_0$, and every solution equals it on $[0,\infty)$.
 \item[(ii)] Let $\varepsilon_N>0$ for every $N\in\mathbb N$, with
 $N\varepsilon_N\to0$, and $h=1/N$. Then, for every $0<t_0<T_0$,
 \[
 \sup_{t_0\le t\le T_0}\Bigl|\varepsilon_N^{-1}
 \norm{\mathcal U(\varepsilon_Nf_h)(t)-\mathring u_h[\varepsilon_Nf_h](t)}_2
 -e^{-\mathcal C_St}\Bigr|\to0\qquad\text{as }h\to0 .
 \]
 \end{enumerate}

\item[(H4)] (\lean{CompanionE.H4}) Fix $t>0$.
 \begin{enumerate}
 \item[(a)] Both hold: $\liminf_{N\to\infty}\norm{\mathring L_h(t)}\ge e^{-\mathcal C_St}$,
 and $\mathring L_h(t)f\to0$ in $L^2$ for every $f\in L^2$ ($h=1/N$). In
 Lean the first part is written as in (H2).
 \item[(b)] For every sequence of integers $N_k\ge1$ with $N_k\to\infty$,
 $h_k=1/N_k$, and every sequence $r_k>0$ with $r_k\to0$, there is a real
 $f\in L^2$ such that the numbers
 $\norm{\mathring L_{h_k}(t)f}_2/r_k$ are unbounded.
 \end{enumerate}

\item[(H6i)] (\lean{H6i})
$\norm{A_h^{-1}\Pi_h-A^{-1}}_{L^2\to L^2}\le2h$ ($h=1/N$, $N\ge1$). For
every $f\in L^2(0,1)$, $\varrho=A^{-1}f$ is the unique $\varrho\in H^3(0,1)$
with $-\varrho'''=f$ and $\varrho(0)=\varrho(1)=\varrho'(1)=0$ (in Lean,
$\varrho\in C^2[0,1]$ with $\varrho''$ absolutely continuous and
$-\varrho'''=f$ a.e.).

\item[(H6ii)] (\lean{H6ii}) Fix $f\in L^2$ and $T_0\in\R$. For every
$\delta>0$ there is $N_0$ such that $\norm{T_h(t)\Pi_hf-T(t)f}_2<\delta$
for all $N\ge\max\{1,N_0\}$ ($h=1/N$) and all $t\in[0,T_0]$. Moreover
$T(0)=I$; $T(t+s)=T(t)T(s)$ for $s,t\ge0$; $t\mapsto T(t)f$ is continuous
on $[0,\infty)$; and $A^{-1}\bigl(T(t)f-f\bigr)=\int_0^tT(s)f\,ds$ for
$t\ge0$.

\item[(H6iii)] (\lean{CompanionB.H6iiiV}) For every $0<\nu<3\pi^2/4$
there are $r_0,R_0$ such that (U)$_{\nu,r_0,R_0}$ holds. Here
(U)$_{\nu,r_0,R_0}$ (\lean{CompanionB.H6iiiAtV}) states: $r_0>0$,
$R_0>0$, and for every real $f\in L^2$ with $\norm f_H<r_0$,
 \begin{enumerate}
 \item[(a)] for every $N\ge1$, $h=1/N$, $u_h=u_h[f]$ solves the W-scheme
 with datum $f$; every solution equals $u_h$ on $[0,\infty)$;
 $u_h\in\mathcal X_0(\nu)$; and $\mathcal E_\nu(u_h)\le R_0$;
 \item[(b)] $\mathcal U(f)\in\mathcal X_0(\nu)$;
 $\mathcal E_\nu(u_h-\mathcal U(f))\to0$ and
 $\sup_{t\ge0}\norm{u_h(t)-\mathcal U(f)(t)}_2\to0$ as $h\to0$;
 \item[(c)] $w=\mathcal U(f)$ satisfies: $w(t)$ is real for $t\ge0$;
 $\mathcal E_\nu(w)\le R_0$; for every $t\ge0$ the function
 $s\mapsto T(t-s)[w(s)\,\partial_xw(s)]$ is Bochner integrable on
 $(0,t]$ and \eqref{eq:mild} holds; if $\tilde w\in\mathcal X_0(\nu)$,
 $\mathcal E_\nu(\tilde w)\le R_0$ and $\tilde w$ satisfies
 \eqref{eq:mild} for every $t\ge0$, then $\tilde w(t)=w(t)$ for every
 $t\ge0$.
 \end{enumerate}
By (a), $u_h[f]$ is not the default. By (b), $\mathcal U(f)(t)$ is the
limit for every $t\ge0$. In the uniqueness clause, the integral in
\eqref{eq:mild} for $\tilde w$ is a genuine Bochner integral by the
internal theorem \lean{CompanionB.probe_mild_integrableV}. It holds for
$0<\nu<3\pi^2/4$ and needs only $\tilde w\in\mathcal X_0(\nu)$ and $t\ge0$.

\item[(B)] (\lean{CompanionFinal.H3bH5V}) For every $0<\nu<3\pi^2/4$
there are $r_0,R_0,\varepsilon_0$ such that (U)$_{\nu,r_0,R_0}$,
(H3b)$_{r_0,\varepsilon_0}$ and (H5)$_{\varepsilon_0}$ hold, with the same
numbers.
 \begin{enumerate}
 \item[(H3b)] (\lean{CompanionH.H3bAtE}) $\varepsilon_0>0$, with $r_0$ as
 in (U)$_{\nu,r_0,R_0}$, and there is $C>0$ such that
 $\sqrt2\varepsilon_0<r_0$ and, for every $0<\varepsilon\le\varepsilon_0$
 and every $t>0$,
 \[
 \limsup_{h\to0}\Bigl|\norm{\mathcal U(\varepsilon f_h)(t)
 -\varepsilon\mathring T_h(t)\mathring\Pi_hf_h}_2-\varepsilon e^{-\mathcal C_St}\Bigr|
 \le C\varepsilon^3\qquad(h=1/N).
 \]
 Lean writes this as: for every $c>C\varepsilon^3$ the quantity is at
 most $c$ for all large $N=1/h$. Here
 $\varepsilon\mathring T_h(t)\mathring\Pi_hf_h=\mathring v_h(t)$ for
 $u_0=\varepsilon f_h$; so (H3b) is Theorem~\ref{thm:total}(b) with
 $\mathcal U(\varepsilon f_h)$ in place of $u$.
 \item[(H5)] (\lean{CompanionH.H5At}) For every $t>0$, both hold:
  \begin{enumerate}
  \item[(a)] for every $N\ge512$, $h=1/N$,
  $\varepsilon^{-2}\mathring\Cd_h(\varepsilon,t)\to\mathring D_h(f_h)(t)$ in
  $L^2$ as $\varepsilon\downarrow0$; moreover
  $\mathring D_h(f_h)(t)\to-E(t)$ in $L^2$ as $h\to0$, and $E(t)\ne0$. In
  Lean the limit is written $P_2(f_h)(t)-\mathring P_{2,h}(f_h)(t)$; since
  $\mathring\Pi_hf_h=f_h$, it is $\mathring D_h(f_h)(t)$. Lean writes the
  limit as $h\to0$ with thresholds: for every $\delta>0$,
  $\norm{\mathring D_h(f_h)(t)+E(t)}_2<\delta$ for all large $N\ge512$;
  \item[(b)] for every $p>0$ and every $0<\varepsilon_1\le\varepsilon_0$,
  with $\varepsilon_0$ as in (H3b), there is no $C$ with
  $\norm{\mathring\Cd_h(\varepsilon,t)}_2\le Ch^p\varepsilon^2$ for all
  $N\ge512$ ($h=1/N$) and all $0<\varepsilon\le\varepsilon_1$.
  \end{enumerate}
 \end{enumerate}
Since $\norm{f_h}_2\le1$, we have $\norm{\varepsilon f_h}_H\le\sqrt2\varepsilon<r_0$
for $0<\varepsilon\le\varepsilon_0$. So, by (U)$_{\nu,r_0,R_0}$,
$\mathcal U(\varepsilon f_h)$ is a limit in (H3b) and (H5b). In (H5a) only
small $\varepsilon$ matter, and in (H3c) only large $N$. There
$\sqrt2\varepsilon$, resp.\ $\sqrt2\varepsilon_N$, is less than the radius
$r_0$ of (H6iii) for any fixed $\nu$, so the same holds. Since
$\varepsilon f_h$ is real, $\mathring u_h[\varepsilon f_h]$ is not the
default, by (H3c)(i).

\item[(F)] (\lean{CompanionFinal.FamilyProps}) For every $N\ge512$,
$h=1/N$: $\varphi_h\in X_h^{\C}$, $\mathring A_h\varphi_h=\zeta_h\varphi_h$,
$\norm{\varphi_h}_2=1$ and
$\bigl|\int_0^1\varphi_h^2\,dx\bigr|\le16h^3$; $\tilde f_h\ne0$ and
$\norm{\tilde f_h}_2^2=1+\Rea\int_0^1\varphi_h^2\,dx$; $f_h\in X_h^{\C}$,
$f_h$ is real and $\norm{f_h}_2=1$. Moreover $\varphi_h\rightharpoonup0$
and $f_h\rightharpoonup0$ as $h\to0$. There is $C>0$ such that
$|\varphi_h(x)|\le C$ and $|\varphi_h'(x)|\le Ch^{-1}$ for all $N\ge512$
and $x\in[0,1]$, for the $C^1$ representative of $\varphi_h$.
\end{enumerate}

\paragraph{The parallel theorem.}
\lean{CompanionH.Abstract.Prop41} reads as follows. Let $\mathcal X$ be a
complex Hilbert space. Let $X_n\subset\mathcal X$ ($n\in\mathbb N$) be
closed subspaces, with orthogonal projections $\Pi_n$, and let $A_n$ be
bounded operators on $X_n$. Let $\zeta_n\in\C$ and $x_n\in X_n$ satisfy
$A_nx_n=\zeta_nx_n$ and $\norm{x_n}=1$. Then:
\begin{enumerate}
\item[(a)] if every $A_n$ is invertible, $K$ is a bounded injective
operator on $\mathcal X$, $\norm{A_n^{-1}\Pi_n-K}\to0$ and
$|\zeta_n|\to\infty$, then $x_n\rightharpoonup0$;
\item[(b)] if $S$ is a $C_0$-semigroup on $\mathcal X$, $S(t_1)$ is
compact for some $t_1>0$, and $x_n\rightharpoonup0$, then both of the
following hold:
 \begin{enumerate}
 \item[(i)] for every $t\ge t_1$ and every $c$ with $e^{t\Rea\zeta_n}\ge c$
 for all large $n$, and every $c'<c$,
 $\norm{e^{tA_n}\Pi_n-S(t)}\ge c'$ for all large $n$; that is,
 $\liminf_n\norm{e^{tA_n}\Pi_n-S(t)}\ge\liminf_ne^{t\Rea\zeta_n}$;
 \item[(ii)] for every $t_0>0$, if $\norm{e^{t_0A_n}\Pi_n-S(t_0)}\to0$,
 then $\Rea\zeta_n\to-\infty$.
 \end{enumerate}
\end{enumerate}
Here $e^{tA_n}\Pi_n$ (\lean{CompanionH.Abstract.expPi}) and
$A_n^{-1}\Pi_n$ (\lean{CompanionH.Abstract.invPi}) are operators on
$\mathcal X$. \lean{CompanionH.Abstract.invPi} uses Mathlib's
\lean{Ring.inverse}, which is $A_n^{-1}$ when $A_n$ is invertible, as
assumed in (a). $S:\R\to\mathcal B(\mathcal X)$ is a $C_0$-semigroup
(\lean{CompanionH.Abstract.IsC0Semigroup}) if $S(0)=I$,
$S(s+t)=S(s)S(t)$ for $s,t\ge0$, and $t\mapsto S(t)x$ is continuous on
$[0,\infty)$ for every $x$. Weak convergence is $\langle y,x_n\rangle\to0$
for every $y$ (\lean{CompanionH.Abstract.WeakNullSeq}). The theorem is
proved for every universe level $u$ of $\mathcal X$.

\paragraph{Files and axioms.}
The Lean development of this subsection has $120$ modules, in the
namespaces \lean{CompanionPilot}, \lean{CompanionA}--\lean{CompanionH} and
\lean{CompanionFinal}. It depends on SLean (Section~\ref{sec:A4}), on the
Lean package \lean{OrderingError} and on the computation package
\lean{Companion}. The command \texttt{\#print axioms}, applied to
\lean{CompanionFinal.main_theorem}, \lean{CompanionFinal.prop318} and the
proofs of the tails, lists only \texttt{propext}, \texttt{Classical.choice}
and \texttt{Quot.sound} (\texttt{CompanionFinal/Audit.lean}). The sources
contain no \texttt{sorry} and no new \texttt{axiom}. All sources, including
the modules of SLean, \lean{OrderingError} and \lean{Companion} that are
used, are in the folder \texttt{lean/} of the ancillary files of
arXiv:2611.12585; \texttt{lean/README.md} gives the build instructions.

\paragraph{Objects used only in Table~\ref{A5:corr}.}
The following definitions serve the internal statements in the table. The
tails do not use them. Let $N\ge1$ and $h=1/N$.
\begin{itemize}
\item With $\{\cdot\}$ the fractional part,
$\mathsf m_h(x)=h\,\mathsf m(\{x/h\})$ and
$\mathsf m_h'(x)=\mathsf m'(\{x/h\})$ (\lean{CompanionC.mfun},
\lean{CompanionC.mder}, \lean{CompanionD.mhD}).
\item $\mathfrak b$ equals the constant \lean{Companion.cellConstant} of
the computation package (\lean{CompanionG.bfrak_eq}).
\item For $q\in L^2$, $z_h[q]=\mathring A_h^{-1}\mathcal G_h[q]$
(\lean{CompanionC.zS}) and,
for $N\ge512$, $z_h=z_h[|\varphi_h|^2]$ (\lean{CompanionG.zh}) and
$\mathcal O_h(t)=\int_0^te^{2\zeta_hs}\mathring T_h(t-s)\mathcal G_h[\varphi_h^2]\,ds$
(\lean{CompanionD.Oh}).
\item For $\vartheta\in L^2$ let $F_1(x)=\int_0^x\vartheta\,d\mu$,
$F_2(x)=\int_0^xF_1$ and $F_3(x)=\int_0^xF_2$, for all real $x$
(\lean{prim1}, \lean{prim2}, \lean{prim3}). The Green operator $G$
(\lean{G}) maps $\vartheta$ to the class of
$-F_3(x)+(2F_3(1)-F_2(1))x+(F_2(1)-F_3(1))x^2$. The dual data are
$\psi_0=F_3-F_3(1)x^2$, $g=F_2-2F_3(1)x$ and $g_1=F_1-2F_3(1)$
(\lean{CompanionC.psi0}, \lean{CompanionC.gD}, \lean{CompanionC.gD1}).
$\mathcal S_h\vartheta\in X_h^{\C}$ (\lean{CompanionC.Sc}) has the nodal
values $\psi_0(x_j)$ and the scaled slopes $h\,p_j$, with
$p_j=g(x_j)-g(1)-\int_{x_j}^1\mathsf m_h\vartheta\,d\mu$
(\lean{CompanionC.nodeU_Sc}, \lean{CompanionC.nodeW_Sc}).
\item $E$ equals \lean{CompanionF.Efun} (\lean{CompanionH.Efun_eq}).
\item Objects defined twice in Lean agree by
\lean{CompanionFinal.unified_objects}: the two definitions of
$\mathring P_{2,h}$ on $X_h^{\C}$ only, and the two definitions of
$\mathcal G_h$ for $q\in L^2$.
\end{itemize}

\paragraph{Correspondence.}
Table~\ref{A5:corr} lists every numbered result and remark of
Sections~\ref{sec:eigfun}--\ref{sec:lineargap} and Section~\ref{sec:gap},
and the displayed facts \eqref{eq:bclosed}, \eqref{eq:bvalue},
\eqref{eq:average}, \eqref{eq:taylor} and the compactness of
Section~\ref{sec:galerkin}. ``Tail'' means a conjunct of
\lean{CompanionFinal.Tails}. ``Internal'' means a theorem with a Lean~4
proof that is not part of the main theorem; its statement is given, with
the objects above and the range of $N$ of the Lean statement. ``P1''
marks internal theorems that assume (P1). ``Remark'' means commentary
without a separate Lean declaration. ``Not in Lean'' means that there is
no Lean statement.

{\small
\begin{longtable}{@{}>{\raggedright\arraybackslash}p{.14\textwidth}>{\raggedright\arraybackslash}p{.12\textwidth}>{\raggedright\arraybackslash}p{.68\textwidth}@{}}
\caption{Results of Sections~\ref{sec:eigfun}--\ref{sec:lineargap} and
Section~\ref{sec:gap} and their Lean declarations.}\label{A5:corr}\\
\toprule
Result & Status & Lean declaration and statement\\\midrule
\endfirsthead
\toprule
Result & Status & Lean declaration and statement\\\midrule
\endhead
\bottomrule
\endfoot
Cor.~\ref{lem:density}(a) & internal & \lean{CompanionD.S23_density}: there is $C$ with $\bigl||\varphi_h(x)|^2-\mathfrak d(x)\mathfrak g(x/h)\bigr|\le Ch$ for $N\ge512$, $x\in[0,1]$ ($C^1$ representative).\\
 & tail (F) & $|\varphi_h|\le C$, $|\varphi_h'|\le Ch^{-1}$ on $[0,1]$ (\lean{CompanionA.S11_sup}).\\
Cor.~\ref{lem:density}(b) & internal & \lean{CompanionD.S23_b}: there is $C$ such that for $N\ge512$, $B\in\R$ and $G,g:\R\to\C$ with $g$ integrable on $[0,1]$, $G(x)=G(0)+\int_0^xg$ and $|G(x)|\le B$ for $x\in[0,1]$, and $|g|\le B$ a.e.\ on $[0,1]$: $\bigl|\int_0^1(|\varphi_h|^2-\mathfrak d)G\,dx\bigr|\le ChB$.\\
Cor.~\ref{lem:density}(c) & internal & \lean{CompanionD.S23_c} with \lean{CompanionG.bfrak_eq}: there is $C$ with $\bigl|\tfrac12\int_0^1|\varphi_h|^2\mathsf m_h'\,dx-\mathfrak b\bigr|\le Ch$ for $N\ge512$.\\
Cor.~\ref{lem:density}(d) & tail (F) & $|\int_0^1\varphi_h^2|\le16h^3$; $\norm{\tilde f_h}_2^2=1+\Rea\int_0^1\varphi_h^2$; $\varphi_h,f_h\rightharpoonup0$.\\
 & internal & \lean{CompanionH.ftilA_norm_sq_sub}: $|\norm{\tilde f_h}_2^2-1|\le16h^3$ for $N\ge512$.\\
\eqref{eq:bclosed} & internal & \lean{Companion.cell_constant_closed} with \lean{CompanionG.bfrak_eq}: for real $c,s,\mathsf v$ with $c^2+s^2=1$, $\tfrac12\int_0^1|p|^2\mathsf m'/\int_0^1|p|^2$ for the profile $p$ with $z_*=c+is$, $\mathcal V_*=i\mathsf v$ equals the right side of \eqref{eq:bclosed}; at $(c_*,s_*,\mathsf v_*)$ this is $\mathfrak b$. \lean{Companion.cellConstant_decimal}: $0.00942<\mathfrak b<0.00943$.\\
\eqref{eq:bvalue} & internal & \lean{CompanionG.probe_bfrak_pos}: $\mathfrak b>0$.\\
\eqref{eq:average} & internal & \lean{CompanionD.S22_average}: for $N\ge1$; $\mathsf c:\R\to\C$ one-periodic, integrable on $[0,1]$, $\bar{\mathsf c}=\int_0^1\mathsf c$, $|\mathsf c-\bar{\mathsf c}|\le M$ a.e.\ on $\R$; $g$ integrable on $[0,1]$ and $G(x)=G(0)+\int_0^xg$ on $[0,1]$: $\bigl|\int_0^1G(x)\mathsf c(x/h)\,dx-\bar{\mathsf c}\int_0^1G\bigr|\le hM\int_0^1|g|$.\\
Lem.~\ref{lem:inverse} & internal & \lean{CompanionC.psi0_solves}: for $\vartheta\in L^2$, $\psi_0(0)=g(0)=\psi_0(1)=0$, $\psi_0'=g$, $g'=g_1$ and $g_1(x)=g_1(0)+\int_0^x\vartheta\,d\mu$ for all real $x$. \lean{CompanionC.dual_identity}: for $N\ge1$, $\vartheta\in L^2$ and $v\in X_h^{\C}$, $\int_0^1v''(\mathcal S_h\vartheta)'\,dx=\int_0^1v\vartheta\,dx$. \lean{CompanionC.AS_injective}: $\mathring A_h$ is injective for $N\ge1$.\\
Lem.~\ref{lem:twoscale} & internal & For $N\ge1$ and $\vartheta,f\in L^2$: \lean{CompanionC.gD_bound}: $|g|,|g'|\le3\norm\vartheta_2$ on $[0,1]$; \lean{CompanionC.two_scale}: $|(\mathcal S_h\vartheta)'-g+g(1)\mathsf m_h'|\le13h\norm\vartheta_2$ on $[0,1]$; \lean{CompanionC.inverse_load}: $\bigl|\int_0^1z_h[f]\vartheta-\tfrac12\int_0^1fg+g(1)\tfrac12\int_0^1f\mathsf m_h'\bigr|\le7h\norm f_1\norm\vartheta_2$; \lean{CompanionC.zS_norm_le}: $\norm{z_h[f]}_2\le10\norm f_1$.\\
Lem.~\ref{lem:slow} & internal & \lean{CompanionG.S24}: there is $C$ with $\norm{z_h-z_\infty}_2\le Ch$ for $N\ge512$.\\
Lem.~\ref{lem:fast} & internal & \lean{CompanionD.S25}, \lean{CompanionD.S25_resolvent}: there is $C$ such that for $N\ge512$, $2\zeta_h-\mathring A_h$ has a two-sided inverse on $X_h^{\C}$ of norm $\le h^3/15$, and $\norm{\mathcal O_h(t)}_2\le Ch^2$ for $t\ge0$.\\
Prop.~\ref{prop:galerkin}(a) & tail (H6i) & $\norm{A_h^{-1}\Pi_h-A^{-1}}\le2h$ for $N\ge1$, $h=1/N$. The proof uses \lean{JW_sub_G} ($\norm{A_h^{-1}\Pi_h-G}\le2/N$) and \lean{Ainv_eq} ($A^{-1}=G$).\\
Prop.~\ref{prop:galerkin}(b) & tail (H6ii) & First clause.\\
Prop.~\ref{prop:galerkin}(c) & tails (H6i), (H6ii) & The boundary value problem; $T(0)=I$, the semigroup law, continuity, $A^{-1}(T(t)f-f)=\int_0^tT(s)f\,ds$. Also \lean{CompanionH.Tsg_C0}.\\
Compactness (Sec.~\ref{sec:galerkin}) & internal & \lean{CompanionH.G_compact}: $A^{-1}=G$ is compact. \lean{CompanionH.Tsg_compact} (P1): $T(t)$ is compact for $t>0$.\\
Lem.~\ref{lem:TK} & internal & \lean{CompanionC.S21_strong}: for $f\in L^2$, $T_0\in\R$ and $\delta>0$ there is $N_0$ with $\norm{\mathring T_h(t)\mathring\Pi_hf-T(t)f}_2<\delta$ for $N\ge\max\{1,N_0\}$, $t\in[0,T_0]$. \lean{CompanionC.S21_rate}: $\norm{(\mathring T_h(t)\mathring\Pi_h-T(t))G^2x}_2\le240(1+t)\norm x_2/N$ for $N\ge1$, $t\ge0$. The fixed-$t$ part is in (H4a).\\
Prop.~\ref{prop:necessary} & parallel & \lean{CompanionFinal.prop318}.\\
Cor.~\ref{cor:principleS}(a) & internal & \lean{CompanionC.S20_rate}: $\norm{\mathring A_h^{-1}\mathring\Pi_h-A^{-1}}\le16/N$ for $N\ge1$.\\
Cor.~\ref{cor:principleS}(b) & tail (F) & $\varphi_h\rightharpoonup0$.\\
 & internal & Via Proposition~\ref{prop:necessary}(a): \lean{CompanionH.phiN_weakNull_of_prop41}, for $x_n=$ \lean{phiN}$(n+512)$ (\lean{CompanionH.xS}), another normalized eigenfunction family of $\mathring A_h$ for $\zeta_h$.\\
Cor.~\ref{cor:principleS}(c) & tail (H2) & \\
Cor.~\ref{cor:lineargap} & tail (H3a) & \\
Cor.~\ref{cor:fixed} & tail (H4) & \\
Rem.~\ref{rem:leanscope} & remark & No separate declaration. Its scope statements are those of this table.\\
Prop.~\ref{prop:Ulimit}(a)--(c) & tails (H6iii), (B) & (U)$_{\nu,r_0,R_0}$ (\lean{CompanionB.H6iiiAtV}).\\
Prop.~\ref{prop:Ulimit}, last sentence & not in Lean & \\
Prop.~\ref{thm:quad}(a) & internal & \lean{CompanionG.S27}: for $T_0>0$ and $\delta>0$ there is $N_0$ with $\norm{\mathring P_{2,h}(\tilde f_h)(t)-E(t)}_2<\delta$ for $N\ge\max\{512,N_0\}$, $t\in[0,T_0]$.\\
Prop.~\ref{thm:quad}(b) & not in Lean & Except $E(0)=0$, which is internal (\lean{CompanionF.probe_Efun_zero}).\\
Prop.~\ref{thm:quad}(c) & tail (B) & (H5a): $E(t)\ne0$ for $t>0$. The proof uses \lean{CompanionF.S28} with \lean{CompanionH.Efun_eq}.\\
Cor.~\ref{cor:Dh} & tail (F) & $\norm{f_h}_2=1$, $f_h\rightharpoonup0$.\\
 & internal & \lean{CompanionE.S16} (P1): for $T_0\ge0$ and $\delta>0$, $\norm{P_2(f_h)(t)}_2<\delta$ for all large $N$ and $t\in[0,T_0]$. \lean{CompanionH.S29} (P1): for $T_0>0$ and $\delta>0$, $\norm{\mathring D_h(f_h)(t)+E(t)}_2<\delta$ for all large $N$ and $t\in[0,T_0]$. \lean{CompanionH.S29_norm} (P1): $\norm{\mathring D_h(f_h)(t)}_2\to\norm{E(t)}_2>0$ for $t>0$. At each fixed $t>0$, the limit $\mathring D_h(f_h)(t)\to-E(t)$ and $E(t)\ne0$ are also in (H5a).\\
Rem.~\ref{rem:cell} & remark & No separate declaration. The load identity is \lean{CompanionC.inverse_load}.\\
Thm.~\ref{thm:total}(a) & tail (H3a) & \\
Thm.~\ref{thm:total}(b) & tail (B) & (H3b)$_{r_0,\varepsilon_0}$, with $\mathcal U(\varepsilon f_h)$ in place of $u$.\\
Thm.~\ref{thm:total}(c) & tail (H3c) & With $\mathcal U(\varepsilon_Nf_h)$ in place of $u$.\\
\eqref{eq:taylor} & internal & \lean{CompanionB.taylor12_pointwise} with \lean{CompanionE.Ucal_eq}: for $\norm f_2<1/72$ and $t\ge0$, $\norm{\mathcal U(f)(t)-T(t)f}_2\le216e^{-t}\norm f_2^2$ and $\norm{\mathcal U(f)(t)-T(t)f-P_2(f)(t)}_2\le7776e^{-t}\norm f_2^3$; with $\mathcal U$ in place of $u$.\\
Cor.~\ref{cor:noproduct}(a) & tail (B) & (H5a). The proof uses \lean{CompanionH.Lambda_N_eq} and \lean{CompanionH.Lam_tendsto}.\\
Cor.~\ref{cor:noproduct}(b) & tail (B) & (H5b).\\
\end{longtable}}

%% file: appB_numerics.tex
\section{Numerical experiments}\label{app:num}

The theorems carry all conclusions. The computations prove nothing. On
the computed data and meshes, they give the actual sizes of the errors and
constants of the two schemes. They also show how these sizes change under
refinement. Unless stated otherwise $L=1$, $h=1/N$, the meshes are
uniform, and all computations are semidiscrete in space. Time integration
is done to high accuracy and checked. Theorem numbers of~\cite{Main}
refer to its version~2.

\begin{remark}[how the computations were made]\label{rem:numerics}
The author set the purpose of each computation.
Section~\ref{app:S2} tests, for the W-scheme, the form of the estimates
of~\cite[Theorems~4.2 and~4.3]{Main}. Their constants are fixed on coarse
meshes and then checked on finer ones. Section~\ref{app:fh} quantifies
the gap of Section~\ref{sec:gap} on the data $f_h$. Its meshes
$N=512,1024,2048$ lie in the proved range $N\ge512$.
Sections~\ref{app:eig}, \ref{app:constants} and~\ref{app:further} answer
narrower questions: the constants of Theorem~\ref{thm:top}, the semigroup
constants of the two schemes, fixed data and running times.

The programs were written and run by GPT-6~Astra and Claude Opus~5.5
under the author's direction. GPT-6~Astra did the computations of
Sections~\ref{app:eig} and~\ref{app:S2}. It also computed the values of the
quadratic coefficient in Section~\ref{app:further}, in an independent check.
Claude Opus~5.5 did the other computations of
Sections~\ref{app:constants}, \ref{app:fh} and~\ref{app:further}. It checked
the values at $N=512$ and $1024$ with a separate program. Before the
computations of Sections~\ref{app:S2} and~\ref{app:fh}, the predicted
values and the pass/fail thresholds were written to files, and their
hashes were recorded. Failed predictions are reported. The author did not
check the programs line by line. The programs, their outputs and logs are
in the ancillary files of arXiv:2611.12585, in one folder for each of
Sections~\ref{app:eig}--\ref{app:further}; the file \texttt{README.md}
there describes them.

After the computations, the design and the matrix computations of
Section~\ref{app:fh} were assessed independently. GPT-6~Astra (in Codex)
and Claude Fable~5.1 each worked in a new session, without access to the
programs. Each received the two papers, a two-page note of questions in
the language of matrices, and the exported matrices and eigenvectors for
$N=512,1024,2048$. Each was asked whether the design is reasonable and
whether the matrix computations are reliable, and to answer without large
computations. Each made only small checks on the exported matrices. The
two assessments agree on the following points.
\begin{enumerate}
\item[(i)] The design suits a test of bounds that are uniform on the unit
ball, on data that depend on the mesh. It gives no ranking of the schemes
for fixed data or for computational cost.
\item[(ii)] The eigenvalue $\zeta_h$ is well conditioned. Both computed
its condition number in the mass norm as $1.00003$, $1.00001$,
$1.000002$ for $N=512,1024,2048$.
\item[(iii)] The S-scheme should not be integrated by Radau~IIA with steps that do
not resolve the spurious mode.
\item[(iv)] The agreement of two eigendecompositions is a cross-check,
not an error bound. P10 remains a failure.
\item[(v)] The W-scheme values carry the uncertainty of the reference check.
\end{enumerate}
The two assessments differ in how far they trust the full
eigendecomposition for the higher-order terms of the S-scheme. Both proposed the
same further check, an evaluation that avoids the eigenvector matrix. It
was not carried out. After these assessments, the statements on precision
in Section~\ref{app:fh} were revised. Neither assessment repeated the
trajectory computations.
\end{remark}

\subsection{The high-frequency spurious eigenpairs}\label{app:eig}
The program \texttt{F\_certify.py}, which uses only the
Python standard library and rational interval arithmetic, gives the
enclosures in Table~\ref{tab:constants}. They are independent numerical
checks. The Lean proof of Theorem~\ref{thm:top} uses its own rational
certificates, which Lean checks.

\begin{table}[H]
\centering\small
\caption{Interval enclosures of the constants of Theorem~\ref{thm:top}.}
\label{tab:constants}
\begin{tabular}{ll}\toprule
Quantity & Enclosure\\\midrule
$\xi_*/\pi$ & $[0.32272944586606,0.32272944586607]$\\
$\varpi_*$ & $[-0.27655651304531,-0.27655651304530]$\\
$\varpi''$ & $[0.372284741795,0.372284741796]$\\
$\Omega_*$ & $[116.1537354790,116.1537354791]$\\
$210\pi^2\varpi''$ & $[771.60365647,771.60365648]$\\
$\Rea\beta$ & $[2.041699771,2.041699772]$\\
$\Ima\beta$ & $[-0.319050216,-0.319050215]$\\
$\mathcal C_S$ & $[1575.383008,1575.383009]$\\
\bottomrule\end{tabular}
\end{table}

We assembled the real generalized pencil $Kv=\eta Mv$ and computed its
spectrum in double precision. Independently, we solved the boundary
determinant at $75$ decimal digits, and at $N=512$ also at $105$ digits.
The two precisions agree to better than $4\cdot10^{-76}$ in relative
eigenvalue error. The normalized determinant residual is below
$10^{-50}$. Table~\ref{tab:pairs} compares the first-mode eigenvalue with
the expansions \eqref{eq:im}--\eqref{eq:re}. Over the five computed
optical modes, we compared the double-precision matrix eigenvalue with the
high-precision boundary root. The largest relative difference is
below $9\cdot10^{-15}$ for the full eigenvalue and below $1.6\cdot10^{-8}$
for its real part.

\begin{table}[H]
\centering\small
\caption{The first spectral-top eigenvalue and the expansions
\eqref{eq:im}--\eqref{eq:re} (with the $O(N^{-4})$ and $O(N^{-1})$ terms
omitted).}\label{tab:pairs}
\begin{tabular}{rcccc}
\toprule
& \multicolumn{2}{c}{$h^3|\Ima\zeta_h|$} & \multicolumn{2}{c}{$L^3\Rea\zeta_h$}\\
\cmidrule(lr){2-3}\cmidrule(lr){4-5}
$N$ & computed & expansion & computed & expansion\\
\midrule
32 & 115.3984463 & 115.3927035 & $-1581.319$ & $-1575.383$ \\
64 & 115.9647737 & 115.9644166 & $-1582.795$ & $-1575.383$ \\
128 & 116.1065454 & 116.1065231 & $-1580.188$ & $-1575.383$ \\
256 & 116.1419485 & 116.1419471 & $-1578.058$ & $-1575.383$ \\
512 & 116.1507903 & 116.1507902 & $-1576.788$ & $-1575.383$ \\
\bottomrule
\end{tabular}
\end{table}

\paragraph{The threshold.}
The script \texttt{S1\_branch\_scan.py} follows the first optical branch
on every integer mesh $3\le N\le128$ and on $N=256,512$, forward from
$N=3$ and backward from the asymptote at $N=512$. Both directions select
the same eigenpair. The bounds \eqref{eq:bounds} fail at $N=3,4,5,6$ and
hold at every tested mesh from $N=7$ on. The largest scaled damping is
$1583.30$, at $N=48$. This is an a posteriori observation on the tested
set, not a proof that $N\ge7$ suffices. The scan does not cover
$129\le N\le255$, $257\le N\le511$ or $N>512$.

\subsection{Semigroup constants}\label{app:constants}
Table~\ref{tab:semigroup} gives two constants of mesh-uniform
regularity. For the S-scheme they are evaluated on the spurious eigenvector, in
$L^2$. The Crandall--Pazy witness is $|\zeta_h|\max_{t>0}t^{3/2}e^{t\Rea\zeta_h}
=|\zeta_h|(3/(2|\Rea\zeta_h|))^{3/2}e^{-3/2}$. The Kato witness is
$\int_0^1\norm{\partial_x\mathring T_h(t)v}_2^2dt/\norm v_2^2$ with
$v=\sqrt2\Rea\varphi_h$. For the W-scheme the table gives the worst cases over all
data. We compute the operator norm at each sampled time. The reported
Crandall--Pazy value is the maximum over these times. It is a numerical
estimate of $\sup_{0<t\le1}t^{3/2}\norm{A_hT_h(t)}_{H\to H}$. The Kato value is the largest eigenvalue of the Kato Gramian in
$L^2$, computed with the integral over $(0,\infty)$. For the S-scheme the worst Kato value agrees with
the witness to three digits at $N=128,256,512$ ($160$, $641$, $2567$).

\begin{table}[H]
\centering\small
\caption{Constants of mesh-uniform regularity. All rows in $L^2$, except the
Crandall--Pazy constant of the W-scheme, which is computed in the norm of $H$.}\label{tab:semigroup}
\begin{tabular}{@{}lcccccc@{}}
\toprule
$N$ & 64 & 128 & 256 & 512 & 1024 & 2048\\
\midrule
Crandall--Pazy, W-scheme (worst case) & $0.33$ & $0.36$ & $0.38$ & $0.39$ & -- & --\\
Crandall--Pazy, S-scheme (witness) & -- & -- & -- & $1.02\times10^5$ & $8.17\times10^5$ & $6.54\times10^6$\\
Kato, W-scheme (worst case) & -- & $0.63$ & $0.64$ & $0.65$ & -- & --\\
Kato, S-scheme (witness) & -- & -- & -- & $2.57\times10^3$ & $1.03\times10^4$ & $4.11\times10^4$\\
\bottomrule
\end{tabular}
\end{table}

The numerical results show that the Crandall--Pazy constant of the W-scheme
increases by $8.2\%$, $4.9\%$ and $4.1\%$ from $N=64$ to $128$, $256$ and
$512$, consistent with the proved mesh-uniform bound. In contrast, the
constants of the S-scheme grow by factors $8$ and $4$ per halving of $h$, as the
witnesses predict. The decay of the two schemes does not differ in this
way. Both have the spectral bound $s\approx-75.86$ on $N=32$--$256$. With
$\nu=0.99|s|$, the numerical estimates of
$\sup_{t\ge0}e^{\nu t}\norm{e^{tA_h}}_{H\to H}$ and
$\sup_{t\ge0}e^{\nu t}\norm{e^{t\mathring A_h}}_{L^2\to L^2}$ are $1.51$
and $1.68$ on all these meshes.

\subsection{The W-scheme on seven data}\label{app:S2}

\paragraph{Quantities.}
For the datum $\varepsilon f$ let $u$, $u_h$ be the nonlinear solutions,
and $v$, $v_h$ the linearized ones for the unit datum $f$. The
corresponding linear solutions for $\varepsilon f$ are $\varepsilon v$ and
$\varepsilon v_h$. This unit-datum convention is used in
Sections~\ref{app:S2}--\ref{app:further}, also for the S-scheme. In this
section the datum is normalized by $\norm f_H=1$ in the continuous norm
of $H$, before projection. In Section~\ref{app:fh} the datum is $f_h$, with
$\norm{f_h}_2=1$. Expand
$u=\varepsilon v+\varepsilon^2q+\varepsilon^3p_3+\cdots$, and likewise
$u_h$, with $q=P_2(f)$ and $q_h=P_{2,h}(\Pi_hf)$. For $\varepsilon>0$ put
\[
 r_\varepsilon=\frac{u-\varepsilon v-\varepsilon^2q}{\varepsilon^3},\qquad
 r_{h,\varepsilon}=\frac{u_h-\varepsilon v_h-\varepsilon^2q_h}{\varepsilon^3}.
\]
Thus the cubic remainder in Table~\ref{tab:B-quantities} contains all
terms of degree at least three. Table~\ref{tab:B-quantities} lists the
error quantities. The full
error is $u-u_h$. The ordering error is $\Cd_h=(u-u_h)-\varepsilon(v-v_h)$.
Its leading term is $\varepsilon^2(q-q_h)$.

\begin{table}[H]
\centering\small
\caption{Error quantities and the normalizing denominators of the
empirical constants; $a=\min\{1,ht_\wedge^{-3/2}\}$. The forms are those
of~\cite[Theorems~4.2 and~4.3]{Main} at $\theta=0$.}\label{tab:B-quantities}
\begin{tabular}{@{}lll@{}}
\toprule
quantity & definition & denominator\\
\midrule
full error & $u-u_h$ & $a\varepsilon+h^{1/6}\varepsilon^2$\\
ordering error & $\Cd_h=(u-u_h)-\varepsilon(v-v_h)$ & $h^{1/6}\varepsilon^2$\\
linear error & $\varepsilon(v-v_h)$ & $a\varepsilon$\\
total error (error of the linear discrete simulation) & $u-\varepsilon v_h$ & $a\varepsilon+\varepsilon^2$\\
continuous correction & $u-\varepsilon v$ & $\varepsilon^2$\\
discrete correction & $u_h-\varepsilon v_h$ & $\varepsilon^2$\\
discretization error of the quadratic Taylor coefficient & $\varepsilon^2(q-q_h)$ & $h^{1/2}t_\wedge^{-5/2}\varepsilon^2$\\
cubic remainder & $\varepsilon^3(r_\varepsilon-r_{h,\varepsilon})$ & $h^{1/6}\varepsilon^3$\\
\bottomrule
\end{tabular}
\end{table}

\paragraph{Norm, integrator, reference.}
The norm is that of $H$ (weight $1+x$), multiplied by $e^{\nu t}$ with
$\nu=3\pi^2/8$, half the rate $3\pi^2/4$. The time integration is by the
five-stage Radau~IIA method (order $9$). The reference is the W-scheme on
$N_{\rm ref}=1024,2048,4096$ ($2048,4096,8192$ for D4 at $N=128$). For
$t\ge10^{-4}$ the linear part is also compared with a continuous
eigenfunction expansion ($70$/$100$-digit coefficients, $32$/$48$/$64$
modes).

\paragraph{Predictions before computation.}
All predictions were written down and hashed before the computations they
concern. The files, hashes and the time stamps of the logs are kept. The
theorems are asymptotic estimates with unknown constants. The prediction
takes the form of the estimate (Table~\ref{tab:B-quantities}) and
determines the constant on coarse meshes. The empirical constant at a
point $(t,\varepsilon)$ is $C_{\rm emp}=e^{\nu t}\norm{\text{error}}_H/\text{denominator}$.
The \emph{frozen constant} $C_{\rm fr}$ is the maximum of $C_{\rm emp}$ over
the accepted points of the original data on $N=16,32$, without safety
factor and without later adjustment. For D4 these are its own meshes
$N=16,32$. The frozen constant is tested on $N=64,128,256$ and on
\emph{held-out data} on all meshes, through the ratio
$\rho=C_{\rm emp}/C_{\rm fr}$. A point exceeds the prediction if
$\rho>1+10^{-12}$. The constants were frozen on 2026-11-01, 10:07~UTC
(SHA-256 \texttt{9034f021\ldots}).

\paragraph{Acceptance.}
A point enters the constants and the tests only if time-step refinement
changes the value by at most $1\%$ and the reference check (two reference
meshes) by at most $10\%$, both with contraction. No threshold was
relaxed. Of $910{,}320$ records, $619{,}586$ ($68.1\%$) were accepted. The
time check failed at $75{,}011$ of the $685{,}504$ eligible records
($10.9\%$).

\paragraph{Data.}
D1: a smooth bump supported in $(0.25,0.75)$. D2: a function with a jump at
$x=0.43$. D3: a function with a kink at $x=0.43$. D4$_h$: a function
supported in one cell, constructed from the fourth derivative of a
compactly supported function divided by the weight, so that
$\Pi_h\mathrm{D4}_h=0$. Held-out data: D1$'$ with support $(0.20,0.65)$;
D2$'$, D3$'$ with the singularity at $x=0.37$. Meshes $N=16,\dots,256$. D4
up to $N=128$, because of cost. Amplitudes $\varepsilon\in\{1/32,1/8,1/2\}$.
The $235$ trajectories took $1.94\times10^5$~s ($54$~h) in total.

\paragraph{Uniform constants.}
Table~\ref{tab:B-frozen} gives the frozen constants of the full and the
ordering error. Table~\ref{tab:B-ratio} gives the maximum of $\rho$ over the
accepted points for each datum and mesh.

\begin{table}[H]
\centering\small
\caption{Frozen constants (maxima over accepted points on $N=16,32$).}\label{tab:B-frozen}
\begin{tabular}{@{}lcccc@{}}
\toprule
 & D1 & D2 & D3 & D4\\
\midrule
full error, $C_{\rm fr}$ & $1.45\times10^{-3}$ & $5.51\times10^{-2}$ & $2.02\times10^{-3}$ & $0.983$\\
ordering error, $C_{\rm fr}$ & $1.07\times10^{-5}$ & $1.21\times10^{-4}$ & $6.41\times10^{-6}$ & $1.46\times10^{-5}$\\
\bottomrule
\end{tabular}
\end{table}

\begin{table}[H]
\centering\small
\caption{Maximum of $\rho=C_{\rm emp}/C_{\rm fr}$ over accepted points
(three amplitudes, $t>0$). Values above $1$ are in bold. For D1--D3 the
value $1$ on $N=16$ or $32$ is the calibration point.}\label{tab:B-ratio}
\begin{tabular}{@{}l*{5}{c}@{\qquad}*{5}{c}@{}}
\toprule
 & \multicolumn{5}{c}{full error} & \multicolumn{5}{c}{ordering error}\\
\cmidrule(r){2-6}\cmidrule(l){7-11}
$N$ & 16 & 32 & 64 & 128 & 256 & 16 & 32 & 64 & 128 & 256\\
\midrule
D1 & 1 & .134 & .0174 & .0023 & .0006 & 1 & .116 & .0315 & .0088 & .0025\\
D1$'$ & \textbf{1.607} & .178 & .0240 & .0029 & .0007 & \textbf{1.393} & .141 & .0327 & .0091 & .0026\\
D2 & 1 & .765 & .530 & .398 & .276 & 1 & .235 & .168 & .0732 & .0285\\
D2$'$ & .939 & .695 & .450 & .354 & .230 & \textbf{1.016} & .276 & .104 & .0491 & .0157\\
D3 & 1 & .424 & .140 & .0562 & .0339 & 1 & .217 & .0371 & .0116 & .0039\\
D3$'$ & .968 & .412 & .163 & .0599 & .0326 & .908 & .210 & .0410 & .0091 & .0030\\
D4 & .998 & 1 & \textbf{1.0019} & \textbf{1.0036} & -- & 1 & .400 & .150 & .0576 & --\\
\bottomrule
\end{tabular}
\end{table}

The numerical results show that, apart from the entries in bold, the
constants fixed on $N=16,32$ suffice on the tested finer meshes and for the
held-out data. For each datum, the observed ordering-error ratio $\rho$
decreases over the tested refinements. For the six fixed data, it is at
most $0.029$ at $N=256$. For D4, which was tested only up to $N=128$, the
last value is $0.0576$. For the
held-out data every ordering prediction holds from $N=32$ on.

\paragraph{All exceedances.}
Table~\ref{tab:B-exceed} lists every (datum, mesh, quantity) with $\rho>1$,
including the auxiliary quantities of Table~\ref{tab:B-quantities}. They
mark the limits of the frozen constants. Two of them are structural. For
D4, $\Pi_h\mathrm{D4}_h=0$ gives $u_h\equiv0$. Hence the full error is
$\norm{u(t)}_H$. The datum changes with the mesh. The ratio $\rho$ stays
near~$1$. For D3$'$ the continuous correction does not depend on the mesh
and exceeds by the same factor $1.083$ on every mesh. This suggests that
the excess comes from the datum itself.

\begin{table}[H]
\centering\small
\caption{All exceedances ($\rho>1$): maximum $\rho$ and number of exceeding points.}\label{tab:B-exceed}
\begin{tabular}{@{}llll@{}}
\toprule
datum & $N$ & quantity & max $\rho$ (points)\\
\midrule
D1 & 64--256 & discrete correction & $1.00005$--$1.00006$ (1 each)\\
D1$'$ & 16 & linear error; full error & $1.607$ (225); $1.607$ (161)\\
 & 16 & ordering error; error of the linear discrete simulation & $1.393$ (162); $1.263$ (48)\\
D2$'$ & 16 & ordering; quadratic; cubic remainder & $1.016$ (18); $1.046$ (36); $1.009$ (12)\\
D3 & 64--256 & error of the linear discrete simulation & $1.006$, $1.010$, $1.010$ (3, 6, 6)\\
 & 64--256 & discrete correction & $1.00004$--$1.00005$ (2 each)\\
D3$'$ & 16--256 & continuous correction & $1.0827$ (106 on each mesh)\\
 & 16--256 & discrete correction & $1.0825$--$1.0827$ (106 on each mesh)\\
 & 16--256 & error of the linear discrete simulation & $1.064$, $1.085$, $1.091$, $1.092$, $1.092$ (9--10)\\
 & 16 & quadratic; cubic remainder & $1.067$ (51); $1.157$ (129)\\
D4 & 64, 128 & full error & $1.0019$ (184), $1.0036$ (122)\\
\bottomrule
\end{tabular}
\end{table}

\paragraph{Full error.}
From $N=32$ on, $\rho$ for the full error decreases for D1, D1$'$, D3, D3$'$
by factors between $1.7$ and $8.3$ per halving of $h$. For D2, D2$'$
(jump) it decreases by factors between $1.3$ and $1.5$. For the jump data
the full error stays closest to the form of the theorem. The numerical results show
no conflict with the form $a\varepsilon+h^{1/6}\varepsilon^2$
of~\cite[Theorem~4.2]{Main}.

\begin{figure}[H]
\centering
\includegraphics[width=\textwidth]{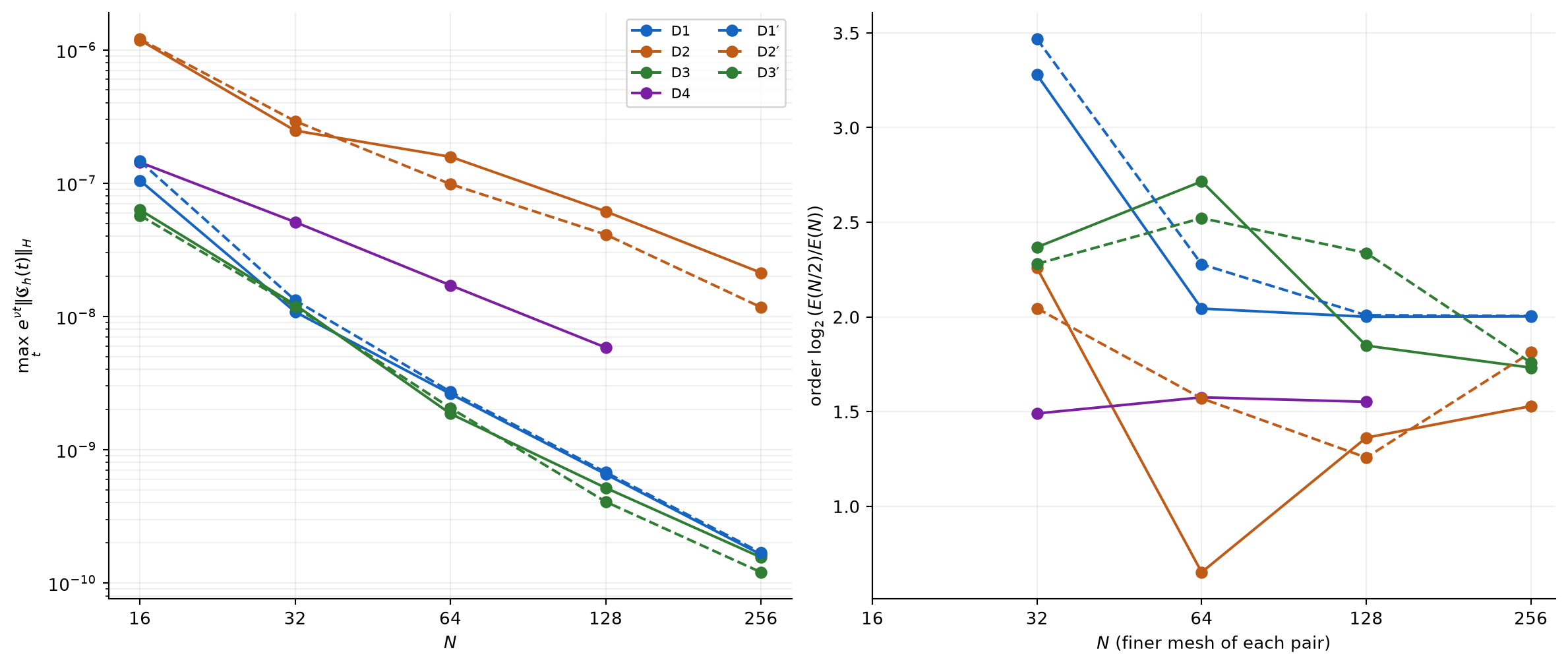}
\caption{W, $\varepsilon=1/8$. Left: maximum of $e^{\nu t}\norm{\Cd_h}_H$
over accepted points; dashed: held-out data. Right: orders
$\log_2(E(N/2)/E(N))$. Data: the accepted samples of this
subsection.}\label{fig:B-mesh}
\end{figure}

\begin{figure}[H]
\centering
\includegraphics[width=\textwidth]{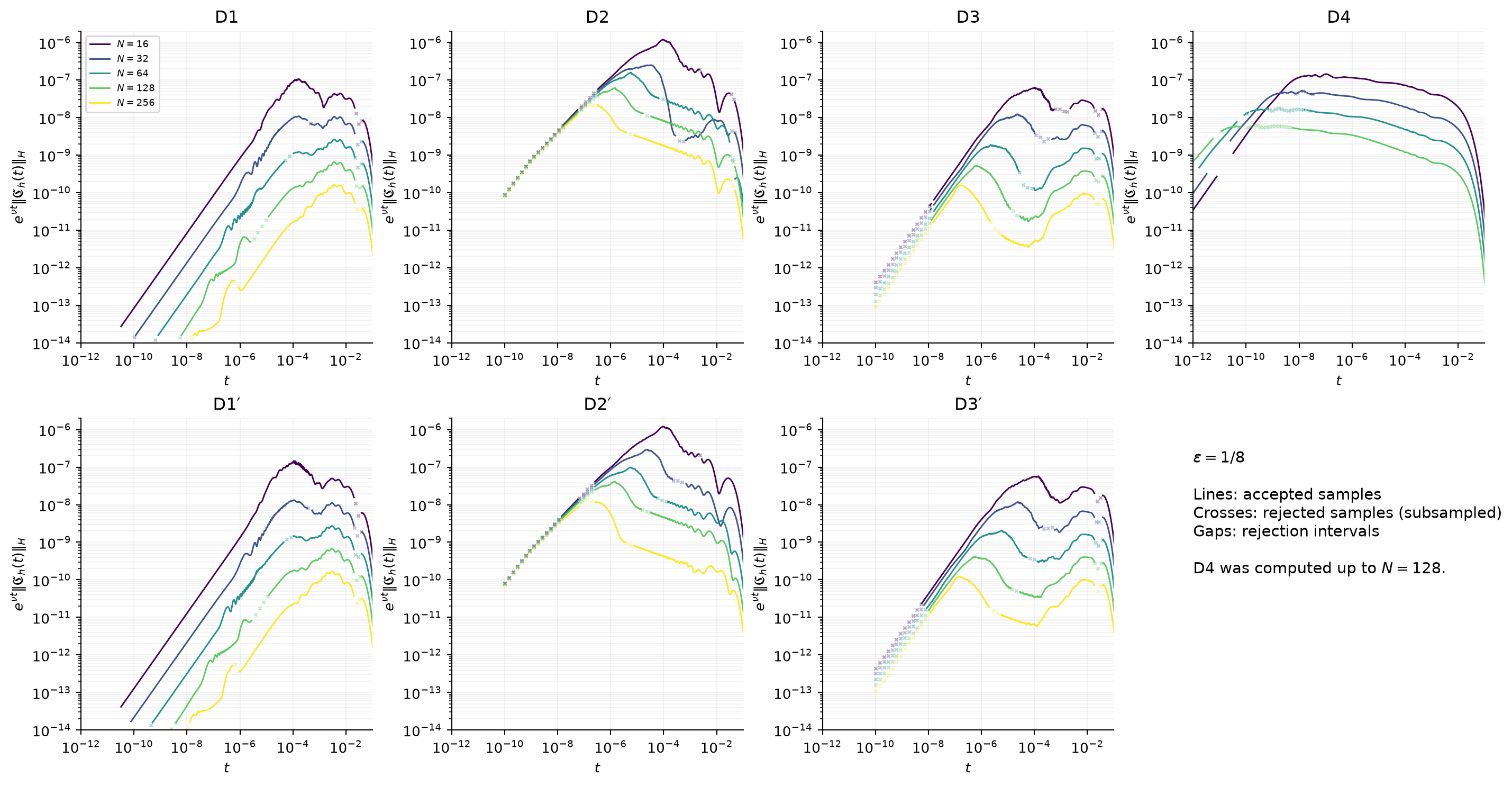}
\caption{W, $\varepsilon=1/8$: $e^{\nu t}\norm{\Cd_h(t)}_H$. Lines: accepted
samples; crosses: rejected samples (subsampled); gaps are rejection
intervals.}\label{fig:B-time}
\end{figure}

\paragraph{Ordering error.}
Figure~\ref{fig:B-mesh} shows the maximum of $e^{\nu t}\norm{\Cd_h}_H$ over
accepted points, $\varepsilon=1/8$, and the orders between neighbouring
meshes. The last orders are about $2$ (D1, D1$'$), $1.73$ (D3), $1.76$ (D3$'$),
and $1.49$, $1.58$, $1.55$ (D4, $N=32,64,128$). For D2 the orders are
$2.26$, $0.65$, $1.36$, $1.53$. For D2$'$ they are $2.05$, $1.57$,
$1.26$, $1.81$. The slowest step has order $0.65$. The order
of~\cite[Theorem~4.3]{Main} is $1/6$, uniformly over the energy ball.
There is no contradiction. Multiplying $\varepsilon$ by $4$ multiplies these
maxima by factors between $15.9966$ and $16.0127$, that is, by $16$ within
$0.08\%$. Figure~\ref{fig:B-time} shows the time profiles. The maximum
over time of the ordering error decreases under refinement for every
datum. Individual time curves of neighbouring meshes cross at some times
(for D2, $N=32$ and $64$ near $t=0.007$, $0.021$ and $0.027$).

\paragraph{Cross-checks.}
The finite-element references and the continuous expansion differ by at
most about $1.5\times10^{-7}$ in $H$. Radau~IIA of orders $9$ and $13$
differ by about $2\times10^{-12}$ to $5\times10^{-4}$ (relative) in four
representative cases. All trajectories have zero rejected steps.

\subsection{The two schemes on the spurious-mode datum}\label{app:fh}
The datum is $f_h$ of \eqref{eq:fh}. It lies in $X_h$. Hence both
schemes start from it without projection error. The norm is that of
$L^2(0,1)$. The $H$ values are recorded in the result files. The times are
$\mathcal C_St\in\{1/4,1,4\}$. The actual decay rates $-\Rea\zeta_h$ are
$1576.8$, $1576.1$, $1575.7$ for $N=512,1024,2048$. At these times
$e^{\nu t}\le1.01$. The decay factor is omitted. The S-scheme is evaluated by
exact exponential sums from its eigendecomposition up to third order in
$\varepsilon$. These formulas evaluate the computed Taylor coefficients.
The reported nonlinear errors of the S-scheme use the expansion truncated
after degree three in $\varepsilon$. P6 measures the effect of the
computed cubic term. It does not bound the omitted terms of degree four
and higher. The reported values should be read as numerical estimates.
The W-scheme is integrated by the seven-stage Radau~IIA method
(order $13$) up to fourth order, with references from the W-scheme on $2N$ and $4N$.
The S-scheme is not integrated in this way. The spurious mode has period about
$4\times10^{-10}$ at $N=512$. The stability function of Radau~IIA tends to
zero at infinity. Hence a stable step that does not resolve this period
damps the mode artificially, and the errors of the S-scheme would appear to
decrease. The predicted values were extrapolations of earlier computations with
$N\le256$. The theory gives only the direction for the linear error and
for the error of the quadratic coefficient. For the S-scheme these errors tend to
nonzero limits. For the W-scheme they tend to zero. The predictions were frozen on
2026-11-03, 05:48~UTC (SHA-256 \texttt{49C95C1A\ldots}). The first
computation finished at 05:58~UTC.

\begin{table}[H]
\centering\small
\caption{$f_h$, $\varepsilon=1/8$, $L^2$ norm. Ratio: ordering error of the
S-scheme divided by that of the W-scheme. $C_W$, $C_S$: constants
required by $\norm{\Cd_h}\le Ch^{1/6}\varepsilon^2$ for this datum and
time.}\label{tab:B-fh}
\setlength{\tabcolsep}{4pt}
\begin{tabular}{@{}cr*{7}{c}@{}}
\toprule
 & & \multicolumn{2}{c}{full error} & \multicolumn{2}{c}{ordering error} & & & \\
\cmidrule(lr){3-4}\cmidrule(lr){5-6}
$\mathcal C_St$ & $N$ & W-scheme & S-scheme & W-scheme & S-scheme & ratio & $C_W$ & $C_S$\\
\midrule
$1/4$ & 512 & $2.3\times10^{-9}$ & $9.733\times10^{-2}$ & $2.59\times10^{-9}$ & $4.634\times10^{-6}$ & $1.8\times10^{3}$ & $4.7\times10^{-7}$ & $8.39\times10^{-4}$\\
  & 1024 & $5.2\times10^{-9}$ & $9.734\times10^{-2}$ & $6.49\times10^{-10}$ & $4.634\times10^{-6}$ & $7.1\times10^{3}$ & $1.3\times10^{-7}$ & $9.42\times10^{-4}$\\
  & 2048 & $7.8\times10^{-10}$ & $9.734\times10^{-2}$ & $1.62\times10^{-10}$ & $4.634\times10^{-6}$ & $2.9\times10^{4}$ & $3.7\times10^{-8}$ & $1.06\times10^{-3}$\\[2pt]
$1$ & 512 & $2.2\times10^{-9}$ & $4.594\times10^{-2}$ & $2.31\times10^{-9}$ & $9.590\times10^{-6}$ & $4.2\times10^{3}$ & $4.2\times10^{-7}$ & $1.74\times10^{-3}$\\
  & 1024 & $2.0\times10^{-9}$ & $4.596\times10^{-2}$ & $5.78\times10^{-10}$ & $9.592\times10^{-6}$ & $1.7\times10^{4}$ & $1.2\times10^{-7}$ & $1.95\times10^{-3}$\\
  & 2048 & $3.2\times10^{-10}$ & $4.597\times10^{-2}$ & $1.44\times10^{-10}$ & $9.593\times10^{-6}$ & $6.6\times10^{4}$ & $3.3\times10^{-8}$ & $2.19\times10^{-3}$\\[2pt]
$4$ & 512 & $1.8\times10^{-9}$ & $2.281\times10^{-3}$ & $1.92\times10^{-9}$ & $9.361\times10^{-6}$ & $4.9\times10^{3}$ & $3.5\times10^{-7}$ & $1.69\times10^{-3}$\\
  & 1024 & $1.1\times10^{-9}$ & $2.285\times10^{-3}$ & $4.79\times10^{-10}$ & $9.367\times10^{-6}$ & $2.0\times10^{4}$ & $9.7\times10^{-8}$ & $1.90\times10^{-3}$\\
  & 2048 & $2.0\times10^{-10}$ & $2.287\times10^{-3}$ & $1.20\times10^{-10}$ & $9.370\times10^{-6}$ & $7.8\times10^{4}$ & $2.7\times10^{-8}$ & $2.14\times10^{-3}$\\
\bottomrule
\end{tabular}
\end{table}

\paragraph{Uniform constants and ordering error.}
The numerical results show that, at all three times, the ordering error of
the W-scheme is divided by about $4$ per halving of $h$. Hence $C_W$
decreases. In contrast, the ordering error of the S-scheme changes by less than
$0.1\%$. Hence $C_S$ grows by about $2^{1/6}$ per halving. The ratio of the
two ordering errors grows by
a factor of about $4$ per halving and is $6.6\times10^{4}$ at $N=2048$,
$\mathcal C_St=1$. The two references differ by up to $5.9\%$ in the W-scheme
values (P8 in Table~\ref{tab:B-pred}). This difference indicates the
uncertainty of the W-scheme values and of the ratio. The numerical results
suggest that the S-scheme has no uniform
constant of the form of~\cite[Theorem~4.3]{Main}. Statements for all
meshes can only come from theorems.

\paragraph{Full error.}
The numerical results show that the full error of the S-scheme is nearly
constant over the three tested meshes. It differs from
$\varepsilon$ times its linear error by less than $10^{-5}$ (relative). It is
carried by the spurious mode. The full error of the W-scheme is at most
$5.2\times10^{-9}$. At these times $ht_\wedge^{-3/2}>1$ on all three
meshes. Hence~\cite[Theorem~4.2]{Main} bounds the linear part of the full
error only by $C\varepsilon$. Table~\ref{tab:B-fh} therefore compares sizes,
not rates. The full error of the W-scheme also varies irregularly with $N$, because
its linear part is about the size of the continuous solution itself, which
depends on $f_h$ (Table~\ref{tab:B-linear}).

\paragraph{Amplitude.}
At $N=2048$, $\mathcal C_St=1$, the ordering errors are, for
$\varepsilon=1/32,1/8,1/2$: S-scheme $6.0\times10^{-7}$, $9.6\times10^{-6}$,
$1.5\times10^{-4}$; W-scheme $9.0\times10^{-12}$, $1.4\times10^{-10}$,
$2.3\times10^{-9}$. The numerical results show that both scale as
$\varepsilon^2$ within $3\%$ at the three computed amplitudes. At these
amplitudes, the ratio of the ordering errors remains close to
$6.6\times10^4$. For the W-scheme at $N=512$, $\varepsilon=1/2$,
the direct nonlinear solution and the fourth-order expansion differ by
$3\times10^{-7}$ times the ordering error of the W-scheme.

\paragraph{Linear errors.}
Table~\ref{tab:B-linear} gives the linear errors. The numerical results
show that the linear error of the S-scheme is within $0.4\%$ of
$e^{-\mathcal C_St}$ on the computed meshes. This is consistent with
Theorem~\ref{thm:total}(a). For the W-scheme, the computed linear error is
within $5\%$ of the norm of the reference solution. These are errors on a
family of data $f_h$ that depends on the mesh. The bound
of~\cite[Theorem~4.2]{Main} implies convergence to zero at each fixed
positive time. The three reported values do not show a regular refinement
rate.

\begin{table}[H]
\centering\small
\caption{Linear errors $\norm{v-v_h}$ and the norm of the reference solution $\norm{v}$, $L^2$.}\label{tab:B-linear}
\footnotesize\setlength{\tabcolsep}{3.5pt}
\begin{tabular}{@{}r*{3}{ccc}@{}}
\toprule
 & \multicolumn{3}{c}{$\mathcal C_St=1/4$} & \multicolumn{3}{c}{$\mathcal C_St=1$} & \multicolumn{3}{c}{$\mathcal C_St=4$}\\
\cmidrule(lr){2-4}\cmidrule(lr){5-7}\cmidrule(l){8-10}
$N$ & W-scheme & S-scheme & $\norm v$ & W-scheme & S-scheme & $\norm v$ & W-scheme & S-scheme & $\norm v$\\
\midrule
512 & $7.2\times10^{-9}$ & $0.7786$ & $6.9\times10^{-9}$ & $2.6\times10^{-9}$ & $0.3676$ & $2.5\times10^{-9}$ & $9.0\times10^{-10}$ & $0.01825$ & $8.6\times10^{-10}$\\
1024 & $3.9\times10^{-8}$ & $0.7787$ & $3.9\times10^{-8}$ & $1.4\times10^{-8}$ & $0.3677$ & $1.4\times10^{-8}$ & $4.8\times10^{-9}$ & $0.01828$ & $4.8\times10^{-9}$\\
2048 & $5.5\times10^{-9}$ & $0.7788$ & $5.5\times10^{-9}$ & $2.0\times10^{-9}$ & $0.3678$ & $2.0\times10^{-9}$ & $6.8\times10^{-10}$ & $0.01830$ & $6.8\times10^{-10}$\\
\bottomrule
\end{tabular}
\end{table}

\paragraph{Check of the predictions.}
Table~\ref{tab:B-pred} lists the eleven predictions with their outcomes.
Ten pass on all three meshes. So do the parts across meshes. P10 fails
in one item on every mesh. The item substitutes the computed quadratic
coefficient of the S-scheme back into its equation. The threshold was $10^{-6}$. The
residual is $1.7\times10^{-6}$, $1.5\times10^{-5}$, $1.1\times10^{-3}$ for
$N=512,1024,2048$. A diagnosis made afterwards shows that it is of the size
of the estimated double-precision rounding error. This size was not
foreseen when the threshold was fixed. The threshold was not changed. A
second eigendecomposition (QZ instead of Cholesky) changes the quadratic
coefficient by $1.7\times10^{-4}$, $1.8\times10^{-5}$, $6.3\times10^{-7}$
(relative) at $\mathcal C_St=1/4,1,4$, $N=2048$. Thus the two
eigendecompositions agree to about four digits at $\mathcal C_St=1,4$ and
to three at $\mathcal C_St=1/4$. This agreement is a cross-check, not an
error bound. The conclusions drawn from Table~\ref{tab:B-fh} need only its
first two digits. The W-scheme values do not use this eigendecomposition.

\begin{table}[H]
\centering\small
\caption{The frozen predictions of the $f_h$ comparison, $N=512,1024,2048$.}\label{tab:B-pred}
\begin{tabular}{@{}l>{\raggedright\arraybackslash}p{0.52\textwidth}>{\raggedright\arraybackslash}p{0.3\textwidth}@{}}
\toprule
 & prediction & outcome\\
\midrule
P1 & linear error of the S-scheme is $0.779$, $0.368$, $0.0183$ on all meshes, within $1\%$ & pass\\
P2 & discretization error of the quadratic Taylor coefficient of the S-scheme is $(3.0,6.1,6.0)\times10^{-4}$, within $2\%$ & pass\\
P3 & ratio of linear errors, S-scheme/W-scheme, $\ge10^5$ & pass; min $3.8\times10^6$\\
P4 & discretization error of the quadratic Taylor coefficient of the W-scheme at least halved per halving of $h$; $\le2\times10^{-8}$ at $N=2048$ & pass; factor about $1/4$; max $1.04\times10^{-8}$\\
P5 & ratio of ordering errors, S-scheme/W-scheme, $\ge10^3$ for $\varepsilon=1/32,1/8$; growing in $N$ at $\mathcal C_St=1,4$ & pass; min $1.8\times10^3$; growth $\approx4$ per halving\\
P6 & S-scheme: third-order term changes the ordering error by $<10\%$ ($\varepsilon=1/32,1/8$) & pass; max $5\times10^{-5}$\\
P7 & W-scheme: fourth-order term $\le1\%$ of the ordering error; direct nonlinear check $\le1\%$ ($N=512$, $\varepsilon=1/2$) & pass; $1.5\times10^{-10}$, $3\times10^{-7}$\\
P8 & reference check ($2N$ against $4N$) $\le10\%$ for the W-scheme quantities & pass; max $5.9\%$\\
P9 & time check (rtol $10^{-9}$ against $10^{-11}$) $\le1\%$ for the W-scheme quantities & pass; max $1.8\times10^{-5}$\\
P10 & S-scheme self-check: eigenpair residual $\le10^{-10}$; bounds of Theorem~\ref{thm:top}; dense against sparse $\varphi_h$ $\le10^{-8}$; residuals of $q$, $p_3$ $\le10^{-6}$ & \textbf{fail} in the $q$ residual; other items pass\\
P11 & full error of the S-scheme within $1\%$ of $\varepsilon$ times its linear error ($\varepsilon\le1/8$) & pass; max $8\times10^{-6}$\\
\bottomrule
\end{tabular}
\end{table}

\paragraph{Computation.}
Twelve jobs, four in parallel, took about $16$~h of wall-clock time. The
longest was the reference on $8192$ cells for $N=2048$.

\paragraph{Scope.}
These computations compare the two schemes on one mesh-dependent datum.
They quantify the gap between the theorems with and without mesh-uniform
regularity. They give no ranking of the schemes for fixed data or for
computational cost (Section~\ref{app:further}).

\subsection{Further computations}\label{app:further}

\paragraph{The limit of the quadratic coefficient.}
The candidate $E$ of Proposition~\ref{thm:quad}, with the term
$\mathfrak b\mathfrak p$, was found after a simpler candidate without this
term failed at $N\le256$. Its values at $N=512,1024$, and the analogue for
the second mode, were recorded before those computations. At
$\mathcal C_St=1/4,1,4$,
$\norm{E(t)}_2=(2.966,6.140,5.998)\times10^{-4}$. At the same times,
$\norm{\mathring P_{2,h}(\tilde f_h)(t)-E(t)}_2$ is
$(1.36,2.02,2.29)\times10^{-6}$ for $N=256$. It is
$(6.75,9.89,11.1)\times10^{-7}$ for $N=512$, and
$(3.36,4.89,5.48)\times10^{-7}$ for $N=1024$. These values are
consistent with first-order convergence.

\paragraph{Fixed data.}
These computations are a posteriori and concern the linear flow on
$N=32$--$256$, with the W-scheme on $N=1024$ or $2048$ as reference. On a smooth datum
and on a datum with a kink, the S-scheme was more accurate than the W-scheme in most cases. For
$t\ge0.003$ the ratio of the errors, S-scheme/W-scheme, lay between $0.0026$ and $0.30$. The
error of the S-scheme behaved like $h^3$ on $N\le128$. At $N=256$ it is at the level
of the reference accuracy. On a datum with a jump, at $t=0.001$, this ratio
was $11$, $22$, $46$ for $N=64,128,256$ (the S-scheme about $h$, the
W-scheme about $h^2$).
For $t\ge0.01$ and $N\ge64$ it lay between $0.18$ and $4.4$. The numerical
results show that at $t=0.001$ almost all of the error of the S-scheme on the jump
datum is carried by the computed spurious modes (for $N=64$,
$1.310\times10^{-3}$ of $1.319\times10^{-3}$). The amplitude of these modes
is roughly proportional to $h$. Their decay rates are close to
$\mathcal C_Sn^2$. For the W-scheme the modes at the top of the spectrum are damped at
rates of about $31h^{-2}$. In these computations, large differences on
fixed data appeared only at short times and for the jump datum.

\paragraph{The coarse-mesh correction.}
The examples of Section~\ref{sec:twomesh} use the stored W-scheme trajectories of
Section~\ref{app:S2}, without new time integration. They use the fine mesh
$N=1024$, coarse meshes $N=16,32,64$, and reference $N=4096$. The
difference of the references on $2048$ and $4096$ elements is
$7.2\times10^{-10}$ to $1.3\times10^{-9}$ at $t=0.01$.
Table~\ref{tab:twomesh} shows a part of the output.\label{app:twomesh}

\paragraph{Running times.}
For D2 and $\varepsilon=1/8$ we computed the four solutions $\varepsilon v_h$,
$\varepsilon v_{\hat h}$, $u_{\hat h}$ and $u_h$ directly, with $N=1024$ and
$64$. We used the five-stage Radau~IIA method. The nonlinear stage
equations were solved by fixed-point iteration with relative tolerance
$3\times10^{-11}$. The same rule generated the time grids for all four
solutions. Its dependence on $N$ was confined to the initial layer. The
fine grid had $96$ more steps there, giving $33{,}800$ steps against
$33{,}704$ up to $t=0.1$. The errors of
$u_{\hat h}$, $\tilde u$ and $u_h$ agree with Table~\ref{tab:twomesh} to
three digits. Table~\ref{tab:timing} gives the times and counts up to
$t=0.1$. Serial times are those of the faster of two runs, measured after
the assembly of the matrices, on one thread. The parallel time is the
wall-clock time of three simultaneous processes, one for each part of
$\tilde u$, from their start to the last result, including process start.
It comes from a separate run. The machine has eight cores (Intel64 Family 6
Model 140, Windows 11, Python 3.12). Five predictions were written in the
programs before the runs:
\begin{enumerate}
\item the errors agree with the stored ones within $10\%$;
\item the time of $u_h$ is at least three times that of $\varepsilon v_h$;
\item the serial time of $\tilde u$ is at most half of that of $u_h$;
\item the parallel time is within $25\%$ of the longest serial part;
\item the parallel time is at most a third of the time of $u_h$.
\end{enumerate}
All five hold. The third holds with a small margin
($95.0$ against $98.6$ seconds at $t=0.01$).

\begin{table}[H]
\centering\small
\caption{Running times and counts up to $t=0.1$, D2, $\varepsilon=1/8$, W-scheme.
Counts: time steps, LU factorizations, LU solves and evaluations of the
nonlinear term.}\label{tab:timing}
\begin{tabular}{@{}lrrrrrr@{}}
\toprule
solution & unknowns & steps & LU fact. & LU solves & nonlinear eval. & seconds\\
\midrule
$\varepsilon v_h$, $N=1024$ & $2047$ & $33{,}800$ & $12{,}276$ & $101{,}400$ & $0$ & $42.4$\\
$\varepsilon v_{\hat h}$, $N=64$ & $127$ & $33{,}704$ & $12{,}228$ & $101{,}112$ & $0$ & $7.4$\\
$u_{\hat h}$, $N=64$ & $127$ & $33{,}704$ & $12{,}228$ & $301{,}815$ & $334{,}505$ & $49.0$\\
$u_h$, $N=1024$ & $2047$ & $33{,}800$ & $12{,}276$ & $303{,}216$ & $336{,}360$ & $207.0$\\
\midrule
$\tilde u$, three parts one after another & & & & & & $98.8$\\
$\tilde u$, three parts simultaneously (wall clock) & & & & & & $47.2$\\
\bottomrule
\end{tabular}
\end{table}

In the parallel run the three parts took $46.1$, $6.7$ and $37.9$ seconds.
Measured as the sum of the times of the processors, $\tilde u$ took
$91$--$99$ seconds against $207$ for $u_h$. Measured as wall-clock time with
three cores, it took $47$ seconds against $207$.

\begin{table}[H]
\centering\small
\caption{Errors in $H$ of the linear simulation $\varepsilon v_h$ ($N=1024$),
of the nonlinear solution on the coarse mesh $u_{\hat h}$, of the combined
solution $\tilde u$ ($N=1024$ and the coarse mesh) and of the nonlinear
solution $u_h$ ($N=1024$); $\varepsilon=1/8$.}\label{tab:twomesh}
\begin{tabular}{@{}llrcccc@{}}
\toprule
datum & $t$ & coarse $N$ & $\varepsilon v_h$ & $u_{\hat h}$ & $\tilde u$ & $u_h$\\
\midrule
D1 & $0.01$ & 16 & $8.19\times10^{-5}$ & $1.03\times10^{-5}$ & $2.84\times10^{-8}$ & $2.17\times10^{-9}$\\
   &        & 64 & $8.19\times10^{-5}$ & $6.37\times10^{-7}$ & $1.74\times10^{-9}$ & $2.17\times10^{-9}$\\
D2 & $0.01$ & 16 & $1.21\times10^{-4}$ & $4.08\times10^{-5}$ & $2.97\times10^{-8}$ & $4.18\times10^{-9}$\\
   &        & 64 & $1.21\times10^{-4}$ & $2.76\times10^{-6}$ & $6.29\times10^{-9}$ & $4.18\times10^{-9}$\\
D2 & $0.1$  & 64 & $3.10\times10^{-7}$ & $1.37\times10^{-8}$ & $1.04\times10^{-10}$ & $6.64\times10^{-11}$\\
\bottomrule
\end{tabular}
\end{table}

\subsection{Limits}
The computations cover finitely many data, meshes and times. Exact
normalized errors at these points would give lower bounds for the optimal
constants in the estimates. The reported values are numerical
approximations to those errors. They provide empirical estimates of the
constants required by the sampled cases. No rigorous bound for the optimal
constants is inferred from them. Only the base norm ($\theta=0$) was
computed. The programs, raw data, logs and hashes are kept for
reproduction.

%% file: appZ_plain.tex
\setcounter{section}{25}
\section{The computations for readers outside numerical analysis}\label{app:plain}

This appendix explains the
design of the computations of Appendix~\ref{app:num} and what they show,
for readers who have not worked in numerical analysis. The numbers are
presented as they are. The judgement is left to the reader. We write
``the numerical results show'' for what can be read off directly, and
``the numerical results suggest'' for what involves an inference.

\subsection{Summary}
All conclusions rest on the theorems. The computations prove nothing. They
do two things. On the computed data and meshes, they give the actual
numbers for the errors and for how they change under refinement. On one
and the same datum, they
measure how far apart the two schemes are, one with and one without
mesh-uniform regularity.

The computations are designed so that a reader can judge the numbers
without checking the programs. Sections~\ref{app:S2} and~\ref{app:fh}
report predictions recorded before the computations they concern.
Section~\ref{app:S2} fixes constants on coarse meshes and tests them on
finer meshes and held-out data. It admits samples only after time-step and
reference checks. Section~\ref{app:fh} uses a separate set of eleven
predictions. Its failed self-check P10 is retained and reported. The
additional computations have the status stated in
Section~\ref{app:further}.

\subsection{What is compared}
The equation is \eqref{eq:pde} with $L=1$. Both schemes are Galerkin
approximations. They use the same trial space: functions that are cubic on
each of $N$ equal subintervals, continuously differentiable, and zero at
both ends together with their slope at the right end, so that the space
has dimension $2N-1$. The approximate solution satisfies a
system of $2N-1$ ordinary differential equations. The mesh size is
$h=1/N$. The two schemes differ only in the inner product used for testing.
\begin{itemize}
\item \emph{The weighted scheme (W-scheme).} The inner product carries the weight
$L+x$, Rosier's multiplier, which gives exponential stability and Kato
smoothing. Our previous paper~\cite{Main} proves that the discrete semigroups of
the W-scheme inherit exponential stability, Crandall--Pazy regularity and
Kato smoothing, with bounds uniform in the mesh size. As a consequence, its
error estimates hold uniformly for all data in a ball.
\item \emph{The unweighted scheme (S-scheme).} It uses the unweighted $L^2$ inner
product. For $N\ge512$ it has a high-frequency spurious mode (Theorem~\ref{thm:top}(a),
whose proof is checked by Lean~4). This mode is an exact eigenpair of the
discrete system whose frequency is at least $100N^3$ and whose decay rate
stays bounded as $N$ grows. The continuous problem has no such mode (see
Section~\ref{sec:mur}). Hence the S-scheme loses
mesh-uniform regularity.
\end{itemize}

\subsection{What is measured}
Let the datum be $\varepsilon f$, where $\varepsilon$ is the amplitude. The
\emph{full error} is the difference between the nonlinear Galerkin
solution and the exact solution. The \emph{linear error} is the same difference for the
linearized problem. A small nonlinear wave can be approximated by a linear
discrete simulation in two orders: linearize first and then discretize,
or discretize first and then linearize. Both orders lead to the same
approximation, the linear discrete solution, and to the same
\emph{total error}, the difference between it and the exact solution.
They divide this error
differently between linearization and discretization. The difference of
the two divisions is the \emph{ordering error}, a notion introduced
in~\cite{Main}. When it is small, the nonlinear part can be computed on a
coarser mesh (Section~\ref{sec:twomesh}). Our previous paper~\cite{Main}
bounds it for the W-scheme by $Ch^{1/6}\varepsilon^2$. The bound holds on a whole
energy ball and for all time. No higher norm of the datum enters. The
exact solution is not available in closed form. It is replaced by
a W-scheme solution on a much finer mesh, called the reference solution.

\subsection{Design of the computations}
Reading the programs line by line is not needed to follow the design of
the computations. The computations were designed so that readers can judge
for themselves how far the numbers can be trusted.
\begin{enumerate}
\item \emph{Predictions first.} Before the computations of
Sections~\ref{app:S2} and~\ref{app:fh}, the predicted values and the
pass/fail thresholds were recorded. The SHA-256
hash of the file was recorded. The file was not changed afterwards. Each
prediction was then checked. The outcome is reported whether it passed
or failed. The prediction files, their hashes and the time stamps of the
logs are kept so that the order can be checked. The theorems are asymptotic
estimates whose constants are not known numerically. Hence the values in the
predictions can only come from computations on coarser meshes. For the
seven data of Section~\ref{app:S2}, the form of the estimate, that is, its
dependence on $h$, $\varepsilon$ and $t$, is taken from the theorem. The
constants were determined on $N=16,32$. For the comparison on $f_h$
(Section~\ref{app:fh}), the predicted values were extrapolated from earlier
computations with $N\le256$. The theory gives only the direction for the
linear error and for the error of the quadratic coefficient. For the S-scheme these
errors tend to nonzero limits. For the W-scheme they tend to zero.
\item \emph{Frozen constants, tested on new cases.} The following
calibration procedure applies to Section~\ref{app:S2}. The constants were
determined on the two coarsest meshes only. Each constant is the largest
value observed there, with no safety factor and no later adjustment. They
were then tested on finer meshes ($N=64,128,256$) and on held-out data.
This is like calibrating an instrument on known samples and then measuring
new ones. It avoids fitting after the fact.
\item \emph{Acceptance of samples in Section~\ref{app:S2}.} A sample is
used for calibration and testing only if the time-step check changes its
value by at most $1\%$ and the reference check by at most $10\%$. The
checks also require contraction under refinement. No threshold was
relaxed. Section~\ref{app:fh} uses the checks listed in
Table~\ref{tab:B-pred}.
\item \emph{Independent cross-checks.} For example, the errors of the S-scheme are
computed from exact formulas based on its eigendecomposition. For the W-scheme at
$N=512$, the expansion in the amplitude was compared with a direct nonlinear
time integration (Section~\ref{app:fh}).
\item \emph{The tests could fail, and every failed entry is listed.} The
predictions fixed on coarse meshes were tested on fine meshes. The
time stamps, hashes and cross-checks are kept. These facts do not exclude
every possible error in the programs. However, they provide checks against
fitting after the fact and against errors that change with the mesh.
Failures did occur. For the seven data, the
held-out data exceed the predictions on the coarsest mesh $N=16$. In the
comparison on $f_h$, one self-check of the S-scheme computation failed. As $N$
grows, the residual of the eigendecomposition method grows and is of the
size of the estimated rounding error. These entries are kept and
reported.
\item \emph{Reproducibility.} Programs, raw data, logs and file hashes are
kept.
\end{enumerate}

\subsection{Result 1: how large the errors of the W-scheme actually are}
The data range from a smooth one to one that the scheme cannot see at all.
They are a smooth bump (D1), a function with a jump (D2), a function with a
kink (D3), and a bump supported in a single mesh cell (D4). The projection
of D4 onto the trial space is zero, so that its whole ordering error comes
from the continuous problem. In addition there are three held-out data.

Figure~\ref{fig:B-mesh} has two panels. In the left panel, each line is
one datum. The solid lines are D1--D4, and the dashed lines are the
held-out data. The horizontal axis is the number of elements, doubled from
one mark to the next. The vertical axis is the largest weighted ordering
error over time, on a logarithmic scale. A steeper line means a faster
decrease. The right panel gives the order between neighbouring meshes.
Order $2$ means division by $4$ per halving. The order $1/6$ of the
theorem lies below all points.

\emph{Decrease under refinement} (Figure~\ref{fig:B-mesh}). Halving the
mesh size divides the ordering error
\begin{itemize}
\item by $4$ for the smooth datum D1 (rate $h^2$);
\item by about $1.6$ to $4.8$ for the jump data (orders $0.65$ to
$2.26$). The factor varies between refinement steps. The smallest, about
$1.6$, occurs for D2 from $N=32$ to $64$;
\item by about $3.4$ for the kink D3 (about $h^{1.75}$);
\item by about $2.8$ for the single cell D4 (about $h^{1.5}$).
\end{itemize}
The bound of the
theorem is $h^{1/6}$, a factor $1.12$ per halving. The numerical results
show that, for single data on the computed meshes, the ordering error
decreases much faster than
the rate of the theorem. The theorem holds uniformly over all data in a
ball. There is no contradiction.

\emph{Amplitude.} For every datum and every mesh, multiplying
$\varepsilon$ by $4$ multiplies the ordering error by about $16$, to within
$0.08\%$. The numerical results show that, at the three computed
amplitudes, the ordering error is close to quadratic in the amplitude.

\emph{The frozen constants.} The ratio of the observed ordering error to
the predicted one decreases under refinement. For the six fixed data, it
is at most $0.03$ at $N=256$. D4 was tested only up to $N=128$, where this
ratio is about $0.058$. Predictions are exceeded on the coarsest mesh, where
each held-out datum exceeds in some quantity, by at most a factor of
about $1.6$. On finer meshes, three auxiliary quantities of the
kink at a different position exceed by up to $9.2\%$, and a few further
quantities by at most about $1\%$ (Table~\ref{tab:B-exceed}).

\emph{Time histories} (Figure~\ref{fig:B-time}). For every datum the
largest ordering error over time decreases under refinement. At some times
the curves of neighbouring meshes cross.

\subsection{Result 2: the W-scheme and the S-scheme on the same datum}
The datum is the real part of the eigenfunction of the spurious mode of the S-scheme, normalized, denoted
$f_h$. Since $f_h$ lies in the trial space of both schemes, they start from
exactly the same datum. The datum $f_h$ depends on $h$. The time is
$\mathcal C_St=1$, that is, $t\approx6.3\times10^{-4}$, about the time at
which the spurious mode has decayed by the factor $1/e$. The meshes are
$N=512,1024,2048$. The norm is that of $L^2$.

\emph{Ordering error} ($\varepsilon=1/8$). The values are in
Table~\ref{tab:intro-numbers}. The numerical results show that halving the
mesh size divides the ordering error of the W-scheme by about $4$. In contrast, the
ordering error of the S-scheme changes by less than $0.1\%$. Their ratio grows by a factor of
about $4$ per halving, from $4.2\times10^3$ to $6.6\times10^4$. The same
holds at $\mathcal C_St=1/4$ and $4$. The W-scheme values are uncertain by up to
about $6\%$, because the reference solutions on two finer meshes differ
by that much.

\emph{Full error.} At $\mathcal C_St=1$, the numerical results show that the
full error of the S-scheme stays close to $4.6\times10^{-2}$ over the three
tested meshes. It comes almost entirely from the spurious mode: it differs from
$\varepsilon$ times the linear error of the S-scheme on the same mesh (about $0.368$)
by less than $10^{-5}$ (relative). The full error of the W-scheme is at
most $2.2\times10^{-9}$ in these computations. It contains the linear part, which is
about the size of the reference solution itself and changes with $f_h$.
Hence the full error of the W-scheme varies irregularly with $N$. What decreases
regularly is the ordering error.

\emph{Uniform constants.} Writing the ordering error in the form
$Ch^{1/6}\varepsilon^2$ of the theorem of~\cite{Main} gives the constant
$C$ required by this datum at this time. For the W-scheme it is
$(4.2,1.2,0.33)\times10^{-7}$. For the S-scheme it is
$(1.74,1.95,2.19)\times10^{-3}$. The constant
required by the S-scheme grows by about $2^{1/6}$ per halving, which is another way of
saying that the ordering error of the S-scheme does not decrease. The numerical
results suggest that the S-scheme has no uniform constant of this form. Statements
for all meshes can only come from the theorems.

\emph{Amplitude} ($N=2048$). Multiplying $\varepsilon$ by $4$ multiplies
the ordering errors of both schemes by about $16$. At the three computed
amplitudes their ratio changes very little.

\emph{Predictions.} Of the eleven predictions, ten pass on all three meshes.
One item of the self-check of the S-scheme computation fails. When the computed
quadratic coefficient of the S-scheme is substituted back into its equation, the
residual should be at most $10^{-6}$. It grows to about $1.1\times10^{-3}$
at $N=2048$. A diagnosis after the fact shows that it is of the size of the
rounding error. The threshold was not changed. A recomputation with a
second eigendecomposition agrees with the first to about four digits
(three at the earliest time). This agreement is a check, not a guarantee.
The conclusions need only the first two digits. The W-scheme values are
obtained by direct time integration and are not affected.

\subsection{Limits}
Numbers cannot replace proofs. Statements for all meshes and all data rest
on the theorems. Each tested case gives a numerical estimate of the
constant needed for that case. These estimates do not determine the
largest value over the whole data ball. Their numerical uncertainty also
prevents us from treating them as certified lower bounds. Only the base
norm was computed.

%% file: kdv_mur_gap.bbl
\begin{thebibliography}{99}

\bibitem{BanksItoWang}
H.~T. Banks, K.~Ito and C.~Wang,
\newblock Exponentially stable approximations of weakly damped wave equations,
\newblock in: {\em Estimation and Control of Distributed Parameter Systems},
Internat. Ser. Numer. Math. 100, Birkh\"auser, Basel, 1991, 1--33.

\bibitem{DWW}
C.~N. Dawson, M.~F. Wheeler and C.~S. Woodward,
\newblock A two-grid finite difference scheme for nonlinear parabolic
equations,
\newblock {\em SIAM J. Numer. Anal.},
35: 435--452, 1998.

\bibitem{ErvedozaZuazua}
S. Ervedoza and E. Zuazua,
\newblock Uniformly exponentially stable approximations for a class of
damped systems,
\newblock {\em J. Math. Pures Appl.},
91: 20--48, 2009.

\bibitem{Glowinski92}
R. Glowinski,
\newblock Ensuring well-posedness by analogy; {S}tokes problem and boundary
control for the wave equation,
\newblock {\em J. Comput. Phys.},
103: 189--221, 1992.

\bibitem{IgnatZuazua}
L.~I. Ignat and E. Zuazua,
\newblock Convergence of a two-grid algorithm for the control of the wave
equation,
\newblock {\em J. Eur. Math. Soc.},
11: 351--391, 2009.

\bibitem{InfanteZuazua}
J.~A. Infante and E. Zuazua,
\newblock Boundary observability for the space semi-discretizations of
the 1-D wave equation,
\newblock {\em M2AN Math. Model. Numer. Anal.},
33: 407--438, 1999.

\bibitem{Main}
K.~Liu,
\newblock Linearize or discretize first? Error estimates under mesh-uniform
regularity and exponential stability for the KdV equation with boundary
dissipation,
\newblock arXiv:2610.39269v2, 2026.

\bibitem{MarionXu}
M. Marion and J. Xu,
\newblock Error estimates on a new nonlinear {G}alerkin method based on
two-grid finite elements,
\newblock {\em SIAM J. Numer. Anal.},
32: 1170--1184, 1995.

\bibitem{RTT}
K. Ramdani, T. Takahashi and M. Tucsnak,
\newblock Uniformly exponentially stable approximations for a class of
second order evolution equations. Application to LQR problems,
\newblock {\em ESAIM Control Optim. Calc. Var.},
13: 503--527, 2007.

\bibitem{TTZ}
L.~R. Tcheugou\'e T\'ebou and E. Zuazua,
\newblock Uniform exponential long time decay for the space
semi-discretization of a locally damped wave equation via an artificial
numerical viscosity,
\newblock {\em Numer. Math.},
95: 563--598, 2003.

\bibitem{Thomee}
V. Thom\'ee,
\newblock {\em Galerkin Finite Element Methods for Parabolic Problems},
2nd ed.,
\newblock Springer Ser. Comput. Math. 25, Springer, Berlin, 2006.

\bibitem{Trefethen}
L.~N. Trefethen,
\newblock Group velocity in finite difference schemes,
\newblock {\em SIAM Rev.},
24: 113--136, 1982.

\bibitem{Xu94}
J. Xu,
\newblock A novel two-grid method for semilinear elliptic equations,
\newblock {\em SIAM J. Sci. Comput.},
15: 231--237, 1994.

\bibitem{Xu96}
J. Xu,
\newblock Two-grid discretization techniques for linear and nonlinear
{PDE}s,
\newblock {\em SIAM J. Numer. Anal.},
33: 1759--1777, 1996.

\bibitem{Zuazua05}
E. Zuazua,
\newblock Propagation, observation, and control of waves approximated by
finite difference methods,
\newblock {\em SIAM Rev.},
47: 197--243, 2005.

\end{thebibliography}
